\documentclass[reqno,10pt]{amsart}
\usepackage{amssymb,amsmath,amsthm}
\usepackage{a4wide}
\usepackage{mathrsfs}
\usepackage{color}
\usepackage{esint}
\usepackage{stmaryrd}
\usepackage[colorlinks,linkcolor=blue,anchorcolor=blue,citecolor=blue]{hyperref}

\newtheorem{theorem}{Theorem}

\newtheorem{lemma}[theorem]{Lemma}
\newtheorem{proposition}[theorem]{Proposition}
\newtheorem{remark}[theorem]{Remark}

\def\beq{\begin{equation}}
\def\eeq{\end{equation}}
\def\beqs{\begin{equation*}}
\def\eeqs{\end{equation*}}
\def\bal#1\eal{\begin{align}#1\end{align}}
\def\bals#1\eals{\begin{align*}#1\end{align*}}
\def\bsp#1\esp{\begin{split}#1\end{split}}

\def\d{{\mathrm{d}}}

\numberwithin{equation}{section}
\numberwithin{theorem}{section}

\allowdisplaybreaks[3]

\begin{document}
\date{}
\title[One-species non-cutoff VMB system]{
\large{Global dynamics and the diffusive limit of the one-species non-cutoff Vlasov-Maxwell-Boltzmann system}}

\author{Shuangqian Liu, Yuan Xu, Fujun Zhou}

\address[Shuangqian Liu]{School of Mathematics and Statistics, and Key Lab NAA-MOE, Central China Normal University, Wuhan 430079, P.R. China}
\email{sqliu@ccnu.edu.cn}

\address[Yuan Xu]{School of Mathematics and Statistics, and Hubei Key Lab--Math. Sci., Central China Normal University, Wuhan 430079, P.R. China}
\email{yuanxu2019@163.com}

\address[Fujun Zhou]{School of Mathematics, South China University of Technology, Guangzhou 510640, P.R. China}
\email{fujunht@scut.edu.cn}

\begin{abstract}
The rigorous derivation of incompressible plasma equations from kinetic plasma models is a fundamental problem in kinetic theory.
For the one-species Vlasov--Maxwell--Boltzmann system without angular cutoff, global well-posedness and the diffusive limit have remained largely open in the soft-potential regime.
We resolve both problems for $0<s<1$ and $\max\{-3,-\frac32-2s\}<\gamma<-2s$.
For sufficiently small initial perturbations, we construct unique global classical solutions with energy and time-decay estimates uniform in the Knudsen number $\varepsilon\in(0,1]$.
The key analytical contribution is a family of coercivity and trilinear estimates for the non-cutoff collision operators with time-dependent quadratic exponential velocity weights.
Exploiting conservation of energy in collisions, these estimates combine with weighted interpolation to control the high-order electromagnetic terms despite the slow decay of the electric field.
For well-prepared initial data, we prove strong convergence to the incompressible Navier--Stokes--Fourier--Maxwell system as $\varepsilon\to0$.
 \\[3mm]
\\[1mm]
\end{abstract}

\maketitle

\tableofcontents

\section{Introduction}

\subsection{The Problem}
\hspace*{\fill}

The Vlasov--Maxwell--Boltzmann (VMB) system is a fundamental model in plasma physics that describes the evolution of dilute charged particles interacting with their self-consistent electromagnetic field.
The rigorous derivation of fluid equations from this system under diffusive scaling remains an important problem in kinetic theory.

In the physical situation, ions are usually much heavier than electrons, so that electrons move faster than ions.
Hence, ions are often modeled as a fixed ionic background $n_b(x)$, while only electron motion is considered.
In this setting, we consider the following one-species VMB system with a fixed ionic background under the diffusive scaling
\begin{equation}\label{GG1}
\left\{\begin{aligned}
& \varepsilon \partial_t F^\varepsilon+v\cdot\nabla_xF^\varepsilon
   +\left(\varepsilon^2 E^\varepsilon+\varepsilon v\times B^\varepsilon\right)\cdot\nabla_v F^\varepsilon =\frac{1}{\varepsilon}Q(F^\varepsilon,F^\varepsilon), \\
& \partial_tE^\varepsilon-\nabla_x\times B^\varepsilon=
    -\frac{1}{\varepsilon}\int_{\mathbb{R}^3}v F^\varepsilon\d v, \\
& \partial_tB^\varepsilon+\nabla_x\times E^\varepsilon=0, \\
& \nabla_x\cdot E^\varepsilon=\int_{\mathbb{R}^3}  F^\varepsilon\d v-n_b, \quad \nabla_x\cdot B^\varepsilon=0.
\end{aligned}\right.
\end{equation}
Here, $F^{\varepsilon}(t,x,v)\geq 0$ represents the density distribution function of the electrons at time $t\geq0$, position $x\in \mathbb{R}^{3}$ and velocity $v\in \mathbb{R}^{3}$, and the uniform ionic background $n_b>0$ is normalized to $1$.
The electromagnetic field $[E^\varepsilon(t, x), B^\varepsilon(t,x)]$ is coupled with the density distribution function $F^{\varepsilon}(t,x,v)$ through the celebrated Maxwell equations in \eqref{GG1}.
The positive parameter $\varepsilon\in (0,1]$ is the Knudsen number, which is the ratio of the mean free path to the macroscopic length scale.
The initial data of the VMB system (\ref{GG1}) are imposed as
\begin{align}\label{GG1-initial}
  F^\varepsilon(0, x,v)= F_{0}^\varepsilon(x,v), \quad E^\varepsilon(0,x)=E^\varepsilon_0(x), \quad B^\varepsilon(0,x)=B^\varepsilon_0(x)
\end{align}
with the compatibility conditions
\begin{align*}
\nabla_x \cdot E_0^\varepsilon= \int_{\mathbb{R}^3} F^\varepsilon_{0}\d v-1,
    \quad \nabla_x \cdot B_0^\varepsilon=0.
\end{align*}
The Boltzmann collision operator $Q(F, G)$ is a bilinear operator, defined as
\begin{align*}
Q(F, G)(v) :=
\iint_{\mathbb{R}^{3} \times \mathbb{S}^2} B(v-v_*,\sigma)[F(v'_{*})G(v')-F(v_*)G(v)]\d v_* \d \sigma.
\end{align*}
Here, $v$, $v_*$ and $v^{\prime}$, $v_*^{\prime}$ are the velocities of a pair of particles before and after collision, respectively. They are connected through the formulas
$$
v^{\prime}=\frac{v+v_*}{2}+\frac{|v-v_*|}{2} \sigma, \quad v_*^{\prime}=\frac{v+v_*}{2}-\frac{|v-v_*|}{2} \sigma, \quad \sigma \in \mathbb{S}^2.
$$

The Boltzmann collision kernel $B(v-v_*, \sigma)$ is a non-negative function, which only depends on the relative velocity $|v-v_*|$ and on the deviation angle $\theta$ given by $\cos \theta=\langle\sigma, k\rangle$, where $k=(v-v_*) /|v-v_*|$ and $\langle\cdot, \cdot\rangle$ is the usual dot product in $\mathbb{R}^3$. As in \cite{GS2011}, without loss of generality, we may suppose that $B(v-v_*, \sigma)$ is supported on $\cos \theta \geq 0$, i.e. $0 \leq \theta \leq \pi/2$.
Throughout this paper, the collision kernel $B(v-v_*, \sigma) $ is supposed to satisfy the following assumptions:
\begin{itemize}
\item
The collision kernel $B(v-v_*,\sigma)$ has the product form
$$
B(v-v_*, \sigma)=\Phi(|v-v_*|) b(\cos \theta).
$$
Here, $\Phi$ and $b$ are non-negative functions.

\item
The angular function $\sigma \rightarrow b(\langle\sigma,(v-v_*) /|v-v_*|\rangle)$ is not locally integrable on $\mathbb{S}^2$, i.e.
$$
\int_{\mathbb{S}^2} b(\cos \theta) \d \sigma=2 \pi \int_0^{\pi / 2} \sin \theta b(\cos \theta) \d \theta=\infty,
$$
and there exists a constant $c_b>0$ such that
$$
\frac{c_b}{\theta^{1+2 s}} \leq \sin \theta b(\cos \theta) \leq \frac{1}{c_b \theta^{1+2 s}}, \quad s \in (0,1).
$$

\item
The kinetic function $z \rightarrow \Phi(|z|)$ satisfies
$$
\Phi(|z|)=C_{\Phi}|z|^\gamma,
$$
where $C_\Phi>0$ is a constant and the exponent $\gamma>-3$ is determined by the intermolecular interaction mechanism.
\end{itemize}

We refer to the cases $-3<\gamma<-2s$ and $\gamma\geq-2s$ as ``soft potentials'' and ``hard potentials'', respectively. The current work is restricted to the case of
\begin{align*}
  \max\left\{-3,-\frac{3}{2}-2s\right\} <\gamma <-2s \; \text{ and }\; 0<s<1.
\end{align*}
Moreover, the ranges $0<s<\frac12$ and $\frac12\leq s<1$ correspond to weak and strong angular singularities, respectively \cite{DL2013}.
When the intermolecular interaction potential takes the inverse power law in the form of $U(|x|)=|x|^{-(\ell-1)}$ with $2 < \ell < \infty$, the collision kernel $B(v-v_*, \sigma)$ in three space dimensions satisfies the above assumptions with $\gamma=\frac{\ell-5}{\ell-1}$ and $s=\frac{1}{\ell-1}$. Notice that $\gamma \to -3$ and $s \to 1$ as $\ell \to 2$ in the limiting case, for which grazing collisions between particles dominate and the Boltzmann collision term needs to be replaced by the classical Landau collision term for the Coulomb potential, cf. \cite{V}.

For the one-species VMB system with angular cutoff and a uniform ionic background, global classical solutions near the global Maxwellian were constructed by Duan \cite{Duan2011} for the hard-sphere model $\gamma=1$ and by Lei--Zhao \cite{LZ2016} for soft potentials $-1<\gamma<0$.
Li--Yang--Zhong \cite{LYZ2016-1} subsequently investigated the spectral structure and optimal time-decay rates for the hard-sphere case.
Under the diffusive scaling considered here, Wu--Zhou--Gong--Xu \cite{WZGX2025SIMA} rigorously justified the incompressible Navier--Stokes--Fourier--Maxwell (NSFM) limit for hard potentials $0\leq\gamma\leq1$ via uniform-in-$\varepsilon$ weighted energy estimates, whereas Chen--Yang--Zhong \cite{CYZ2026JFA} derived convergence rates and precise initial-layer estimates for the hard-sphere model under an alternative electromagnetic scaling.
Regarding the two-species counterpart with angular cutoff, global existence and time-decay estimates near the global Maxwellian were established in \cite{Guo2003,Strain2006CMP,DS2011,DLYZ2017}.
Moreover, the corresponding diffusive limit toward the incompressible NSFM system with Ohm's law was proved via the energy method in \cite{JL2019} and via Hilbert expansions in \cite{JLZ2023ARMA}, with convergence rates further obtained through spectral analysis in \cite{YZ24}.

Substantial progress has also been made on the two-species VMB system without angular cutoff.
Building on the theory for the non-cutoff Boltzmann equation \cite{GS2011,AMUXY2011AA,AMUXY2012JFA}, Duan--Liu--Yang--Zhao \cite{DLYZ2013} obtained global classical solutions for soft potentials $\max\{-3,-\frac32-2s\}<\gamma<-2s$ with strong angular singularity $\frac12\leq s<1$, while Fan--Lei--Liu--Zhao \cite{FLLZ2018} extended the result to weak angular singularity $0< s<\frac12$.
Concerning the diffusive limit toward the two-fluid incompressible NSFM system with Ohm's law, Jiang--Lei \cite{JL2023ARXIV} addressed soft potentials with $\frac12\leq s<1$, whereas Xu--Zhou--Gong--Wu \cite{XZGW2026SIMA} resolved the full angular singularity range $0<s<1$ for soft potentials as well as hard potentials.

Beyond the incompressible NSFM regime, several other hydrodynamic limits have been explored.
Specifically, Jang \cite{J2009} derived the Vlasov--Navier--Stokes--Fourier limit, while in the hyperbolic scaling, the compressible Euler--Maxwell limit for the one-species VMB system was justified by Duan--Yang--Yu \cite{DYY2023M3AS} for the hard-sphere model and by Lei--Liu--Xiao--Zhao \cite{LLXZ2023ARXIV} for non-cutoff soft potentials.
Within the framework of renormalized solutions, Ars\'{e}nio--Saint-Raymond \cite{AS2019} comprehensively investigated various limits toward incompressible viscous electro-magneto-hydrodynamics.
For further background on kinetic equations and their hydrodynamic limits, we refer the reader to the monographs \cite{CIP,Saint} and the related works \cite{Ca,DEL,BGL91,BU91,BGL93,Lions1994-1,Lions1994-2,LM,GS,guo2006,guo2010cmp,GJJ,JXZ2018}, as well as the references therein.

Despite these advances, to the best of our knowledge, global well-posedness and the incompressible NSFM limit for the one-species non-cutoff VMB system remain open in the perturbative framework.
The present paper resolves both problems for soft potentials $\max\{-3,-\frac32-2s\}<\gamma<-2s$ in the full angular singularity range $0< s<1$.
Notably, the weighted energy method developed for the two-species non-cutoff system in \cite{XZGW2026SIMA} does not directly extend to the one-species case.
In that work, time-dependent exponential weights with exponents linear in $\langle v\rangle$ are used to generate additional velocity dissipation.
Thanks to the faster electric-field decay in the two-species model, this dissipation is strong enough to control the weighted microscopic terms in the high-order electromagnetic estimates.
In the one-species case, however, the slower decay of the electric field does not provide the uniform smallness needed for this absorption argument.
To compensate for this slower decay, it is natural to use quadratic exponential velocity weights to produce stronger dissipation at large velocities.
However, as pointed out by Lei--Liu--Xiao--Zhao \cite{LLXZ2023ARXIV}, establishing suitable non-cutoff collision estimates with such weights is a major technical difficulty.

Our primary technical contribution in this paper is to establish new weighted coercivity and trilinear estimates for the non-cutoff collision operators $L$ and $\Gamma$ with a time-dependent quadratic exponential velocity weight.
This weight produces stronger dissipation at large velocities, and the new collision estimates allow us to apply it in the non-cutoff setting.
Coupled with weighted interpolation for the electromagnetic terms and a carefully designed hierarchy of polynomial velocity weights, our approach yields global energy and time-decay estimates uniform in $\varepsilon\in(0,1]$.
These uniform bounds not only imply global well-posedness, but also enable us to rigorously justify the diffusive limit toward the incompressible NSFM system when combined with estimates for the acoustic and microscopic components.
\medskip

\subsection{Reformulation}
\hspace*{\fill}

We investigate the diffusive limit of the one-species VMB system (\ref{GG1}) in the perturbative framework around the global Maxwellian
$$
\mu=\mu(v) :=(2 \pi)^{-3 / 2} e^{-|v|^{2} / 2}.
$$
To this end, we write
$$
F^{\varepsilon}(t,x,v)=\mu+\varepsilon \mu^{1/2}f^{\varepsilon}(t,x,v),\quad
F^\varepsilon_{0}(x,v)=\mu+\varepsilon \mu^{1/2}f_{0}^{\varepsilon}(x,v),
$$
where $f^\varepsilon$ and $f_0^\varepsilon$ denote the fluctuations.
The VMB system \eqref{GG1} then becomes the following equivalent perturbed system
\begin{equation}\label{rVMB}
\left\{\begin{aligned}
& \partial_{t} f^{\varepsilon}+\frac{1}{\varepsilon}v \cdot \nabla_{x} {f}^{\varepsilon}
+(\varepsilon E^\varepsilon+v\times B^\varepsilon)\cdot \nabla_{v} {f}^{\varepsilon}
- E^{\varepsilon} \cdot v \mu^{1/2}+\frac{1}{\varepsilon^{2}} L {f}^{\varepsilon} \\
& \qquad\qquad\qquad\qquad\qquad\qquad= \frac{\varepsilon}{2}E^{\varepsilon} \cdot v {f}^{\varepsilon}+\frac{1}{\varepsilon} \Gamma({f}^{\varepsilon}, {f}^{\varepsilon}),  \\
& \partial_tE^\varepsilon-\nabla_x\times B^\varepsilon=-\int_{\mathbb{R}^3}v \mu^{1/2}f^\varepsilon\d v, \\
& \partial_tB^\varepsilon+\nabla_x\times E^\varepsilon=0, \\
& \nabla_x\cdot E^\varepsilon=\varepsilon \int_{\mathbb{R}^3}  \mu^{1/2}f^{\varepsilon}\d v,\quad \nabla_x\cdot B^\varepsilon=0,
\end{aligned}\right.
\end{equation}
with initial data
\begin{align}\label{rVMBinitial}
f^\varepsilon(0,x,v)=f^\varepsilon_{0}(x,v),\qquad
E^\varepsilon(0,x)=E^\varepsilon_{0}(x), \qquad  B^\varepsilon(0,x)=B^\varepsilon_{0}(x),
\end{align}
where  the linearized Boltzmann collision operator $Lf^{\varepsilon}$ and the nonlinear collision term $\Gamma(f^{\varepsilon}, f^{\varepsilon})$ are respectively defined by
\begin{align*}
Lf^{\varepsilon}\;&:=-\mu^{-1/2}\left[Q(\mu,\mu^{1/2}f^\varepsilon)+Q(\mu^{1/2}f^\varepsilon,\mu)\right],\\
\Gamma(f^{\varepsilon}, f^{\varepsilon})\;&:=\mu^{-1/2} Q(\mu^{1/2}f^\varepsilon,\mu^{1/2}f^\varepsilon).
\end{align*}

It is well known that the operator $L$ is non-negative with null space
\begin{align*}
\mathcal{N}(L) = \text{span} \Big\{ \mu^{1/2},   v \mu^{1/2} ,   (|v|^{2}-3) \mu^{1/2} \Big \},
\end{align*}
cf. \cite{CIP}.
For a given function $f(t,x,v)$, the following macro-micro decomposition was introduced in \cite{Guo2003}
\begin{equation}\label{f decomposition}
f=\mathbf{P}f+\mathbf{P}^{\perp}f.
\end{equation}
Here, $\mathbf{P}$ denotes the orthogonal projection from $L^2(\mathbb{R}_{v}^3)$ to $\mathcal{N}(L)$, defined by
\begin{equation}\label{Pf define}
\mathbf{P}f:=\big\{a(t,x)+ v\cdot b(t,x) +(|v|^2-3)c(t,x)\big\}\mu^{1/2},
\end{equation}
where
\begin{align*}
a(t,x)=\langle\mu^{1/2},f\rangle,\qquad
b(t,x)=\langle v\mu^{1/2},f\rangle,\qquad
c(t,x)=\frac{1}{6}\langle (|v|^2-3)\mu^{1/2},f\rangle.
\end{align*}
The complementary projection is defined by $\mathbf{P}^{\perp}:=\mathbf{I}-\mathbf{P}$, where $\mathbf{I}$ denotes the identity operator.
We call $\mathbf{P}f$ and $\mathbf{P}^{\perp}f$ the macroscopic and microscopic components of $f(t,x,v)$, respectively.

\smallskip

\subsection{Notations}
\hspace*{\fill}

Throughout the paper, $C$ denotes a generic positive constant independent of $\varepsilon$.
We use $X \lesssim Y$ to denote $X \leq CY$, where $C$ is a constant independent of $X$ and $Y$.
We also use the notation $X \approx Y$ to represent $X\lesssim Y$ and $Y\lesssim X$.
The notation $X \ll 1$ means that $X$ is a sufficiently small positive constant.

The multi-indexes $\alpha = [\alpha_1, \alpha_2, \alpha_3]\in\mathbb{N}^3$ and $\beta = [\beta_1, \beta_2, \beta_3]\in\mathbb{N}^3$ will be used to record space and velocity derivatives, respectively.
The mixed derivatives are denoted by $\partial^\alpha_\beta=\partial^\alpha_x\partial^\beta_v=\partial_{x_1}^{\alpha_1}\partial_{x_2}^{\alpha_2}\partial_{x_3}^{\alpha_3}\partial_{v_1}^{\beta_1}\partial_{v_2}^{\beta_2}\partial_{v_3}^{\beta_3}$.
For pure spatial and velocity derivatives, we write $\partial^\alpha=\partial_x^\alpha=\partial_{x_1}^{\alpha_1}\partial_{x_2}^{\alpha_2}\partial_{x_3}^{\alpha_3}$ and $\partial_\beta=\partial_v^\beta=\partial_{v_1}^{\beta_1}\partial_{v_2}^{\beta_2}\partial_{v_3}^{\beta_3}$, respectively.
The length of $\alpha$ is denoted by $|\alpha|=\alpha_1+\alpha_2+\alpha_3$.
If each component of $\theta$ is no greater than the corresponding component of $\bar{\theta}$, we write $\theta\leq\bar{\theta}$.
We write $\theta<\bar{\theta}$ when $\theta\leq\bar{\theta}$ and $|\theta|<|\bar{\theta}|$.

We use $\|\cdot\|_{L^p}$ to denote the $L_x^p$ norm for functions depending only on $x$, and the $L_{x,v}^p$ norm for functions depending on both $x$ and $v$.
When $p=2$, we simply write $\|\cdot\|$.
More generally, for a function space $X$ in the velocity variable $v$, we use $|\cdot|_X$ to denote its norm and set
\begin{align*}
\|f\|_{L_x^pX}:=\big\||f(x,\cdot)|_X\big\|_{L_x^p}.
\end{align*}
In particular, for the weighted and dissipation norms defined below, we write $\|f\|_{L_\ell^2}:=\|f\|_{L_x^2L_\ell^2}$ and $\|f\|_D:=\|f\|_{L_x^2L_D^2}$.
The brackets $\langle\cdot,\cdot\rangle$ denote the $L_v^2$ inner product, while $(\cdot,\cdot)$ denotes the $L_x^2$ or $L_{x,v}^2$ inner product, as appropriate.
For any integer $N\geq0$, we define
\begin{align*}
\|f\|_{H_x^NL_v^2}^2:=\sum_{|\alpha|\leq N}\|\partial^\alpha f\|^2, \qquad
\|f\|_{H_x^NL_D^2}^2:=\sum_{|\alpha|\leq N}\|\partial^\alpha f\|_D^2.
\end{align*}
Other integer-order spatial Sobolev norms are defined analogously.

The operator $\Lambda^{-\varrho}$ is defined by
\begin{align*}
\Lambda^{-\varrho} f(x):= {(2\pi)^{-3/2}} \int_{\mathbb{R}^3} |y|^{-\varrho} \widehat{f}(y) e^{ix \cdot y} \d y,
\end{align*}
where $\widehat{f}(y) :=  {(2\pi)^{-3/2}} \int_{\mathbb{R}^3} f(x) e^{-ix \cdot y}\d x$ represents the Fourier transform of $f(x)$.

As described in \cite{AMUXY2012JFA, DLYZ2013}, we introduce the following norm for the non-cutoff VMB system
\begin{align}
\begin{split}\label{D-norm}
|f|_{L^2_D}^2=|f|_{D}^2:= & \iiint_{\mathbb{R}^6 \times \mathbb{S}^2} B\left(v-v_*, \sigma\right) \mu(v_*)\Big(f(v^\prime)-f(v)\Big)^2 \d v_* \d v \d \sigma \\
& +\iiint_{\mathbb{R}^6 \times \mathbb{S}^2} B\left(v-v_*, \sigma\right) f(v_*)^2\left(\mu(v^{\prime})^{1/2}- \mu(v)^{1/2}\right)^2 \d v_* \d v \d \sigma.
\end{split}
\end{align}
Let $\langle v\rangle:=(1+|v|^2)^{\frac{1}{2}}$. For $\ell \in \mathbb{R}$, $L^2_\ell$ denotes the weighted space with the norm
\begin{align*}
|f|_{L_{\ell}^2}^2:=\int_{\mathbb{R}^3}\langle v\rangle^{2\ell}|f(v)|^2 \d v.
\end{align*}
And the weighted fractional Sobolev norm $\left|f(v)\right|_{H^s_\ell}^2=\left|\langle v \rangle^\ell f(v)\right|_{H^s}^2$ is given by
\begin{align*}
|f|_{H_{\ell}^s}^2:=|f|_{L_{\ell}^2}^2+\iint_{\mathbb{R}^3 \times \mathbb{R}^3} \frac{\big[\langle v\rangle^{\ell}  f(v)-\langle v^{\prime}\rangle^{\ell}  f(v^{\prime})\big]^2}{|v-v^{\prime}|^{3+2 s}} \chi_{|v-v^{\prime}| \leq 1} \d v \d v^{\prime},
\end{align*}
which is equivalent to
\begin{align*}
| f|_{H_{\ell}^s}^2:=\int_{\mathbb{R}^3}\langle v\rangle^{2\ell} \big|(1-\Delta_v)^{\frac{s}{2}} f(v)\big|^2 \d v.
\end{align*}
Notice that
\begin{align}\label{norm inequality}
| f |_{L_{\gamma/2 + s}^2}^2+| f |_{H_{\gamma/2}^s}^2 \lesssim | f |_{D}^2
\lesssim | f |_{H_{\gamma/2+  s}^s}^2,
\end{align}
cf. \cite{AMUXY2012JFA}.

\subsection{Main Results}
\hspace*{\fill}

To state the main results of this paper, we introduce the following fundamental energy functional and the corresponding dissipation rate functional
\begin{align}
\label{without weight energy functional}
{\mathcal{E}}_{N} (t) \sim\;
&\sum_{|\alpha|\leq N} \left\|\partial^{\alpha}f^{\varepsilon}\right\|^2
+\sum_{|\alpha|\leq N} \left\|\partial^{\alpha}E^{\varepsilon}\right\|^2
+\sum_{|\alpha|\leq N} \left\|\partial^{\alpha}B^{\varepsilon}\right\|^2,  \\
\label{without weight dissipation functional}
{\mathcal{D}}_{N} (t)\sim \;
&\sum_{1 \leq |\alpha| \leq N} \left \| \partial^{\alpha} \mathbf{P}f^\varepsilon \right \|^2
+ \frac{1}{\varepsilon^2} \sum_{|\alpha| \leq N}\left\|\partial^{\alpha}
\mathbf{P}^{\perp}f^\varepsilon\right\|^2_D + \varepsilon^2 \|a^\varepsilon\|^2
 \nonumber\\
\;&+ \sum_{1\leq |\alpha| \leq N-1} \left\|\partial^{\alpha}E^{\varepsilon}\right\|^2
+ \sum_{2 \leq |\alpha| \leq N-1} \left\|\partial^{\alpha}B^{\varepsilon}\right\|^2,
\end{align}
respectively, where the integer $N\in\mathbb{Z}$ will be determined later.

Because the linearized Boltzmann operator $L$ provides weaker dissipation for soft potentials than for hard potentials, we introduce the following time-velocity weight:
\begin{align}\label{weight function}
w_{l}(\alpha,\beta):= \overline{w}_{l}(\alpha,\beta)
e^\frac{q \langle v \rangle^2}{(1+t)^\vartheta}.
\end{align}
This weight is used to handle the force term generated by the self-consistent electromagnetic field.
Here, $0<q\ll1$, and the constants $l$ and $\vartheta$ will be determined later.
The polynomial factor $\overline w_l(\alpha,\beta)$ is defined by
\begin{align}\label{weight function2}
\overline{w}_{l}(\alpha,\beta):=
\langle v \rangle^{l-(m+\gamma+2s)|\alpha|-m|\beta|}
\quad\; \text{with} \;\;\;
m=1+\frac{1-2(\gamma+2s)}{\min\{2s,1\}}.
\end{align}
To overcome the regularity loss of the electromagnetic field, we use the following notation
\begin{align}\label{weight function1}
\widetilde{w}_{l}(\alpha,\beta):=
\begin{cases}
w_{l}(\alpha,\beta), &\quad \text{ if } |\alpha|+|\beta| \leq N-1, \\[1mm]
(1+t)^{-\frac{1+\varepsilon_0}{4}}w_{l}(\alpha,\beta), &\quad\text{ if } |\alpha|+|\beta| =N,
\end{cases}
\end{align}
where $\varepsilon_0> 0$ is small enough.
The notation $\widetilde w_l(\alpha,\beta)$ is introduced only for convenience.
For $|\alpha|+|\beta|=N$, we can first derive the energy estimate using $w_l(\alpha,\beta)$ and then multiply the resulting estimate by $(1+t)^{-(1+\varepsilon_0)/2}$.
Moreover, sometimes we write $\widetilde{w}_{l}(|\alpha|,|\beta|)$ to denote $\widetilde{w}_{l}(\alpha, \beta)$ for simplicity.
Then the weighted energy functional $\widetilde{\mathcal{E}}_{N,l} (t)$ and the corresponding dissipation rate functional $\widetilde{\mathcal{D}} _{N,l} (t)$ are defined as
\begin{align}\label{energy functional}
\widetilde{\mathcal{E}}_{N,l} (t) \sim\;
& {\mathcal{E}}_{N} (t)  + \sum_{\substack{|\alpha|+|\beta| \leq N \\ |\alpha| \leq N-1}}\left\|\widetilde{w}_l (\alpha, \beta) \partial_{\beta}^{\alpha} \mathbf{P}^{\perp}f^{\varepsilon}\right\|^2
+\varepsilon\sum_{|\alpha|=N}\left\|\widetilde{w}_l (\alpha, 0) \partial^{\alpha} f^{\varepsilon}\right\|^2 , \\
\label{dissipation functional}
\widetilde{\mathcal{D}}_{N,l} (t)\sim \;& {\mathcal{D}}_{N} (t)
+\frac{1}{\varepsilon^{2}} \sum_{\substack{|\alpha|+|\beta| \leq N \\
|\alpha| \leq N-1}}\left\|\widetilde{w}_l (\alpha, \beta) \partial_{\beta}^{\alpha}
\mathbf{P}^{\perp} f^{\varepsilon}\right\|^2_D
+\frac{1}{\varepsilon}\sum_{|\alpha|=N} \left\|\widetilde{w}_l (\alpha, 0) \partial^{\alpha}
\mathbf{P}^{\perp} f^{\varepsilon}\right\|^2_D \nonumber\\
&+ \frac{1}{(1+t)^{1+\vartheta}} \left\{ \sum_{\substack{|\alpha|+|\beta| \leq N \\
|\alpha| \leq N-1}}\left\|\langle v \rangle \widetilde{w}_l (\alpha, \beta) \partial_{\beta}^{\alpha}
\mathbf{P}^{\perp} f^{\varepsilon}\right\|^2 +  \varepsilon \sum_{|\alpha|= N} \left\|\langle v \rangle \widetilde{w}_l (\alpha, 0) \partial^{\alpha} f^{\varepsilon}\right\|^2 \right\} \nonumber \\
&+ \frac{1}{1+t} \left\{ \sum_{\substack{|\alpha|+|\beta| = N \\
|\alpha| \leq N-1}}\left\| \widetilde{w}_l (\alpha, \beta) \partial_{\beta}^{\alpha}
\mathbf{P}^{\perp} f^{\varepsilon}\right\|^2 +  \varepsilon \sum_{|\alpha|= N} \left\| \widetilde{w}_l (\alpha, 0) \partial^{\alpha} f^{\varepsilon}\right\|^2 \right\},
\end{align}
respectively.

To obtain the time-decay rate needed to close the energy estimate, we define the following energy functional $\mathcal E^k_{N_0}(t)$ and the corresponding dissipation rate functional $\mathcal D^k_{N_0}(t)$
\begin{align}
\label{low k energy}
\mathcal{E}^k_{N_0}(t)\sim\;
& \sum_{k \leq |\alpha| \leq N_0} \left\|\partial^{\alpha}f^{\varepsilon}\right\|^2
+\sum_{k \leq |\alpha| \leq N_0} \left\|\partial^{\alpha}E^{\varepsilon}\right\|^2
+\sum_{k \leq |\alpha| \leq N_0} \left\|\partial^{\alpha}B^{\varepsilon}\right\|^2,  \\
\label{low k dissipation}
{\mathcal{D}}^k_{N_0} (t)\sim \;
& \sum_{k+1 \leq |\alpha| \leq N_0} \left \| \partial^{\alpha} \mathbf{P}f^\varepsilon \right \|^2
+ \frac{1}{\varepsilon^2} \sum_{k \leq |\alpha| \leq N_0} \left\|\partial^{\alpha}
\mathbf{P}^{\perp}f^\varepsilon\right\|^2_D
+ \varepsilon^2 \sum_{ |\alpha| = k } \left \| \partial^{\alpha} a^\varepsilon \right \|^2 \nonumber\\
\;& +\sum_{k+1 \leq |\alpha| \leq N_0-1} \left\|\partial^{\alpha} E^{\varepsilon}\right\|^2
+  \sum_{k+2 \leq |\alpha| \leq N_0-1} \left\|\partial^{\alpha} B^{\varepsilon}\right\|^2,
\end{align}
respectively, where the integer $N_0\in\mathbb{Z}$ will be determined later.

Moreover, the singular dependence on $\varepsilon$ makes it difficult to derive uniform nonlinear time-decay estimates directly from the Duhamel formula and linear semigroup bounds.
We therefore resort to the Sobolev spaces with negative index.
For this, we define the following energy functional $\overline{\mathcal{E}}_{N,l}(t)$ and the corresponding dissipation rate functional $\overline{\mathcal{D}}_{N,l}(t)$
\begin{align}
\label{negative sobolev energy}
\overline{\mathcal{E}}_{N,l} (t) \sim\;
& \widetilde{\mathcal{E}}_{N,l} (t) +\left\| \Lambda^{-\varrho}f^\varepsilon\right\|^2
+\left\| \Lambda^{-\varrho}E^\varepsilon\right\|^2
+\left\| \Lambda^{-\varrho}B^\varepsilon\right\|^2,\\
\label{negative sobolev dissipation}
\overline{\mathcal{D}}_{N,l} (t) \sim\;&\widetilde{\mathcal{D}}_{N,l}(t)
+\left\| \Lambda^{1-\varrho} \mathbf{P}f^\varepsilon\right\|^2
+\frac{1}{\varepsilon^2}\left\| \Lambda^{-\varrho} \mathbf{P}^{\perp}f^\varepsilon\right\|^2_D \nonumber\\
\;& + \varepsilon^2 \left\| \Lambda^{-\varrho} a^\varepsilon \right\|^2
+  \left\| \Lambda^{1-\varrho}E^\varepsilon\right\|^2
+  \left\| \Lambda^{2-\varrho}B^\varepsilon\right\|^2,
\end{align}
respectively, where the positive constant $\varrho>0$ will be determined later.
\medskip

Our first main result provides energy and time-decay estimates uniform in $\varepsilon\in(0,1]$ for global solutions to the VMB system \eqref{rVMB} with $\max\{-3,-\frac32-2s\}<\gamma<-2s$ and $0<s<1$.

\begin{theorem}\label{mainth1}
Let $0 < s <1$, $\max\{-3, -\frac{3}{2}-2s\}<\gamma<-2s$ and $0<\varepsilon \leq 1$. Let $m$ be as defined in \eqref{weight function2}, and introduce the following constants in sequence:
\begin{align}\label{mainth1 assumption}
\begin{split}
&N_0\geq5\; ( N_0\in\mathbb Z),\quad N=N_0+2, \quad l\geq mN,
\quad \max\left\{1,\frac{2-5(\gamma+2s)}{4-2(\gamma+2s)}\right\}<\varrho<\frac32,\\
&0<\vartheta\leq \frac{(3+2\varrho)(2-\gamma-2s)}{8(1-\gamma-2s)}-1,
\quad 0<q\ll\frac{1}{4},
\quad 0<\varepsilon_0\ll \varrho+\frac{1}{2}-2\vartheta.
\end{split}
\end{align}
There exists a sufficiently small constant $\delta_1>0$, independent of $\varepsilon$, such that, if the initial data satisfy $F_0^\varepsilon=\mu+\varepsilon\mu^{1/2}f_0^\varepsilon\geq0$ and
\begin{align*}
\overline{\mathcal E}_{N,l}(0)\leq\delta_1^2,
\end{align*}
then the VMB system \eqref{rVMB}--\eqref{rVMBinitial} admits a unique global solution $(f^\varepsilon,E^\varepsilon,B^\varepsilon)$, which satisfies $F^\varepsilon=\mu+\varepsilon\mu^{1/2}f^\varepsilon\geq0$ and
\begin{align}
&\overline{\mathcal E}_{N,l}(t)+ \int_0^t\overline{\mathcal D}_{N,l}(\tau)\d\tau\leq C\overline{\mathcal E}_{N,l}(0),
\label{thm1 estimate 1}\\
&\mathcal E_{N_0}^k(t)\leq C\overline{\mathcal E}_{N,l}(0)(1+t)^{-\frac{k+\varrho}{2}},
\qquad k=0,1,2, \label{thm1 estimate 2}
\end{align}
for all $t\geq0$, where $C>0$ is independent of $t$ and $\varepsilon$.
\end{theorem}
\medskip

The second main result concerns the hydrodynamic limit from the one-species VMB system \eqref{rVMB} to the incompressible NSFM system.

\begin{theorem}\label{mainth2}
Let $(f^\varepsilon,E^\varepsilon,B^\varepsilon)$ be the global solution to the VMB system \eqref{rVMB}--\eqref{rVMBinitial} constructed in Theorem \ref{mainth1}.
Suppose that there exist scalar functions $\rho_0,\theta_0\in H^N_x$ and vector-valued functions $u_0,E_0,B_0\in H^N_x$ satisfying
\begin{align*}
\nabla_x\cdot u_0=0,\qquad
\rho_0+\theta_0=0,\qquad
\nabla_x\cdot E_0=\nabla_x\cdot B_0=0,
\end{align*}
such that
\begin{align*}
{f}_0^\varepsilon \to {f}_0\; \text{ strongly in } H^{N}_x L^2_v, \;\;\;\;\;
(E_0^\varepsilon, B_0^\varepsilon) \to (E_0, B_0)\;   \text{ strongly in } H^{N}_x
\end{align*}
as $\varepsilon \to 0$, where $f_0=f_0(x,v)$ has the form
\begin{align*}
f_0=\Big\{\rho_0+u_0\cdot v+\theta_0\frac{|v|^2-3}{2}\Big\}\mu^{1/2}.
\end{align*}
Then the following convergences hold
\begin{equation*}
\begin{aligned}
{f}^\varepsilon &\to {f}
\;&& \text{ weakly}\!-\!* \text{ in }  L^\infty(\mathbb{R}^+; H^N_x L^2_v),
&&\text{ strongly in } C([0,\infty); H_x^{N-1} L^2_v), \\
(E^\varepsilon, B^\varepsilon)  &\to (E, B)
\;&& \text{ weakly}\!-\!*  \text{ in }   L^\infty(\mathbb{R}^+;H^{N}_{x}),
&&\text{ strongly in } C([0,\infty); H_x^{N-1}), \\
\end{aligned}
\end{equation*}
as $\varepsilon \to 0$, where $f={f}(t, x, v)$ has the form
\begin{align*}
f=\Big\{\rho+u\cdot v+\theta\frac{|v|^2-3}{2}\Big\}\mu^{1/2}.
\end{align*}
Moreover, we have
\begin{equation*}
\begin{aligned}
\left(\mathcal P\big\langle f^\varepsilon,v\mu^{1/2}\big\rangle,
\frac15\big\langle f^\varepsilon,(|v|^2-5)\mu^{1/2}\big\rangle\right)
&\to(u,\theta)&\quad&\text{strongly in }C([0,\infty);H_x^{N-1}),\\
\left(\frac13\big\langle f^\varepsilon,|v|^2\mu^{1/2}\big\rangle,
\mathcal P^\perp\big\langle f^\varepsilon,v\mu^{1/2}\big\rangle\right)
&\to(0,0) &\quad&\text{strongly in }C([0,\infty);H_x^{N-1})
\end{aligned}
\end{equation*}
as $\varepsilon \to 0$. The limiting functions satisfy
\begin{align*}
(\rho,u,\theta,E,B)\in C([0,\infty);H^{N-1}_x) \cap L^\infty(\mathbb R^+;H^N_x)
\end{align*}
and form the unique global solution to the incompressible NSFM system
\begin{equation}\label{INSFM limit}
\left\{
\begin{aligned}
&\partial_t u+u\cdot\nabla_xu-\nu\Delta_xu+\nabla_xp=E+u\times B, &&\nabla_x\cdot u=0,\\
&\partial_t\theta+u\cdot\nabla_x\theta-\kappa\Delta_x\theta=0, &&\rho+\theta=0,\\
&\partial_tE-\nabla_x\times B=-u, &&\nabla_x\cdot E=0,\\
&\partial_tB+\nabla_x\times E=0, &&\nabla_x\cdot B=0,
\end{aligned}
\right.
\end{equation}
with initial data
\begin{align}\label{INSFM limit initial}
u(0,x)=u_0(x),\qquad
\theta(0,x)=\theta_0(x),\qquad
E(0,x)=E_0(x),\qquad
B(0,x)=B_0(x).
\end{align}
Here, $\mathcal{P}$ is the Leray projection, while the viscosity coefficient $\nu$ and the heat conductivity coefficient $\kappa$ are defined by \eqref{coefficient-define}.
\end{theorem}
\medskip

\begin{remark}\label{remark-0}
The convergence in Theorem \ref{mainth2}, combined with the uniform time-decay estimates in Theorem \ref{mainth1}, implies that
\begin{align*}
&\sup_{t\geq0}\left\|F^\varepsilon(t,x,v)-\mu(v)-\varepsilon\mu^{1/2}(v)f(t,x,v)\right\|_{H^{N-1}_x L_v^2}=o(\varepsilon),\\
&\sup_{t\geq0}\Big\{\|E^\varepsilon(t,x)-E(t,x)\|_{H^{N-1}_x}+\|B^\varepsilon(t,x)-B(t,x)\|_{H^{N-1}_x}\Big\}=o(1)
\end{align*}
as $\varepsilon\to0$.
Thus, the incompressible NSFM system \eqref{INSFM limit} provides a first-order approximation of the original VMB system \eqref{GG1}.
\end{remark}
\medskip

\begin{remark}\label{remark-1.4}
To the best of our knowledge, this is the first work to establish global well-posedness together with the incompressible NSFM limit for the one-species VMB system with non-cutoff potentials in the perturbative framework.

More importantly, the main technical contribution of this work lies in the new weighted estimates for $L$ and $\Gamma$ established in Section \ref{new-weighted-estimate}.
These estimates may also help relax restrictions on collision parameters and velocity weights in related problems involving non-cutoff collision operators.
They also cover hard potentials $\gamma+2s\geq0$, providing a basis for extending the present approach to this case through suitable modifications of the polynomial velocity weights.
We restrict our analysis to soft potentials, where the weaker collision dissipation creates the main difficulties.
\end{remark}
\medskip

\begin{remark}
The extended collision estimates for $\gamma>-3$ in Remark \ref{remark2.5} suggest that the present results may be extended to the wider soft-potential range $\max\{-3,-2-2s\}<\gamma<-2s$.
Indeed, for $-2<\gamma+2s<0$, the condition in \eqref{E-new-estimate} can be satisfied by choosing $\varrho<3/2$ sufficiently close to $3/2$ and $\vartheta>0$ sufficiently small.
Since the extended trilinear estimates in Remark \ref{remark2.5} involve additional velocity derivatives, closing the uniform energy estimates may require higher regularity of the initial data and suitable adjustments to the polynomial velocity weights.
\end{remark}
\medskip

\begin{remark}\label{remark-1.5}
Compared with earlier approaches that require a slightly larger exponent in the polynomial velocity weights, our estimates allow the endpoint $l=mN$, with $m$ defined in \eqref{weight function2}.
For $l\geq mN$ and $|\alpha|+|\beta|\leq N$, the polynomial weight exponents satisfy
\begin{align*}
l-(m+\gamma+2s)|\alpha|-m|\beta| \geq m(N-|\alpha|-|\beta|)-(\gamma+2s)|\alpha| \geq0.
\end{align*}
The additional velocity factors are controlled by the gaps between polynomial weight exponents at successive derivative orders and the dissipation generated by the time-dependent exponential factor, so that the condition $l\geq mN$ suffices.
\end{remark}
\medskip

\begin{remark}
If the strong convergence of the initial data in Theorem \ref{mainth2} is given with a prescribed rate, the estimates in Section \ref{Limit section} yield quantitative convergence
bounds on every finite time interval.
Combined with the uniform time-decay estimates \eqref{thm1 estimate 2}, these bounds also give a convergence rate uniform for all $t\geq0$.
The convergence rate obtained in this way is not optimal.
\end{remark}
\medskip

\subsection{Difficulties and Strategies}
\hspace*{\fill}

In this subsection, we first explain the main difficulty in the weighted energy estimates and describe three key techniques for obtaining global bounds uniform in $\varepsilon\in(0,1]$.
We then outline the additional arguments needed to establish strong convergence in the diffusive limit.

The main difficulty in the weighted energy estimates is to control the velocity growth in $E^\varepsilon\cdot vf^\varepsilon$.
In the Leibniz expansion of $\partial_\beta^\alpha(E^\varepsilon\cdot vf^\varepsilon)$ with $|\alpha|+|\beta|=N$, if at least one spatial derivative acts on $E^\varepsilon$ or a velocity derivative acts on $v$, the resulting terms involve lower-order derivatives of $f^\varepsilon$ and can be controlled using the stronger polynomial weights at those orders.
However, when all derivatives act on $f^\varepsilon$, both microscopic factors in the weighted energy estimate carry the same polynomial weight.
Controlling the factor $v$ then requires an additional weight $\langle v\rangle^{1/2}$ in the weighted $L^2$ norm of each microscopic factor.
The exponential factor $e^{q\langle v\rangle/(1+t)^\vartheta}$ used in earlier non-cutoff analyses produces dissipation with the same additional velocity weight $\langle v\rangle^{1/2}$ and a time coefficient proportional to $(1+t)^{-1-\vartheta}$.
Thus, this dissipation can control the weighted microscopic contribution in the electric field estimate if $\varepsilon(1+t)^{1+\vartheta}\|E^\varepsilon(t)\|_{L^\infty}$ remains sufficiently small for all $t\geq0$.
While the electric field decay in the two-species case is fast enough to meet this condition, the available decay estimate in the one-species case is
\begin{align*}
\|E^\varepsilon(t)\|_{L^\infty} \lesssim \delta_0(1+t)^{-\frac{3+2\varrho}{8}}.
\end{align*}
Since $\frac{3+2\varrho}{8}<1+\vartheta$, this estimate does not provide the uniform smallness required to absorb the microscopic contribution into the dissipation.

We overcome this difficulty by introducing the quadratic exponential factor $e^{q\langle v\rangle^2/(1+t)^\vartheta}$ in \eqref{weight function}.
Its time derivative produces an additional factor $\langle v\rangle$ in the dissipative norm.
Interpolation between this stronger velocity control and the collision dissipation then reduces the decay required of the electric field.
To use the quadratic exponential weight in the non-cutoff setting, we establish new weighted coercivity and trilinear estimates for $L$ and $\Gamma$, including estimates with velocity derivatives, in Section \ref{new-weighted-estimate}.
The key to these new weighted estimates is to move the exponential factor inside the collision operator using conservation of energy.
Writing $w_0=e^{q\langle v\rangle^2/(1+t)^\vartheta}$, we have
\begin{align*}
w_0(v)w_0(v_*)=w_0(v')w_0(v_*').
\end{align*}
Consequently, the trilinear operator $\mathcal T$ defined in \eqref{mathcalT-define} satisfies
\begin{align*}
w_0\mathcal T(f,g,h)=\mathcal T(w_0f,w_0g,w_0^{-1}h),
\end{align*}
as stated in \eqref{weighted T identity}.
For the bilinear collision operator $\Gamma$, this identity places the exponential weights on the two inputs and replaces $\mu^{1/2}$ by $w_0^{-1}\mu^{1/2}$.
The estimates for these modified Gaussian factors in Lemma \ref{lem:profiletransfer} then allow us to avoid the additional velocity-weight loss arising from direct estimates of exponential weight differences, as in \cite[Lemma 2.3 and Remark 2.3]{DLYZ2013}.
The same argument applies to velocity derivatives by first using \eqref{mathcalT-derivative} and then transferring the exponential weight.
For the linearized collision operator $L$, the remaining commutator is controlled by a small multiple of the weighted collision dissipation without additional velocity weights on the unknown; see \eqref{exponential-L-commutator}.

\smallskip

To derive the uniform energy estimates, we use the following three techniques:
\smallskip

\noindent{\bf 1. Choice of the polynomial velocity weights.}

To control the transport commutator generated by $\frac{1}{\varepsilon} v\cdot\nabla_x f^\varepsilon$, we choose the polynomial weights $\langle v\rangle^{l-(m+\gamma+2s)|\alpha|-m|\beta|}$ in \eqref{weight function2}.
The difference $\gamma+2s$ between the spatial and velocity derivative penalties ensures that
\begin{align*}
w_l(\alpha,\beta)=\langle v\rangle^{\gamma+2s} w_l(\alpha+e_i,\beta-e_i), \quad \text{for}\;\; |\beta|\geq1.
\end{align*}
This identity matches the velocity weights in the transport commutator with those provided by the collision dissipation, allowing us to control these terms through a suitable linear combination of the weighted energy estimates.
\smallskip

\noindent{\bf 2. Interpolation for the electric field term.}

The dissipation in \eqref{dissipation functional} controls the weighted unknown with an additional factor $\langle v\rangle$  through the time derivative of the quadratic exponential weight, and with a factor $\langle v\rangle^{(\gamma+2s)/2}$ through the collision operator.
Interpolating between these two bounds gives the velocity weight $\langle v\rangle^{1/2}$ needed to control the microscopic contribution in the electric field estimate.
In particular, for mixed derivatives, \eqref{E-new-estimate} gives
\begin{align*}
\varepsilon\|E^\varepsilon\|_{L^\infty}\left\|\langle v\rangle^{\frac{1}{2}}\widetilde w_l(\alpha,\beta)\partial_\beta^\alpha\mathbf P^\perp f^\varepsilon\right\|^2
\lesssim (1+t)^{\frac{(1+\vartheta)(1-\gamma-2s)}{2-\gamma-2s}}\|E^\varepsilon\|_{L^\infty}\widetilde{\mathcal D}_{N,l}(t)
\lesssim \delta_0\widetilde{\mathcal D}_{N,l}(t).
\end{align*}
Here, the last inequality follows from the available electric field decay and the parameter choice in \eqref{mainth1 assumption}, which ensures that
\begin{align*}
\frac{(1+\vartheta)(1-\gamma-2s)}{2-\gamma-2s}
\leq \frac{3+2\varrho}{8}.
\end{align*}
Thus, the quadratic exponential weight allows us to close this estimate under the slower electric field decay in the one-species case.
\smallskip

\noindent{\bf 3. Matching velocity derivatives and polynomial weights.}

The magnetic term $(v\times B^\varepsilon)\cdot\nabla_v f^\varepsilon$ contains an additional velocity derivative, while higher velocity derivatives carry weaker polynomial weights.
To control this term, we combine Sobolev duality in $v$ with weighted interpolation.
Let $\sigma=\min\{s,1/2\}$. For $0<\sigma<1/2$ and $|\alpha|\geq1$, we use the interpolation estimate
\begin{align*}
\left|\langle v\rangle^{1-\frac{\gamma+2s}{2}+\sigma}w_l(\alpha,\beta)g\right|_{H_v^{1-\sigma}}
\lesssim
\left|\langle v\rangle^{\frac{\gamma+2s}{2}-\sigma}w_l(|\alpha|-1,|\beta|)g\right|_{H_v^\sigma}^{2\sigma}
\left|\langle v\rangle^{\frac{\gamma+2s}{2}-\sigma}w_l(|\alpha|-1,|\beta|+1)g\right|_{H_v^{1+\sigma}}^{1-2\sigma}.
\end{align*}
To match the polynomial weights on both sides, we choose them to satisfy
\begin{align*}
\langle v\rangle^{\frac{\gamma+2s}{2}-\sigma}
\overline w_l(|\alpha|-1,|\beta|)^{2\sigma}
\overline w_l(|\alpha|-1,|\beta|+1)^{1-2\sigma}
=\langle v\rangle^{1-\frac{\gamma+2s}{2}+\sigma}\overline w_l(\alpha,\beta).
\end{align*}
By \eqref{weight function2}, this identity is equivalent to $2\sigma(m-1)=1-2(\gamma+2s)$, which leads to the choice
\begin{align*}
m=1+\frac{1-2(\gamma+2s)}{\min\{2s,1\}}.
\end{align*}
When $\sigma=1/2$ (equivalently, $s\geq1/2$), the required fractional velocity regularity follows directly from the collision dissipation.
Together with the polynomial weights above, these estimates control the additional velocity derivative in the magnetic term;  see \eqref{weighted:top-mixed-magnetic-interpolation1} and \eqref{weighted:top-mixed-magnetic-interpolation}.
\smallskip

The uniform estimates obtained using the three techniques above provide the starting point for the diffusive limit.
To establish strong convergence, we derive additional estimates for the acoustic components $\rho^\varepsilon+\theta^\varepsilon$ and $\mathcal P^\perp u^\varepsilon$, and the microscopic component
$\mathbf P^\perp f^\varepsilon$. These estimates are established in Section \ref{Limit section} and supplement the convergence argument in \cite{WZGX2025SIMA}.

For well-prepared initial data, energy estimates for the corrected acoustic variables defined in \eqref{limit:corrected-acoustic-variables} give
\begin{align*}
(\rho^\varepsilon+\theta^\varepsilon,\mathcal P^\perp u^\varepsilon)
\to(0,0) \quad\text{strongly in }C([0,T];H_x^{N-1})
\end{align*}
as $\varepsilon\to0$; see \eqref{limit:acoustic-strong-C}.
A separate microscopic energy estimate yields
\begin{align*}
\mathbf P^\perp f^\varepsilon\to0
\quad\text{strongly in }C([0,T];H_x^{N-1}L_v^2)
\end{align*}
as $\varepsilon\to0$; see \eqref{limit:microscopic-strong-C}.
We then estimate the difference between the macroscopic variables and the limiting solution to obtain strong convergence in the whole space $\mathbb{R}^3$, as in \eqref{limit:strong-global-space}.
Finally, the uniform time-decay estimates make the convergence uniform for all $t\geq0$; see \eqref{limit:strong-global-time}.
\smallskip

The rest of this paper is organized as follows. In Section \ref{new-weighted-estimate}, we establish new weighted estimates for the collision operators $L$ and $\Gamma$ associated with the quadratic exponential weight.
In Section \ref{The a Priori Estimate}, we derive a priori estimates uniform in $\varepsilon$ using the weighted energy method.
In Section \ref{Global Existence}, we obtain time-decay estimates, close the a priori assumption, and prove Theorem \ref{mainth1}.
In Section \ref{Limit section}, we use these uniform estimates to justify the incompressible NSFM limit \eqref{INSFM limit} and prove Theorem \ref{mainth2}.
Finally, Appendices \ref{appendix:auxiliary} and \ref{appendix:inequalities} collect auxiliary estimates for the collision operators and the functional inequalities used in the proofs, respectively.
\medskip

\section{New Weighted Estimates for the Collision Operators}\label{new-weighted-estimate}

In this section, we establish new weighted estimates for the collision operators $L$ and $\Gamma$.
We then apply these estimates, together with the new time-velocity weight function introduced in \eqref{weight function1}, to control the nonlinear collision term $\Gamma(f^\varepsilon,f^\varepsilon)$ in the energy estimates for the VMB system \eqref{rVMB}.

\subsection{Weighted Estimates for $L$ and $\Gamma$}
\hspace*{\fill}

We first introduce the following time-velocity weight function
\begin{align*}
  w_\ell=w_\ell(t,v):=\langle v \rangle^\ell e^{\frac{q\langle v \rangle^2}{(1+t)^\vartheta}},
\end{align*}
where $\ell\geq0$, $\vartheta>0$ and $0\leq q \ll 1/4$.
To distinguish the polynomial factor from the exponential factor, we sometimes write $w_\ell=\langle v\rangle^{\ell}w_0$.

It is convenient to express both collision operators in terms of the trilinear operator
\begin{equation}\label{mathcalT-define}
\mathcal T(f,g,h)(v)
=\iint_{\mathbb R^3\times\mathbb S^2}B(v-v_*,\sigma)h(v_*)\big(f(v_*')g(v')-f(v_*)g(v)\big)\d v_*\d\sigma.
\end{equation}
With this notation, the bilinear collision operator $\Gamma$ and the linearized collision operator $L$ can be written as
\begin{align*}
\Gamma(f,g)=\mathcal T(f,g,\mu^{1/2}),\quad
Lf=-\Gamma(\mu^{1/2},f)-\Gamma(f,\mu^{1/2})
=-\mathcal T(\mu^{1/2},f,\mu^{1/2})-\mathcal T(f,\mu^{1/2},\mu^{1/2}).
\end{align*}
The translation invariance of $\mathcal T$ further gives
\begin{align}\label{mathcalT-derivative}
&\partial_\beta\Gamma(f,g)
=\sum_{\beta_1+\beta_2+\beta_3=\beta}C_{\beta}^{\beta_1,\beta_2,\beta_3}
\mathcal T\left(\partial_{\beta_1}f,\partial_{\beta_2}g,\partial_{\beta_3}\mu^{1/2}\right),\\
&\partial_\beta Lf
=-\!\!\!\sum_{\beta_1+\beta_2+\beta_3=\beta}\!\!C_{\beta}^{\beta_1,\beta_2,\beta_3}
\Big[\mathcal T\left(\partial_{\beta_1}\mu^{1/2},\partial_{\beta_2}f,\partial_{\beta_3}\mu^{1/2}\right)
+\mathcal T\left(\partial_{\beta_1}f,\partial_{\beta_2}\mu^{1/2},\partial_{\beta_3}\mu^{1/2}\right)\Big].\label{L-derivative}
\end{align}

These representations motivate the following weighted estimates for $\mathcal T$ with Gaussian profiles in its third argument.
\begin{lemma}\label{lem:profiletransfer}
Let $0<s<1$,
$\gamma>\max\{-3,-\frac32-2s\}$,
$\ell\geq0$, and $0<q_0<1/4$.
Suppose that $0\leq q_1(t)\leq q_0$. Then
\begin{align}
\begin{split}\label{profile zero}
\left|\left\langle\langle v\rangle^\ell\mathcal T(f,g,e^{-q_1(t)\langle v\rangle^2}\mu^{1/2}),h\right\rangle\right|
\lesssim\;& \Big\{|f|_{L^2_{\gamma/2+s}}|\langle v\rangle^\ell g|_D
+|f|_D|\langle v\rangle^\ell g|_{L^2_{\gamma/2+s}}\Big\}|h|_D \\
&+\min\Big\{|f|_{L^2}|\langle v\rangle^\ell g|_{L^2_{\gamma/2+s}},|f|_{L^2_{\gamma/2+s}}|\langle v\rangle^\ell g|_{L^2}\Big\}|h|_D.
\end{split}
\end{align}
More generally, for any multi-index $\beta$, we have
\begin{align}
\begin{split}\label{profile derivative}
\left|\left\langle\langle v\rangle^\ell\mathcal T(f,g,e^{-q_1(t)\langle v\rangle^2}\partial_\beta\mu^{1/2}),h\right\rangle\right|
\lesssim\;&
\Big\{|f|_{L^2_{\gamma/2+s}}|\langle v\rangle^\ell g|_D
+|f|_D|\langle v\rangle^\ell g|_{L^2_{\gamma/2+s}}\Big\}|h|_D \\
&+\min\Big\{|f|_{L^2}|\langle v\rangle^\ell g|_{L^2_{\gamma/2+s}},|f|_{L^2_{\gamma/2+s}}|\langle v\rangle^\ell g|_{L^2}\Big\}|h|_D.
\end{split}
\end{align}
The implicit constants are independent of $t$ and $q_1(t)$.
\end{lemma}

\begin{proof}
Fix $t\geq0$ and write $q_1=q_1(t)$. Define
\begin{align*}
\lambda:=\sqrt{1+4q_1},
\qquad \widetilde\mu(v):=e^{-2q_1\langle v\rangle^2}\mu(v),
\qquad \widetilde f(y):=f(y/\lambda),
\end{align*}
and define $\widetilde g$ and $\widetilde h$ similarly. Clearly, we have
\begin{align*}
1\leq\lambda\leq\sqrt{1+4q_0}<\sqrt2,
\qquad \widetilde\mu(v)=e^{-2q_1}(2\pi)^{6q_1}\mu(v)^{1+4q_1}=e^{-2q_1}\mu(\lambda v).
\end{align*}
Let $|\cdot|_{\widetilde D}$ denote $|\cdot|_D$ with $\mu$ replaced by $\widetilde\mu$.
Lemma \ref{lem:Gaussian-norm-comparison} gives
\begin{equation*}
|f|_{\widetilde D}\approx |f|_D.
\end{equation*}
A change of variables yields
\begin{equation*}
|f|_{\widetilde D}^2=e^{-2q_1}\lambda^{-\gamma-6}|\widetilde f|_D^2.
\end{equation*}
Consequently, one has
\begin{align}\label{D scaling}
|\widetilde f|_D\approx |f|_D.
\end{align}
For any fixed $\ell\in\mathbb R$, we also have
\begin{align}\label{poly-scaling-norm}
|\widetilde f|_{L^2_\ell}^2=\int_{\mathbb R^3}\langle y\rangle^{2\ell}|\widetilde f(y)|^2\d y
=\lambda^3\int_{\mathbb R^3}\langle\lambda v\rangle^{2\ell}|f(v)|^2\d v\approx |f|_{L^2_\ell}^2.
\end{align}
Define
\begin{align*}
m_\lambda(y) :=\frac{\langle y/\lambda \rangle^\ell}{\langle y\rangle^\ell}.
\end{align*}
For $1\leq\lambda<\sqrt2$, a direct calculation shows that the $W^{1,\infty}$ norms of $m_\lambda$ and $m_\lambda^{-1}$ are bounded by a constant depending only on $\ell$ and independent of $\lambda$.
Lemma \ref{lem:multiplier} therefore gives
\begin{align*}
\big|m_\lambda\widetilde h\big|_D \lesssim\big|\widetilde h\big|_D,
\qquad
\left|\langle y\rangle^\ell\widetilde g\right|_D
\lesssim \left|\langle y/\lambda\rangle^\ell\widetilde g\right|_D.
\end{align*}
Combining these estimates with \eqref{polynomial-Gamma}, \eqref{D scaling} and \eqref{poly-scaling-norm}, we obtain
\begin{align}\label{scaled polynomial estimate}
\begin{split}
\left|\left\langle\langle y/\lambda\rangle^\ell\Gamma(\widetilde f,\widetilde g),\widetilde h\right\rangle_{L^2_y}\right|
=\;&\left|\left\langle\langle y\rangle^\ell\Gamma(\widetilde f,\widetilde g),m_\lambda(y)\widetilde h\right\rangle_{L^2_y}\right|\\
\lesssim\;&
\Big\{\big|\widetilde f\big|_{L^2_{\gamma/2+s}}\big|\langle y/\lambda\rangle^\ell\widetilde g\big|_D
+\big|\widetilde f\big|_D\big|\langle y/\lambda\rangle^\ell\widetilde g\big|_{L^2_{\gamma/2+s}}\Big\}\big|\widetilde h\big|_D\\
&+ \min\Big\{ \big|\widetilde f\big|_{L^2} \big|\langle y/\lambda\rangle^\ell \widetilde g\big|_{L^2_{\gamma/2+s}},
\big|\widetilde f\big|_{L^2_{\gamma/2+s}} \big|\langle y/\lambda\rangle^\ell \widetilde g\big|_{L^2} \Big\}\big|\widetilde h\big|_D \\
\approx\;& \Big\{|f|_{L^2_{\gamma/2+s}}|\langle v\rangle^\ell g|_D
+|f|_D|\langle v\rangle^\ell g|_{L^2_{\gamma/2+s}}\Big\}|h|_D \\
&+\min\Big\{|f|_{L^2}|\langle v\rangle^\ell g|_{L^2_{\gamma/2+s}},|f|_{L^2_{\gamma/2+s}}|\langle v\rangle^\ell g|_{L^2}\Big\}|h|_D.
\end{split}
\end{align}
Hence, \eqref{profile zero} follows from \eqref{scaled polynomial estimate} and the identity
\begin{align*}
\left\langle\langle v\rangle^\ell\mathcal T(f,g,e^{-q_1\langle v\rangle^2}\mu^{1/2}),h\right\rangle
=e^{-q_1}\lambda^{-\gamma-6}\left\langle\langle y/\lambda\rangle^\ell\Gamma(\widetilde f,\widetilde g),\widetilde h\right\rangle_{L^2_y}.
\end{align*}

To prove \eqref{profile derivative}, we first note that
\begin{align*}
\mu^{1/4}(v)=(2\pi)^{3/8}\mu^{1/2}(v/\sqrt2).
\end{align*}
A scaling argument similar to that used to derive \eqref{profile zero}, together with Lemma \ref{lem:Gaussian-norm-comparison}, yields
\begin{align}\label{profile quarter}
\begin{split}
\left|\left\langle\langle v\rangle^\ell\mathcal T(f,g,\mu^{1/4}),h\right\rangle\right|
\lesssim\;&
\Big\{|f|_{L^2_{\gamma/2+s}}|\langle v\rangle^\ell g|_D
+|f|_D|\langle v\rangle^\ell g|_{L^2_{\gamma/2+s}}\Big\}|h|_D \\
&+\min\Big\{|f|_{L^2}|\langle v\rangle^\ell g|_{L^2_{\gamma/2+s}},|f|_{L^2_{\gamma/2+s}}|\langle v\rangle^\ell g|_{L^2}\Big\}|h|_D.
\end{split}
\end{align}
By the definition of $\mathcal T$ in \eqref{mathcalT-define} and the Fourier inversion formula, we have
\begin{align*}
&\mathcal T \left(f,g,e^{-q_1\langle v\rangle^2}\partial_\beta\mu^{1/2}\right)(v)\\
=\;&\int_{\mathbb R^3\times\mathbb S^2}
B(v-v_*,\sigma)e^{-q_1\langle v_*\rangle^2}(\partial_\beta\mu^{1/2})(v_*)\big(f(v_*')g(v')-f(v_*)g(v)\big)\d v_*\d\sigma\\
=\;&\frac1{(2\pi)^{3/2}}
\int_{\mathbb R^3\times\mathbb S^2}B(v-v_*,\sigma)\mu^{1/4}(v_*)\left[\int_{\mathbb R^3}\widehat{\mathcal{A}}(\xi)e^{iv_*\cdot\xi}\d\xi\right]
\big(f(v_*')g(v')-f(v_*)g(v)\big)\d v_*\d\sigma\\
=\;&\frac1{(2\pi)^{3/2}}
\int_{\mathbb R^3}\widehat{\mathcal{A}}(\xi)\mathcal T(f,g,e^{iv\cdot\xi}\mu^{1/4})(v)\d\xi,
\end{align*}
where
\begin{align*}
\widehat{\mathcal{A}}(\xi)
=\mathcal F_v\left(e^{-q_1\langle v\rangle^2}\mu^{-1/4}\partial_\beta\mu^{1/2}\right)(\xi).
\end{align*}
Taking the $L_v^2$ inner product with $\langle v\rangle^\ell h$, we obtain
\begin{align}
\begin{split}\label{mathcalTequality}
&\left\langle \langle v\rangle^{\ell} \mathcal T \left(f,g,e^{-q_1\langle v\rangle^2}\partial_\beta\mu^{1/2}\right), h \right\rangle\\
=\;& \frac1{(2\pi)^{3/2}} \int_{\mathbb R^3} \widehat{\mathcal{A}}(\xi) \left\langle \langle v\rangle^\ell\mathcal T(f,g,e^{iv\cdot\xi}\mu^{1/4}),h\right\rangle \d \xi \\
=\;&\frac1{(2\pi)^{3/2}} \int_{\mathbb R^3} \widehat{\mathcal{A}}(\xi) \left\langle
\langle v\rangle^\ell\mathcal T(e^{iv\cdot\xi}f,e^{iv\cdot\xi}g,\mu^{1/4}),e^{iv\cdot\xi}h\right\rangle \d \xi.
\end{split}
\end{align}
The last equality follows from conservation of momentum, $v+v_*=v'+v_*'$, which gives
\begin{equation*}
\mathcal T(f,g,e^{iv\cdot\xi}\mu^{1/4})=e^{-iv\cdot\xi}\mathcal T(e^{iv\cdot\xi}f,e^{iv\cdot\xi}g,\mu^{1/4}).
\end{equation*}
Applying Lemma \ref{lem:multiplier} with $m(v)=e^{iv\cdot\xi}$, we obtain
\begin{equation}\label{D modulation}
|e^{iv\cdot\xi}f|_D\lesssim\langle\xi\rangle^s|f|_D, \qquad
|e^{iv\cdot\xi}f|_{L^2_\ell}=|f|_{L^2_\ell}.
\end{equation}
It follows from \eqref{profile quarter} and \eqref{D modulation} that
\begin{align}\label{modulated profile estimate}
\begin{split}
&\Big|\left\langle\langle v\rangle^\ell\mathcal T(e^{iv\cdot\xi}f,e^{iv\cdot\xi}g,\mu^{1/4}),e^{iv\cdot\xi}h\right\rangle\Big|\\
\lesssim\;&
\langle\xi\rangle^{2s}\Big\{|f|_{L^2_{\gamma/2+s}}|\langle v\rangle^\ell g|_D
+|f|_D|\langle v\rangle^\ell g|_{L^2_{\gamma/2+s}}\Big\}|h|_D \\
&+\langle\xi\rangle^{2s} \min\Big\{|f|_{L^2}|\langle v\rangle^\ell g|_{L^2_{\gamma/2+s}},|f|_{L^2_{\gamma/2+s}}|\langle v\rangle^\ell g|_{L^2}\Big\}|h|_D.
\end{split}
\end{align}
Writing
\begin{align*}
\partial_\beta\mu^{1/2}(v)
=\sum_{\beta'\leq\beta} C_{\beta,\beta'}v^{\beta'}\mu^{1/2}(v),
\end{align*}
we obtain by a direct computation
\begin{equation}\label{profile Fourier bound}
\int_{\mathbb R^3}\langle\xi\rangle^{2s}\left| \widehat{\mathcal{A}}(\xi)\right|\d\xi
=\frac{1}{(2\pi)^{3/8}}\frac{e^{-q_1}}{(2q_1+1/4)^{3/2}}
\int_{\mathbb R^3}\langle\xi\rangle^{2s}
\Big|\sum_{\beta'\leq\beta}C_{\beta,\beta'}i^{|\beta'|}\partial_\xi^{\beta'} e^{-\frac{|\xi|^2}{4q_1+1/2}}\Big|\d\xi
\lesssim 1.
\end{equation}
Combining \eqref{mathcalTequality}, \eqref{modulated profile estimate} and \eqref{profile Fourier bound} proves \eqref{profile derivative}.
\end{proof}
\smallskip

The following consequence of the preceding lemma will be used in the proof of Lemma \ref{lem:weighted-L}.
\begin{lemma}\label{lem:Gaussian-linear-T}
Let $0<s<1$,
$\gamma>\max\{-3,-\frac32-2s\}$,
$\ell\geq0$, and $0<q_0<1/4$.
Let $p_1,p_2$ be fixed polynomials.
For $0\leq q_1(t)\leq q_0$, define
\begin{align*}
\phi_{q_1}(v)=e^{q_1(t)\langle v\rangle^2}p_1(v)\mu^{1/2}(v),
\qquad
\psi_{q_1}(v)=e^{-q_1(t)\langle v\rangle^2}p_2(v)\mu^{1/2}(v).
\end{align*}
Then, we have
\begin{align}\label{Gaussian-linear-T}
&\left|\left\langle\langle v\rangle^\ell\mathcal T(\phi_{q_1},g,\psi_{q_1}),h\right\rangle\right|
+\left|\left\langle\langle v\rangle^\ell\mathcal T(g,\phi_{q_1},\psi_{q_1}),h\right\rangle\right|
\lesssim |\langle v\rangle^\ell g|_D|h|_D.
\end{align}
The implicit constant is independent of $t$ and $q_1(t)$.
\end{lemma}

\begin{proof}
Since $0\leq q_1(t)\leq q_0<1/4$, Gaussian decay gives
\begin{align*}
|\phi_{q_1}|_{L^2}+|\phi_{q_1}|_D+|\langle v\rangle^\ell\phi_{q_1}|_{L^2} +|\langle v\rangle^\ell\phi_{q_1}|_D
\lesssim1.
\end{align*}
Moreover,
\begin{align*}
\psi_{q_1}(v)=e^{-q_1(t)\langle v\rangle^2} \sum_\beta C_\beta\partial_\beta \mu^{1/2}(v),
\end{align*}
where the sum is finite and the coefficients depend only on $p_2$.
Thus, applying \eqref{profile derivative} and using \eqref{norm inequality} and Lemma \ref{lem:multiplier}, we obtain
\begin{align*}
&\left|\left\langle\langle v\rangle^\ell\mathcal T(\phi_{q_1},g,\psi_{q_1}),h\right\rangle\right|
+\left|\left\langle\langle v\rangle^\ell\mathcal T(g,\phi_{q_1},\psi_{q_1}),h\right\rangle\right|
\lesssim \big(|\langle v\rangle^\ell g|_D+|g|_D\big)|h|_D
\lesssim |\langle v\rangle^\ell g|_D|h|_D.
\end{align*}
This proves \eqref{Gaussian-linear-T}.
\end{proof}

Based on the preceding lemmas, we establish the following new weighted estimates for $\Gamma$ and $L$.

\begin{lemma}\label{lem:Gamma-weighted}
Let $0<s<1$,
$\gamma>\max\{-3,-\frac32-2s\}$ and $\ell\geq0$.
Then
\begin{align}
\begin{split}\label{Gamma-weighted-1}
\left|\left\langle\Gamma(f,g),w_\ell^2h\right\rangle\right|
\lesssim\;&
\Big\{|w_0f|_{L^2_{\gamma/2+s}}|w_\ell g|_D+|w_0f|_D|w_\ell g|_{L^2_{\gamma/2+s}}\\
&+\min\big\{|w_0f|_{L^2}|w_\ell g|_{L^2_{\gamma/2+s}},|w_0f|_{L^2_{\gamma/2+s}}|w_\ell g|_{L^2}\big\}\Big\}|w_\ell h|_D,
\end{split}
\end{align}
and for $|\beta|\geq 1$,
\begin{align}\label{Gamma-weighted-2}
\left|\left\langle\partial_\beta\Gamma(f,g),w_\ell^2\partial_\beta h\right\rangle\right|
\lesssim\;&
\sum_{\beta_1+\beta_2\leq\beta}
\Big\{|w_0\partial_{\beta_1}f|_{L^2_{\gamma/2+s}}|w_\ell\partial_{\beta_2}g|_D
+|w_0\partial_{\beta_1}f|_D|w_\ell\partial_{\beta_2}g|_{L^2_{\gamma/2+s}} \\
& +\min\big\{|w_0\partial_{\beta_1}f|_{L^2}|w_\ell\partial_{\beta_2}g|_{L^2_{\gamma/2+s}},
|w_0\partial_{\beta_1}f|_{L^2_{\gamma/2+s}}|w_\ell\partial_{\beta_2}g|_{L^2}\big\}\Big\}|w_\ell\partial_\beta h|_D. \nonumber
\end{align}
\end{lemma}

\begin{proof}
The proof uses the estimates for $\mathcal T$ established in Lemma \ref{lem:profiletransfer}.
Conservation of energy $|v|^2+|v_*|^2=|v'|^2+|v'_*|^2$ gives
\begin{equation*}
w_0(v)w_0(v_*)=w_0(v')w_0(v_*').
\end{equation*}
Recalling the definition of $\mathcal T$ in \eqref{mathcalT-define}, we have
\begin{align}\label{weighted T identity}
w_0\mathcal T(f,g,h)=\mathcal T(w_0f,w_0g,w_0^{-1}h).
\end{align}
Therefore, by \eqref{profile zero} and the identity $\Gamma(f,g)=\mathcal{T}(f,g,\mu^{1/2})$, one has
\begin{align*}
\left|\left\langle\Gamma(f,g),w_\ell^2h\right\rangle\right|
=\;&\left|\left\langle\langle v\rangle^\ell\mathcal T\left(w_0f,w_0g,e^{-\frac{q\langle v\rangle^2}{(1+t)^{\vartheta}}}\mu^{1/2}\right),w_\ell h\right\rangle\right|\\
\lesssim\;&
\Big\{|w_0f|_{L^2_{\gamma/2+s}}|\langle v\rangle^\ell w_0g|_D
+|w_0f|_D|\langle v\rangle^\ell w_0g|_{L^2_{\gamma/2+s}}\\
&+\min\big\{|w_0f|_{L^2}|\langle v\rangle^\ell w_0 g|_{L^2_{\gamma/2+s}},
|w_0f|_{L^2_{\gamma/2+s}}|\langle v\rangle^\ell w_0 g|_{L^2}\big\} \Big\}|w_\ell h|_D,
\end{align*}
which proves \eqref{Gamma-weighted-1}.

For $|\beta|\geq1$, by \eqref{mathcalT-derivative} and \eqref{weighted T identity}, we have
\begin{align*}
w_0\partial_\beta\Gamma(f,g)
=\sum_{\beta_1+\beta_2+\beta_3=\beta}C_\beta^{\beta_1,\beta_2,\beta_3}
\mathcal T\left(w_0\partial_{\beta_1}f,w_0\partial_{\beta_2}g,e^{-\frac{q\langle v\rangle^2}{(1+t)^{\vartheta}}}\partial_{\beta_3}\mu^{1/2}\right).
\end{align*}
Applying \eqref{profile derivative}, we obtain
\begin{align*}
\left|\left\langle\partial_\beta\Gamma(f,g),w_\ell^2\partial_\beta h\right\rangle\right|
\lesssim\;& \sum_{\beta_1+\beta_2+\beta_3=\beta}
\left|\left\langle\langle v\rangle^\ell\mathcal T\left(w_0\partial_{\beta_1}f,w_0\partial_{\beta_2}g,
e^{-\frac{q\langle v\rangle^2}{(1+t)^{\vartheta}}}\partial_{\beta_3}\mu^{1/2}\right),
w_\ell\partial_\beta h\right\rangle\right|\\
\lesssim\;&
\sum_{\beta_1+\beta_2\leq\beta}
\Big\{|w_0\partial_{\beta_1}f|_{L^2_{\gamma/2+s}}|w_\ell\partial_{\beta_2}g|_D
+|w_0\partial_{\beta_1}f|_D|w_\ell\partial_{\beta_2}g|_{L^2_{\gamma/2+s}} \\
& +\min\big\{|w_0\partial_{\beta_1}f|_{L^2}|w_\ell\partial_{\beta_2}g|_{L^2_{\gamma/2+s}},
|w_0\partial_{\beta_1}f|_{L^2_{\gamma/2+s}}|w_\ell\partial_{\beta_2}g|_{L^2}\big\}\Big\}|w_\ell\partial_\beta h|_D.
\end{align*}
This proves \eqref{Gamma-weighted-2}.
\end{proof}
\smallskip

\begin{lemma}\label{lem:weighted-L}
Let $0<s<1$ and $\gamma>\max\{-3,-\frac32-2s\}$.
Then
\begin{equation}\label{L-coercive-weighted-1}
\left\langle w_\ell^2Lf,f\right\rangle
\gtrsim|w_\ell f|_D^2-C|f|_{L^2(B_C)}^2,
\end{equation}
and for $|\beta|\geq1$,
\begin{align}
\left\langle w_\ell^2\partial_\beta Lf,\partial_\beta f\right\rangle
\gtrsim|w_\ell\partial_\beta f|_D^2-C\sum_{|\beta'|<|\beta|}|w_\ell\partial_{\beta'}f|_D^2. \label{L-coercive-weighted-2}
\end{align}
\end{lemma}

\begin{proof}
We first prove \eqref{L-coercive-weighted-1}.
Although $q_1=q(1+t)^{-\vartheta}$, we treat $q_1$ as an independent parameter in $[0,q_0]$ in the following calculation and write $w_0=e^ {q_1\langle v\rangle^2}$.
For any function $F(v)$, held fixed as $q_1$ varies, the identity
\begin{align*}
LF=-\Gamma(\mu^{1/2},F)-\Gamma(F,\mu^{1/2})=-\mathcal T(\mu^{1/2},F,\mu^{1/2})-\mathcal T(F,\mu^{1/2},\mu^{1/2})
\end{align*}
and \eqref{weighted T identity} yield
\begin{align*}
w_0 L(w_0^{-1}F)=-\mathcal T(w_0\mu^{1/2},F,w_0^{-1}\mu^{1/2})-\mathcal T(F,w_0\mu^{1/2},w_0^{-1}\mu^{1/2}).
\end{align*}
Differentiating this identity with respect to $q_1$ gives
\begin{align*}
\frac{\d}{\d q_1}
\left[ w_0L(w_0^{-1}F)\right]
=\;&
-\mathcal T(\langle v\rangle^2w_0\mu^{1/2},F,w_0^{-1}\mu^{1/2})
+\mathcal T(w_0\mu^{1/2},F,\langle v\rangle^2w_0^{-1}\mu^{1/2})\\
&-\mathcal T(F,\langle v\rangle^2w_0\mu^{1/2},w_0^{-1}\mu^{1/2})
+\mathcal T(F,w_0\mu^{1/2},\langle v\rangle^2w_0^{-1}\mu^{1/2}),
\end{align*}
which, combined with Lemma \ref{lem:Gaussian-linear-T}, implies
\begin{align*}
&\left|\frac{\d}{\d q_1}\left\langle\langle v\rangle^\ell w_0L(w_0^{-1}F),
\langle v\rangle^\ell F\right\rangle\right|
\lesssim |\langle v\rangle^\ell F|_D^2.
\end{align*}
Integrating with respect to $q_1$ from $0$ to $q(1+t)^{-\vartheta}$, and then taking $F=w_0f$, we obtain
\begin{equation}\label{exponential-L-commutator}
\left|\left\langle\langle v\rangle^\ell\big(w_0Lf-L(w_0f)\big),w_\ell f\right\rangle\right|
\lesssim \frac{q}{(1+t)^\vartheta}|w_\ell f|_D^2.
\end{equation}
By \eqref{polynomial-L-coercivity} and \eqref{exponential-L-commutator}, we have
\begin{align*}
\left\langle Lf,w_\ell^2f\right\rangle
=\;& \left\langle L(w_0f),\langle v\rangle^{2\ell}w_0f
\right\rangle
+\left\langle\langle v\rangle^\ell\big(w_0Lf-L(w_0f)\big),w_\ell f\right\rangle\\
\geq\;& (c_0-Cq)|w_\ell f|_D^2-C|w_0f|_{L^2(B_C)}^2\\
\geq\;& (c_0-Cq)|w_\ell f|_D^2-C|f|_{L^2(B_C)}^2.
\end{align*}
Choosing $q$ sufficiently small gives \eqref{L-coercive-weighted-1}.

We now prove \eqref{L-coercive-weighted-2}.
Let $|\beta|\geq1$. By \eqref{L-derivative}, \eqref{weighted T identity}, and Lemma \ref{lem:Gaussian-linear-T}, we have, for any $\eta>0$,
\begin{align*}
\left|\left\langle w_\ell\big(\partial_\beta Lf-L\partial_\beta f\big),w_\ell\partial_\beta f\right\rangle\right|
\lesssim\;&\sum_{|\beta'|<|\beta|}|w_\ell\partial_{\beta'}f|_D|w_\ell\partial_\beta f|_D\\
\lesssim\;& \eta |w_\ell\partial_\beta f|_D^2+\sum_{|\beta'|<|\beta|}|w_\ell\partial_{\beta'}f|_D^2.
\end{align*}
By local Sobolev interpolation and \eqref{norm inequality}, for any $\eta>0$, we have
\begin{align*}
|\partial_\beta f|_{L^2(B_C)}^2
\lesssim \eta|w_\ell\partial_\beta f|_D^2
+\sum_{|\beta'|<|\beta|} |w_\ell\partial_{\beta'}f|_D^2.
\end{align*}
Applying \eqref{L-coercive-weighted-1} to $\partial_\beta f$, using the preceding estimates, and choosing $\eta$ sufficiently small, we obtain
\begin{align*}
\left\langle w_\ell^2\partial_\beta Lf,\partial_\beta f\right\rangle
=\;&\left\langle L\partial_\beta f,w_\ell^2\partial_\beta f\right\rangle
+\left\langle w_\ell(\partial_\beta Lf-L\partial_\beta f),w_\ell\partial_\beta f\right\rangle\\
\gtrsim\;& |w_\ell\partial_\beta f|_D^2-C\sum_{|\beta'|<|\beta|}|w_\ell\partial_{\beta'}f|_D^2.
\end{align*}
This proves \eqref{L-coercive-weighted-2}.
\end{proof}
\smallskip

\begin{remark}\label{remark2.5}
By adapting the estimates in \cite[Proposition 6.1]{GS2011}, the weighted coercivity estimates \eqref{L-coercive-weighted-1}--\eqref{L-coercive-weighted-2} can be extended to $\gamma>-3$ and $0<s<1$ for sufficiently small $q>0$, with constants independent of $t$.

The trilinear estimates can also be extended to this range by allowing up to two additional velocity derivatives on either input factor.
This follows by adapting the arguments in Sections 4 and 6 of \cite{GS2011}, including Lemma 6.1, and combining \eqref{weighted T identity} with the corresponding weighted derivative estimates.
More precisely, for any $\ell,m\geq0$ and any multi-index $\beta$, we have
\begin{align*}
\left|\left\langle\partial_\beta\Gamma(f,g),w_\ell^2\partial_\beta h\right\rangle\right|
\lesssim\;&
\sum_{\beta_1+\beta_2\leq\beta}\big\{|w_0\partial_{\beta_1}f|_{L_v^2}|w_\ell\partial_{\beta_2}g|_{N^s_\gamma}
+|w_\ell\partial_{\beta_1}f|_{L_v^2}|w_0\partial_{\beta_2}g|_{N^s_\gamma}\big\} |w_\ell\partial_\beta h|_{N^s_\gamma}\\
&+\sum_{\beta_1+\beta_2\leq\beta}\min\bigg\{\sum_{|\beta'|\leq2}|\langle v\rangle^{-m}w_0\partial_{\beta_1+\beta'}f|_{L_v^2}|w_\ell\partial_{\beta_2}g|_{N^s_\gamma},\\
&\qquad \qquad \qquad \quad |\langle v\rangle^{-m}w_0 \partial_{\beta_1}f|_{L_v^2}\sum_{|\beta'|\leq2}|w_\ell\partial_{\beta_2+\beta'}g|_{N^s_\gamma}\bigg\} |w_\ell\partial_\beta h|_{N^s_\gamma},
\end{align*}
where $|\cdot|_{N^s_\gamma}$ denotes the anisotropic norm defined in \cite[(1.8), p.~774]{GS2011}, and the implicit constant is independent of $t$.
\end{remark}
\smallskip

\subsection{Mixture Estimates of the Nonlinear Collision Term}
\hspace*{\fill}

We now apply Lemma \ref{lem:Gamma-weighted} to estimate the nonlinear collision term $\Gamma(f^\varepsilon,f^\varepsilon)$.
Recall the definitions of $w_l(\alpha,\beta)$, $\overline w_l(\alpha,\beta)$, $\widetilde w_l(\alpha,\beta)$, $\widetilde{\mathcal E}_{N,l}(t)$, and $\widetilde{\mathcal D}_{N,l}(t)$ in \eqref{weight function}, \eqref{weight function2}, \eqref{weight function1}, \eqref{energy functional} and \eqref{dissipation functional}, respectively.
In the proof of the following lemma, we use the decomposition
\begin{align}\label{gamma decomposition}
\Gamma(f^\varepsilon,f^\varepsilon)
=\Gamma(\mathbf Pf^\varepsilon,\mathbf Pf^\varepsilon)
+\Gamma(\mathbf Pf^\varepsilon,\mathbf P^\perp f^\varepsilon)
+\Gamma(\mathbf P^\perp f^\varepsilon,\mathbf Pf^\varepsilon)
+\Gamma(\mathbf P^\perp f^\varepsilon,\mathbf P^\perp f^\varepsilon).
\end{align}

\begin{lemma}\label{soft Gamma}
Let $(f^\varepsilon,E^\varepsilon,B^\varepsilon)$ be the solution to \eqref{rVMB}, and assume that all the parameters involved are fixed and satisfy \eqref{mainth1 assumption}.
\begin{itemize}
\item[(1)] For $|\alpha|\le N$, there holds
\begin{align}\label{soft gamma3}
 \frac1\varepsilon \Big(\partial^\alpha\Gamma(f^\varepsilon,f^\varepsilon), \partial^\alpha f^\varepsilon\Big)
\lesssim \sqrt{\mathcal E_{N}(t)} \mathcal D_N(t).
\end{align}

\item[(2)] For $|\alpha|+|\beta|\le N-1$, there holds
\begin{align}\label{soft gamma4}
\frac1\varepsilon\Big(\partial_\beta^\alpha\Gamma(f^\varepsilon,f^\varepsilon),
\widetilde w_l^2(\alpha,\beta)\partial_\beta^\alpha\mathbf P^\perp f^\varepsilon\Big)
\lesssim \sqrt{\widetilde{\mathcal E}_{N,l}(t)}\widetilde{\mathcal D}_{N,l}(t).
\end{align}

\item[(3)] For $|\alpha|=N$, there holds
\begin{align}\label{soft gamma5}
\Big(\partial^\alpha\Gamma(f^\varepsilon,f^\varepsilon),\widetilde w_l^2(\alpha,0)\partial^\alpha f^\varepsilon\Big)
\lesssim \sqrt{\widetilde{\mathcal E}_{N,l}(t)}\widetilde{\mathcal D}_{N,l}(t).
\end{align}

\item[(4)] For $|\alpha|+|\beta|=N$, $|\beta|\ge1$ and $|\alpha|\leq N-1$, there holds
\begin{align}\label{soft gamma6}
\frac1\varepsilon\Big(\partial_\beta^\alpha\Gamma(f^\varepsilon,f^\varepsilon),
\widetilde w_l^2(\alpha,\beta)\partial_\beta^\alpha\mathbf P^\perp f^\varepsilon\Big)
\lesssim \sqrt{\widetilde{\mathcal E}_{N,l}(t)}\widetilde{\mathcal D}_{N,l}(t).
\end{align}
\end{itemize}
\end{lemma}

\begin{proof}
For the sake of brevity, we present only the proofs of \eqref{soft gamma3} and \eqref{soft gamma6}.
The proofs of \eqref{soft gamma4} and \eqref{soft gamma5} follow by analogous arguments.

We first prove \eqref{soft gamma6}.
By the decomposition \eqref{gamma decomposition}, we write
\begin{align}\label{Gamma-mixed-decomposition}
\frac1\varepsilon\Big(\partial_\beta^\alpha\Gamma(f^\varepsilon,f^\varepsilon),
\widetilde w_l^2(\alpha,\beta)\partial_\beta^\alpha\mathbf P^\perp f^\varepsilon\Big)
:=\sum_{i=1}^{4}I_i,
\end{align}
where the terms $I_{i}$ $(i=1,\cdots,4)$ correspond to the four terms on the right-hand side of \eqref{gamma decomposition}, respectively.

For the term $I_{1}$, by \eqref{Pf define}, Lemma \ref{lem:Gamma-weighted}, H\"older's inequality, and \eqref{Sobolev-ineq}, we obtain
\begin{align*}
I_{1}
\lesssim\;& \frac1\varepsilon \sum_{\alpha_1\leq\alpha}
\int_{\mathbb R^3} |\partial^{\alpha_1}\mathbf Pf^\varepsilon|_{L^2}
|\partial^{\alpha-\alpha_1}\mathbf Pf^\varepsilon|_{L^2}
|\widetilde w_l(\alpha,\beta)\partial_\beta^\alpha\mathbf P^\perp f^\varepsilon|_D\d x\\
\lesssim\;& \frac1\varepsilon\sum_{\alpha_1\leq\alpha}
\|\partial^{\alpha_1}\mathbf Pf^\varepsilon\|_{L_x^6L_v^2}
\|\partial^{\alpha-\alpha_1}\mathbf Pf^\varepsilon\|_{L_x^3L_v^2}
\|\widetilde w_l(\alpha,\beta)\partial_\beta^\alpha\mathbf P^\perp f^\varepsilon\|_D\\
\lesssim\;& \sqrt{\mathcal E_N(t)}\widetilde{\mathcal D}_{N,l}(t).
\end{align*}
Here we have used the inequality $\widetilde w_l(\alpha,\beta)\mu^\delta\lesssim\mu^{\frac{\delta}{2}}$ provided that $q>0$ is sufficiently small, together with the fact that $1\leq|\alpha_1+\alpha'|\leq N$ for $|\alpha'|=1$, since $|\alpha_1|\leq|\alpha|\leq N-1$.
Here and below, $\alpha'$ denotes the additional spatial multi-index introduced by Sobolev embeddings \eqref{Sobolev-ineq}.

For the term $I_{2}$, taking the first term in the minimum in \eqref{Gamma-weighted-2} and applying H\"{o}lder's inequality, \eqref{Pf define} and \eqref{norm inequality}, we obtain
\begin{align*}
I_{2}
\lesssim\;&
\frac1\varepsilon\sum_{\substack{\alpha_1\leq\alpha\\\beta_2\leq\beta}}
\int_{\mathbb R^3}|\partial^{\alpha_1}\mathbf Pf^\varepsilon|_{L^2}
\left|(1+t)^{-\frac{1+\varepsilon_0}{4}}w_l(\alpha,\beta)\partial_{\beta_2}^{\alpha-\alpha_1}\mathbf P^\perp f^\varepsilon\right|_D
|\widetilde w_l(\alpha,\beta)\partial_\beta^\alpha\mathbf P^\perp f^\varepsilon|_D\d x\\
\lesssim\;&
\frac1\varepsilon \sum_{\beta_2\le\beta}
\|\mathbf Pf^\varepsilon\|_{L_x^\infty L_v^2}
\left\|(1+t)^{-\frac{1+\varepsilon_0}{4}}w_l(\alpha,\beta)\partial_{\beta_2}^{\alpha}\mathbf P^\perp f^\varepsilon\right\|_D
\|\widetilde w_l(\alpha,\beta)\partial_\beta^\alpha\mathbf P^\perp f^\varepsilon\|_D \\
&+ \frac1\varepsilon
\sum_{\substack{1\leq|\alpha_1|\leq |\alpha| \\ \beta_2\le\beta}}
\|\partial^{\alpha_1}\mathbf Pf^\varepsilon\|_{L_x^6L_v^2}
\left\|(1+t)^{-\frac{1+\varepsilon_0}{4}}w_l(\alpha,\beta)\partial_{\beta_2}^{\alpha-\alpha_1}\mathbf P^\perp f^\varepsilon\right\|_{L_x^3L_D^2}
\|\widetilde w_l(\alpha,\beta)\partial_\beta^\alpha\mathbf P^\perp f^\varepsilon\|_D \\
\lesssim\;& \sqrt{\mathcal E_N(t)} \widetilde{\mathcal D}_{N,l}(t).
\end{align*}
Here, we have used
\begin{equation*}
\begin{aligned}
(1+t)^{-\frac{1+\varepsilon_0}{4}} w_l(\alpha,\beta)
&\leq \widetilde w_l(\alpha,\beta_2),
&\quad &\text{for}\; |\beta_2|\leq |\beta|,\\
(1+t)^{-\frac{1+\varepsilon_0}{4}} w_l(\alpha,\beta)
&\leq \widetilde w_l(\alpha-\alpha_1+\alpha',\beta_2),
&\quad &\text{for}\;\; |\alpha_1|\geq 1, \; |\alpha'|\leq 1 \;\text{and}\; |\beta_2|\leq |\beta|.
\end{aligned}
\end{equation*}

Similarly, for the term $I_{3}$, choosing the second bound in the minimum of \eqref{Gamma-weighted-2} gives
\begin{align*}
I_{3}\lesssim\sqrt{\mathcal E_N(t)}\widetilde{\mathcal D}_{N,l}(t).
\end{align*}

For the term $I_{4}$, since $\gamma+2s<0$, applying Lemma \ref{lem:Gamma-weighted} in conjunction with \eqref{norm inequality} gives
\begin{align*}
I_{4}
\lesssim\;&\frac{1}{\varepsilon} \sum \int_{\mathbb R^3} (1+t)^{-\frac{1+\varepsilon_0}{4}}
\left|e^{\frac{q \langle v \rangle^2}{(1+t)^\vartheta}}\partial_{\beta_1}^{\alpha_1}\mathbf P^\perp f^\varepsilon\right|_{L^2}
|w_l(\alpha,\beta)\partial_{\beta_2}^{\alpha_2}\mathbf P^\perp f^\varepsilon|_D
|\widetilde w_l(\alpha,\beta)\partial_\beta^\alpha\mathbf P^\perp f^\varepsilon|_D\d x \\
&+\frac{1}{\varepsilon} \sum \int_{\mathbb R^3} (1+t)^{-\frac{1+\varepsilon_0}{4}}
\left|e^{\frac{q \langle v \rangle^2}{(1+t)^\vartheta}}\partial_{\beta_1}^{\alpha_1}\mathbf P^\perp f^\varepsilon\right|_D
|w_l(\alpha,\beta)\partial_{\beta_2}^{\alpha_2}\mathbf P^\perp f^\varepsilon|_{L^2}
|\widetilde w_l(\alpha,\beta)\partial_\beta^\alpha\mathbf P^\perp f^\varepsilon|_D\d x\\
:=\;&I_{4,1}+I_{4,2},
\end{align*}
where the summation $\sum$ is taken over $\alpha_1+\alpha_2=\alpha$ and $\beta_1+\beta_2 \leq \beta$.
We estimate $I_{4,1}$ by considering two cases. If $|\alpha_2|+|\beta_2|\leq N/2$, then $|\alpha_2+\alpha'|+|\beta_2|\leq N-1$ for $1\leq|\alpha'|\leq2$.
Hence, we use H\"{o}lder's inequality, \eqref{Sobolev-ineq}, \eqref{minkowski} and
\begin{equation}
\begin{aligned}
&(1+t)^{-\frac{1+\varepsilon_0}{4}} e^{\frac{q \langle v \rangle^2}{(1+t)^\vartheta}} \leq \widetilde{w}_l( \overline{\alpha}, \overline{\beta}),
&\quad& \text{ for }  |\overline{\alpha}|+|\overline{\beta}|\leq N, \\
&w_l(\alpha, \beta) \leq \widetilde{w}_l( \overline{\alpha}, \overline{\beta}),
&\quad& \text{ for }\; | \overline{\alpha} | + |\overline{\beta}| \leq | \alpha |+ |\beta|-1 \text{ and } | \overline{\beta} | \leq |\beta|,\label{weight inequality2}
\end{aligned}
\end{equation}
to obtain
\begin{align*}
I_{4,1} \lesssim \frac{1}{\varepsilon}
\left\| \widetilde{w}_{l}(\alpha_1,\beta_1) \partial^{\alpha_1}_{\beta_1} \mathbf{P}^{\perp}f^{\varepsilon} \right\|
\left\| w_l(\alpha, \beta) \partial^{\alpha_2}_{\beta_2} \mathbf{P}^{\perp} f^{\varepsilon}\right\|_{L^\infty_x L^2_D}
\left\| \widetilde{w}_{l}(\alpha,\beta) \partial^{\alpha}_{\beta} \mathbf{P}^{\perp} f^{\varepsilon}\right\|_{D}
\lesssim \sqrt{\widetilde{\mathcal{E}}_{N,l}(t)} \widetilde{\mathcal{D}}_{N,l}(t).
\end{align*}
Here, we assign the time factor $(1+t)^{-\frac{1+\varepsilon_0}{4}}$ to the first factor.
In fact, after applying \eqref{Sobolev-ineq}, we always attach the factor $(1+t)^{-\frac{1+\varepsilon_0}{4}}$ to whichever of the first two factors has total derivative order $N$.
If neither factor has total derivative order $N$, we bound the time factor by $1$.
For the case $N/2 < |\alpha_2|+|\beta_2| \leq N$, observing that $|\alpha_1+ \alpha^\prime|+|\beta_1| \leq N-1$ for $1 \leq |\alpha^\prime| \leq 2$, by H\"{o}lder's inequality, one has
\begin{align*}
I_{4,1}\lesssim \frac{1}{\varepsilon}
\left\| e^{\frac{q \langle v \rangle^2}{(1+t)^\vartheta}} \partial^{\alpha_1}_{\beta_1} \mathbf{P}^{\perp}f^{\varepsilon} \right\|_{L^\infty_x L^2_v}
\left\| \widetilde{w}_l(\alpha,\beta) \partial^{\alpha_2}_{\beta_2} \mathbf{P}^{\perp} f^{\varepsilon}\right\|_{D}
\left\| \widetilde{w}_{l}(\alpha,\beta) \partial^{\alpha}_{\beta} \mathbf{P}^{\perp} f^{\varepsilon}\right\|_{D},
\end{align*}
which is further bounded by $\sqrt{\widetilde{\mathcal{E}}_{N,l}(t)} \widetilde{\mathcal{D}}_{N,l}(t)$ in view of \eqref{Sobolev-ineq}, \eqref{minkowski}, \eqref{weight inequality2} and
\begin{align*}
\widetilde{w}_l(\alpha, \beta) \leq \widetilde{w}_l( \overline{\alpha}, \overline{\beta}) \quad
\text{ for } | \overline{\alpha} | + |\overline{\beta}| \leq | \alpha |+ |\beta|\; \text{ and }| \overline{\beta} | \leq |\beta|.
\end{align*}
An analogous argument yields the same bound for the term $I_{4,2}$.

Therefore, we obtain the desired estimate \eqref{soft gamma6} by substituting all the above estimates into \eqref{Gamma-mixed-decomposition}.

We next prove \eqref{soft gamma3}. By the collision invariants, we have
\begin{align*}
\Big(\partial^\alpha\Gamma(f^\varepsilon,f^\varepsilon),\partial^\alpha f^\varepsilon\Big)
=\Big(\partial^\alpha\Gamma(f^\varepsilon,f^\varepsilon),\partial^\alpha\mathbf P^\perp f^\varepsilon\Big).
\end{align*}
Applying \eqref{Gamma-unweighted}, H\"older's inequality, \eqref{Sobolev-ineq}, and \eqref{norm inequality}, we obtain
\begin{align*}
\frac1\varepsilon
\left|\Big(\partial^\alpha\Gamma(f^\varepsilon,f^\varepsilon),\partial^\alpha f^\varepsilon\Big)\right|
\lesssim\;&\frac1\varepsilon
\bigg\{\|f^\varepsilon\|_{L_x^\infty L_D^2}\|\partial^\alpha f^\varepsilon\|
+\chi_{|\alpha|\geq1}\|f^\varepsilon\|_{L_x^\infty L_v^2}\|\partial^\alpha f^\varepsilon\|_D\\
&\quad+\sum_{1\leq|\alpha_1|\leq|\alpha|-1}\|\partial^{\alpha_1}f^\varepsilon\|_{L_x^6L_D^2}
\|\partial^{\alpha-\alpha_1}f^\varepsilon\|_{L_x^3L_v^2}\bigg\}
\|\partial^\alpha\mathbf P^\perp f^\varepsilon\|_D\\
\lesssim\;& \sqrt{\mathcal E_N(t)}\mathcal D_N(t).
\end{align*}
This proves \eqref{soft gamma3}.
\end{proof}
\medskip

\section{A Priori Estimates}\label{The a Priori Estimate}

In this section, we establish a priori energy estimates for the one-species VMB system \eqref{rVMB} in the soft potential case, uniformly in $\varepsilon\in(0,1]$.

Let $0<s<1$, $\max\{-3,-\frac32-2s\}<\gamma<-2s$, and $0<\varepsilon\leq1$.
All the remaining parameters are fixed to satisfy \eqref{mainth1 assumption}.
Suppose that the solution $(f^\varepsilon,E^\varepsilon,B^\varepsilon)$ to the VMB system \eqref{rVMB} on $[0,T]$, where $T>0$, satisfies
\begin{align}\label{priori-assumptions-decay}
\sup_{0\leq t\leq T}\Big\{\overline{\mathcal E}_{N,l}(t)+(1+t)^{\frac{3+2\varrho}{4}}\|E^\varepsilon(t)\|_{L^\infty}^2
+(1+t)^{\frac{2+\varrho}{2}}\|\nabla_x(E^\varepsilon,B^\varepsilon)(t)\|_{L^\infty}^2\Big\}\leq\delta_0^2,
\end{align}
where $\delta_0>0$ is sufficiently small.
\smallskip

\subsection{Macroscopic Estimate}
\hspace*{\fill}

In this subsection, we derive macroscopic energy estimates for the one-species VMB system \eqref{rVMB} that are uniform in $\varepsilon\in (0,1]$.

We first define the macroscopic dissipation by
\begin{align}\label{macro:dissipation}
\mathcal D_{mac}^{N}(t):=
\sum_{1\leq|\alpha|\leq N}\|\partial^\alpha\mathbf Pf^\varepsilon\|^2
+\varepsilon^2\|a^\varepsilon\|^2
+\sum_{1\leq|\alpha|\leq N-1}\|\partial^\alpha E^\varepsilon\|^2
+\sum_{2\leq|\alpha|\leq N-1}\|\partial^\alpha B^\varepsilon\|^2.
\end{align}

\begin{lemma}\label{macroscopic estimate}
Let $(f^\varepsilon,E^\varepsilon,B^\varepsilon)$ be the solution to the VMB system \eqref{rVMB}. There exists an interactive energy functional $\mathcal E_{int}^{N}(t)$ defined in \eqref{macro:int-definition} satisfying
\begin{align*}
|\mathcal E_{int}^{N}(t)|
\lesssim \sum_{|\alpha|\leq N}\left(\|\partial^\alpha f^\varepsilon\|^2
+\|\partial^\alpha E^\varepsilon\|^2 +\|\partial^\alpha B^\varepsilon\|^2\right),
\end{align*}
such that for any $t\geq 0$, there holds
\begin{align}\label{macro estimate}
\frac{\d}{\d t}\mathcal E_{int}^{N}(t)+\mathcal D_{mac}^{N}(t)
\lesssim\;&
\frac1{\varepsilon^2}\sum_{|\alpha|\leq N}\|\partial_x^\alpha\mathbf P^\perp f^\varepsilon\|_D^2+\mathcal E_N(t)\mathcal D_N(t).
\end{align}
\end{lemma}

\begin{proof}
Following \cite{LZ2016}, we apply the macro-micro decomposition \eqref{f decomposition} introduced in \cite{Guo2003}, and define the moment functions
\begin{align*}
\Theta_{ij}(f):=\int_{\mathbb{R}^3}(v_i v_j-1)\mu^{1/2}f\d v,
\qquad
\Lambda_i(f):=\frac{1}{10}\int_{\mathbb{R}^3}(|v|^2-5)v_i\mu^{1/2}f\d v.
\end{align*}
Then we can derive from \eqref{rVMB} the following two fluid-type systems of equations
\begin{equation}\label{macro equation 1}
\left\{\begin{aligned}
&\partial_t a^\varepsilon+\frac{1}{\varepsilon}\nabla_x\cdot b^\varepsilon=0,\\
&\partial_t b_i^\varepsilon+\frac{1}{\varepsilon}\partial_i(a^\varepsilon+2c^\varepsilon)-E_i^\varepsilon
+\frac{1}{\varepsilon}\sum_{j=1}^3\partial_j\Theta_{ij}(\mathbf{P}^{\perp}f^\varepsilon)
=\varepsilon a^\varepsilon E_i^\varepsilon+(b^\varepsilon\times B^\varepsilon)_i,\\
&\partial_t c^\varepsilon+\frac{1}{3\varepsilon}\nabla_x\cdot b^\varepsilon
+\frac{5}{3\varepsilon}\sum_{i=1}^3\partial_i\Lambda_i(\mathbf{P}^{\perp}f^\varepsilon)
=\frac{\varepsilon}{3}E^\varepsilon\cdot b^\varepsilon,
\end{aligned}\right.
\end{equation}
and
\begin{equation}\label{macro equation 2}
\left\{\begin{aligned}
&\partial_t\Big(\Theta_{ij}(\mathbf{P}^{\perp}f^\varepsilon)+2c^\varepsilon \delta_{ij}\Big)
+\frac{1}{\varepsilon}\partial_i b_j^\varepsilon
+\frac{1}{\varepsilon}\partial_j b_i^\varepsilon
=\Theta_{ij}(g^\varepsilon+h^\varepsilon),\\
&\partial_t\Lambda_i(\mathbf{P}^{\perp}f^\varepsilon)
+\frac{1}{\varepsilon}\partial_i c^\varepsilon
=\Lambda_i(g^\varepsilon+h^\varepsilon).
\end{aligned}\right.
\end{equation}
Here $1\leq i,j\leq3$ and $\delta_{ij}$ denotes the Kronecker delta.
The source terms $g^\varepsilon$ and $h^\varepsilon$ are defined by
\begin{align}
\begin{split}\label{g,h define}
g^\varepsilon:=\;&\frac{\varepsilon}{2}E^\varepsilon\cdot vf^\varepsilon
-\big(\varepsilon E^\varepsilon+v\times B^\varepsilon\big)\cdot\nabla_v f^\varepsilon
+\frac{1}{\varepsilon}\Gamma(f^\varepsilon,f^\varepsilon),\\
h^\varepsilon:=\;&-\frac{1}{\varepsilon}v\cdot\nabla_x\mathbf{P}^{\perp}f^\varepsilon
-\frac{1}{\varepsilon^2}Lf^\varepsilon.
\end{split}
\end{align}

Next, applying an argument analogous to that in \cite{LZ2016} to the fluid-type systems \eqref{macro equation 1} and \eqref{macro equation 2}, we obtain, for any $\eta>0$,
\begin{align}\label{macro:abc-estimate}
&\varepsilon\frac{\d}{\d t}
\mathcal{E}_{a^\varepsilon,b^\varepsilon,c^\varepsilon}^{N}(t)
+ \sum_{|\alpha|\leq N-1}\left\|\partial^\alpha\nabla_x(a^\varepsilon,b^\varepsilon,c^\varepsilon)\right\|^2
+\varepsilon^2\|a^\varepsilon\|^2 \nonumber\\
\lesssim\;&
\eta\varepsilon^2 \sum_{1\leq|\alpha|\leq N-1}\|\partial^\alpha E^\varepsilon\|^2
+\frac{1}{\varepsilon^2}\sum_{|\alpha|\leq N}\left\|\partial^\alpha\mathbf{P}^{\perp}f^\varepsilon\right\|_D^2
+\mathcal{E}_N(t)\mathcal{D}_N(t).
\end{align}
Here, $\mathcal{E}_{a^\varepsilon,b^\varepsilon,c^\varepsilon}^{N}(t)$ is defined by
\begin{align}\label{macro:abc-functional}
\mathcal{E}_{a^\varepsilon,b^\varepsilon,c^\varepsilon}^{N}(t)
:=\;& \overline C_1 \sum_{|\alpha|\leq N-1}
\Bigg\{\sum_{i,j=1}^3
\Big(\partial^\alpha\Theta_{ij}(\mathbf{P}^{\perp}f^\varepsilon),
\partial^\alpha(\partial_i b_j^\varepsilon+\partial_j b_i^\varepsilon)\Big) \nonumber\\
&\qquad \qquad
+\sum_{i=1}^3 \Big(\partial^\alpha\Lambda_i(\mathbf{P}^{\perp}f^\varepsilon),
\partial^\alpha\partial_i c^\varepsilon \Big)
+\kappa_1 \Big(\partial^\alpha b^\varepsilon,
\partial^\alpha\nabla_xa^\varepsilon\Big) \Bigg\},
\end{align}
where $0<\kappa_1\ll1$ is a fixed constant and $\overline C_1>0$ is sufficiently large.
For brevity, we omit the details of the derivation.

To estimate $E^\varepsilon$, we use the Maxwell system in \eqref{rVMB} and \eqref{macro equation 1} to obtain, for $1\leq|\alpha|\leq N-1$,
\begin{align}\label{macro:E-identity}
&\|\partial^\alpha E^\varepsilon\|^2=\Big(\partial^\alpha E^\varepsilon,\partial^\alpha E^\varepsilon\Big)
\nonumber\\
=\;& \Big(\partial^\alpha\partial_t b^\varepsilon
+\frac1\varepsilon\partial^\alpha\nabla_x(a^\varepsilon+2c^\varepsilon)
+\frac1\varepsilon\partial^\alpha\nabla_x\cdot \Theta(\mathbf P^\perp f^\varepsilon)
-\partial^\alpha (\varepsilon a^\varepsilon E^\varepsilon
+b^\varepsilon\times B^\varepsilon), \partial^\alpha E^\varepsilon \Big)\nonumber\\
=\;& \frac{\d}{\d t}
\Big(\partial^\alpha b^\varepsilon, \partial^\alpha E^\varepsilon\Big)
-\Big(\partial^\alpha b^\varepsilon, \partial^\alpha\partial_t E^\varepsilon\Big)\nonumber\\
&+ \Big(\frac1\varepsilon\partial^\alpha\nabla_x(a^\varepsilon+2c^\varepsilon)
+\frac1\varepsilon\partial^\alpha\nabla_x\cdot \Theta(\mathbf P^\perp f^\varepsilon)
-\partial^\alpha (\varepsilon a^\varepsilon E^\varepsilon +b^\varepsilon\times B^\varepsilon),
\partial^\alpha E^\varepsilon\Big) \nonumber\\
=\;& \frac{\d}{\d t}
\Big(\partial^\alpha b^\varepsilon, \partial^\alpha E^\varepsilon\Big)
+\|\partial^\alpha b^\varepsilon\|^2-\Big(\partial^\alpha b^\varepsilon, \partial^\alpha\nabla_x\times B^\varepsilon\Big)
-\|\partial^\alpha a^\varepsilon\|^2 -2\Big(\partial^\alpha a^\varepsilon, \partial^\alpha c^\varepsilon\Big)\nonumber\\
&+\frac1\varepsilon
\Big( \partial^\alpha\nabla_x\cdot \Theta(\mathbf P^\perp f^\varepsilon), \partial^\alpha E^\varepsilon \Big)
-\Big( \partial^\alpha (\varepsilon a^\varepsilon E^\varepsilon +b^\varepsilon\times B^\varepsilon), \partial^\alpha E^\varepsilon \Big).
\end{align}
To estimate the third term on the right-hand side of \eqref{macro:E-identity}, we first consider $|\alpha|=1$. By the Cauchy--Schwarz inequality, we have
\begin{align*}
\left|\Big(\partial^\alpha\nabla_x\times B^\varepsilon, \partial^\alpha b^\varepsilon\Big)\right|
 \lesssim \eta\|\partial^\alpha\nabla_x B^\varepsilon\|^2 +\|\partial^\alpha b^\varepsilon\|^2.
\end{align*}
For $2\leq|\alpha|\leq N-1$, integration by parts yields
\begin{align*}
\left|\Big(\partial^\alpha\nabla_x\times B^\varepsilon, \partial^\alpha b^\varepsilon\Big)\right|
=\left|\Big(\partial^\alpha B^\varepsilon, \partial^\alpha\nabla_x\times b^\varepsilon\Big)\right|
\lesssim \eta\|\partial^\alpha B^\varepsilon\|^2 +\|\partial^\alpha\nabla_x b^\varepsilon\|^2.
\end{align*}
Combining these two cases and summing over $1\leq|\alpha|\leq N-1$, we obtain
\begin{align}
\sum_{1\leq|\alpha|\leq N-1} \left|\Big(\partial^\alpha\nabla_x\times B^\varepsilon, \partial^\alpha b^\varepsilon\Big)\right|
\lesssim \eta \sum_{2\leq|\alpha|\leq N-1} \|\partial^\alpha B^\varepsilon\|^2
+\sum_{1\leq |\alpha|\leq N}\|\partial^\alpha b^\varepsilon\|^2.
\end{align}
The last three terms on the right-hand side of \eqref{macro:E-identity} can be controlled by
\begin{align*}
\eta \Big(\|\partial^\alpha a^\varepsilon\|^2+\|\partial^\alpha E^\varepsilon\|^2\Big)
+\|\partial^\alpha c^\varepsilon\|^2
+\frac{1}{\varepsilon^2} \|\partial^\alpha \nabla_x \mathbf P^\perp f^\varepsilon\|_D^2
+\mathcal E_N(t)\mathcal D_N(t).
\end{align*}
Summing \eqref{macro:E-identity} over $1\leq|\alpha|\leq N-1$, applying the above estimates, and choosing $\eta>0$ sufficiently small, we obtain
\begin{align}\label{macro:E-estimate}
&\frac{\d}{\d t}\mathcal E_{E^\varepsilon}^{N}(t)
+\frac12\sum_{1\leq|\alpha|\leq N-1}\Big(\|\partial^\alpha E^\varepsilon\|^2+\|\partial^\alpha a^\varepsilon\|^2\Big)\\
\lesssim\;&\eta\sum_{2\leq|\alpha|\leq N-1}\|\partial^\alpha B^\varepsilon\|^2
+\sum_{1\leq|\alpha|\leq N}\|\partial^\alpha(b^\varepsilon,c^\varepsilon)\|^2+\frac{1}{\varepsilon^2}\sum_{2\leq |\alpha|\leq N} \|\partial^\alpha \mathbf P^\perp f^\varepsilon\|_D^2 +\mathcal E_N(t)\mathcal D_N(t), \nonumber
\end{align}
where $\mathcal E_{E^\varepsilon}^{N}(t)$ is given by
\begin{align*}
\mathcal E_{E^\varepsilon}^{N}(t):=-\sum_{1\leq|\alpha|\leq N-1} \Big( \partial^\alpha E^\varepsilon,\partial^\alpha b^\varepsilon\Big).
\end{align*}

For the dissipation of $B^\varepsilon$, we use the Maxwell system in \eqref{rVMB} to obtain, for $1\leq|\alpha|\leq N-2$,
\begin{align}\label{macro:B-identity}
&\|\partial^\alpha\nabla_x B^\varepsilon\|^2
=\|\partial^\alpha\nabla_x\times B^\varepsilon\|^2
\nonumber\\
=\;& \Big(\partial^\alpha\partial_t E^\varepsilon
+\partial^\alpha b^\varepsilon, \partial^\alpha\nabla_x\times B^\varepsilon\Big) \nonumber\\
=\;& \frac{\d}{\d t}
\Big(\partial^\alpha E^\varepsilon, \partial^\alpha\nabla_x\times B^\varepsilon\Big)
-\Big(\partial^\alpha E^\varepsilon, \partial^\alpha\nabla_x\times\partial_t B^\varepsilon\Big)
+\Big(\partial^\alpha b^\varepsilon, \partial^\alpha\nabla_x\times B^\varepsilon\Big) \nonumber\\
=\;& \frac{\d}{\d t}
\Big(\partial^\alpha E^\varepsilon, \partial^\alpha\nabla_x\times B^\varepsilon\Big)
+\Big(\partial^\alpha E^\varepsilon, \partial^\alpha\nabla_x\times(\nabla_x\times E^\varepsilon)\Big)
+\Big(\partial^\alpha b^\varepsilon, \partial^\alpha\nabla_x\times B^\varepsilon\Big) \nonumber\\
=\;& \frac{\d}{\d t}
\Big(\partial^\alpha E^\varepsilon, \partial^\alpha\nabla_x\times B^\varepsilon\Big)
+\|\partial^\alpha\nabla_x\times E^\varepsilon\|^2
+\Big(\partial^\alpha b^\varepsilon, \partial^\alpha\nabla_x\times B^\varepsilon\Big).
\end{align}
The last term can be estimated by the Cauchy--Schwarz inequality as follows
\begin{align*}
\left|\Big(\partial^\alpha b^\varepsilon, \partial^\alpha\nabla_x\times B^\varepsilon\Big)\right|
\lesssim \eta \|\partial^\alpha\nabla_x B^\varepsilon\|^2+\|\partial^\alpha b^\varepsilon\|^2.
\end{align*}
Summing \eqref{macro:B-identity} over $1\leq|\alpha|\leq N-2$, applying the above estimate and choosing $\eta>0$ sufficiently small, we have
\begin{align}\label{macro:B-estimate}
\frac{\d}{\d t}\mathcal E_{B^\varepsilon}^{N}(t) +\frac12\sum_{2\leq|\alpha|\leq N-1} \|\partial^\alpha B^\varepsilon\|^2
\lesssim \sum_{2\leq|\alpha|\leq N-1} \|\partial^\alpha E^\varepsilon\|^2
+\sum_{1\leq|\alpha|\leq N-2} \|\partial^\alpha b^\varepsilon\|^2,
\end{align}
where $\mathcal E_{B^\varepsilon}^{N}(t)$ is given by
\begin{align*}
\mathcal E_{B^\varepsilon}^{N}(t)
:=-\sum_{1\leq|\alpha|\leq N-2} \Big(\partial^\alpha E^\varepsilon, \partial^\alpha\nabla_x\times B^\varepsilon\Big).
\end{align*}

Finally, we define
\begin{align}\label{macro:int-definition}
\mathcal E_{int}^{N}(t)
:=\overline C_2\Big( \varepsilon\mathcal E_{a^\varepsilon,b^\varepsilon,c^\varepsilon}^{N}(t)
+\kappa_2\mathcal E_{E^\varepsilon}^{N}(t)
+\kappa_3\mathcal E_{B^\varepsilon}^{N}(t)\Big),
\end{align}
where $0<\kappa_3\ll\kappa_2\ll 1$ are suitably chosen fixed constants and $\overline C_2>0$ is sufficiently large.
Taking $\eta>0$ sufficiently small in \eqref{macro:abc-estimate} and \eqref{macro:E-estimate} and combining these estimates with \eqref{macro:B-estimate}, we arrive at \eqref{macro estimate}.
The asserted bound for $\mathcal E_{int}^{N}(t)$ follows from the Cauchy--Schwarz inequality.
This completes the proof.
\end{proof}

\subsection{Energy Estimate Without Weight}\label{ss:unweighted-subsection}

\begin{lemma}\label{ss:unweighted}
Let $(f^\varepsilon,E^\varepsilon,B^\varepsilon)$ be the solution to the VMB system \eqref{rVMB}, and assume that all the parameters involved are fixed and satisfy \eqref{mainth1 assumption}.
Then, there holds
\begin{align}\label{ss:unweighted-estimate}
&\frac{1}{2}\frac{\d}{\d t}
\Bigg\{\sum_{|\alpha|\leq N} \|\partial^\alpha (f^\varepsilon, E^\varepsilon,B^\varepsilon)\|^2
-\varepsilon\int_{\mathbb R^3}(a^\varepsilon+2c^\varepsilon)|b^\varepsilon|^2\d x\Bigg\}
+\frac{\lambda}{\varepsilon^2}\sum_{|\alpha|\leq N}\|\partial^\alpha\mathbf P^\perp f^\varepsilon\|_D^2\nonumber\\
\lesssim\;&
\Big\{\sqrt{\widetilde{\mathcal E}_{N,l}(t)}+\mathcal E_N(t)+\delta_0\Big\}\widetilde{\mathcal D}_{N,l}(t).
\end{align}
\end{lemma}

\begin{proof}
Applying $\partial^\alpha$ with $|\alpha|\leq N$ to the first equation in \eqref{rVMB}, taking the inner product with $\partial^\alpha f^\varepsilon$, and using the Maxwell equations in \eqref{rVMB},
integration by parts and \eqref{L-coercive-unweighted}, we obtain
\begin{align}\label{unweighted-energy-basic}
&\frac12\frac{\d}{\d t}
\left\|\partial^\alpha(f^\varepsilon,E^\varepsilon,B^\varepsilon)\right\|^2
+\frac{\lambda}{\varepsilon^2}\left\|\partial^\alpha\mathbf P^\perp f^\varepsilon\right\|_D^2\nonumber\\
\leq\;&\frac{\varepsilon}{2}\Big(\partial^\alpha(E^\varepsilon\cdot v f^\varepsilon),\partial^\alpha f^\varepsilon\Big)
-\chi_{|\alpha|\geq1}\Big(\partial^\alpha\big((v\times B^\varepsilon)\cdot\nabla_v f^\varepsilon\big),
\partial^\alpha f^\varepsilon\Big)\nonumber\\
&-\varepsilon\chi_{|\alpha|\geq1}
\Big(\partial^\alpha(E^\varepsilon\cdot\nabla_v f^\varepsilon),\partial^\alpha f^\varepsilon\Big)
+\frac1\varepsilon\Big(\partial^\alpha\Gamma(f^\varepsilon,f^\varepsilon),\partial^\alpha f^\varepsilon\Big)
:=\sum_{i=1}^{4} I_i.
\end{align}

For the term $I_1$, we distinguish the following two cases.

\noindent{\bf Case 1.} $|\alpha|=0$.
By the macro-micro decomposition \eqref{f decomposition}, we have
\begin{align*}
I_1=\;&\frac{\varepsilon}{2}
\Big(E^\varepsilon\cdot v\mathbf Pf^\varepsilon,\mathbf Pf^\varepsilon\Big)
+\varepsilon\Big(E^\varepsilon\cdot v\mathbf Pf^\varepsilon,\mathbf P^\perp f^\varepsilon\Big)
+\frac{\varepsilon}{2}\Big(E^\varepsilon\cdot v\mathbf P^\perp f^\varepsilon,\mathbf P^\perp f^\varepsilon\Big)
:=I_{1,1}+I_{1,2}+I_{1,3}.
\end{align*}
By \eqref{Pf define}, \eqref{macro equation 1}, integration by parts, and $b^\varepsilon\cdot(b^\varepsilon\times B^\varepsilon)=0$, we obtain
\begin{align*}
I_{1,1}
=\;&\varepsilon\int_{\mathbb R^3} (a^\varepsilon+2c^\varepsilon)b^\varepsilon \cdot E^\varepsilon \d x\\
=\;&\frac{\varepsilon}{2}\frac{\d}{\d t}
\int_{\mathbb R^3}(a^\varepsilon+2c^\varepsilon)|b^\varepsilon|^2 \d x
-\frac{\varepsilon}{2}\int_{\mathbb R^3}\partial_t(a^\varepsilon+2c^\varepsilon) |b^\varepsilon|^2 \d x\\
&+\int_{\mathbb R^3}(a^\varepsilon+2c^\varepsilon)b^\varepsilon\cdot \Big(\nabla_x(a^\varepsilon+2c^\varepsilon)
+\nabla_x\cdot\Theta(\mathbf P^\perp f^\varepsilon)-\varepsilon^2 a^\varepsilon E^\varepsilon\Big)\d x\\
=\;&\frac{\varepsilon}{2}\frac{\d}{\d t}
\int_{\mathbb R^3}(a^\varepsilon+2c^\varepsilon) |b^\varepsilon|^2 \d x
+\frac56\int_{\mathbb R^3}\nabla_x\cdot b^\varepsilon|b^\varepsilon|^2 \d x
-\frac12\int_{\mathbb R^3}\nabla_x\cdot b^\varepsilon(a^\varepsilon+2c^\varepsilon)^2 \d x\\
&+\frac53\int_{\mathbb R^3}\nabla_x\cdot\Lambda(\mathbf P^\perp f^\varepsilon)|b^\varepsilon|^2\d x
-\sum_{i,j=1}^3\int_{\mathbb R^3}\Theta_{ij}(\mathbf P^\perp f^\varepsilon)
\partial_j\Big((a^\varepsilon+2c^\varepsilon)b_i^\varepsilon\Big) \d x\\
&-\frac{\varepsilon^2}{3}\int_{\mathbb R^3}(E^\varepsilon\cdot b^\varepsilon)|b^\varepsilon|^2 \d x
-\varepsilon^2\int_{\mathbb R^3}a^\varepsilon(a^\varepsilon+2c^\varepsilon)E^\varepsilon\cdot b^\varepsilon \d x.
\end{align*}
By H\"older's inequality, \eqref{Sobolev-ineq} and \eqref{norm inequality}, we deduce that
\begin{align*}
I_{1,1} \leq \frac{\varepsilon}{2}\frac{\d}{\d t}
\int_{\mathbb R^3}(a^\varepsilon+2c^\varepsilon) |b^\varepsilon|^2 \d x
+ C\Big\{\sqrt{\mathcal E_N(t)}+\mathcal E_N(t)\Big\} \mathcal D_N(t).
\end{align*}
An application of H\"{o}lder's inequality, together with \eqref{Sobolev-ineq} and \eqref{norm inequality}, yields
\begin{align*}
I_{1,2}+I_{1,3}
\lesssim\;&
\varepsilon\|E^\varepsilon\|_{L^3}
\|\mathbf Pf^\varepsilon\|_{L_x^6L_v^2}
\|\mathbf P^\perp f^\varepsilon\|_D
+\varepsilon\|E^\varepsilon\|_{L^\infty}
\left\|\langle v\rangle^{1-\frac{\gamma+2s}{2}} \mathbf P^\perp f^\varepsilon\right\|
\left\|\langle v\rangle^{\frac{\gamma+2s}{2}} \mathbf P^\perp f^\varepsilon\right\|\\
\lesssim\;&
\varepsilon\|E^\varepsilon\|_{H^1}
\|\nabla_x\mathbf Pf^\varepsilon\|
\|\mathbf P^\perp f^\varepsilon\|_D
+\varepsilon\|\nabla_xE^\varepsilon\|_{H^1}
\left\|\langle v\rangle^{1-\frac{\gamma+2s}{2}} \mathbf P^\perp f^\varepsilon\right\|
\|\mathbf P^\perp f^\varepsilon\|_D\\
\lesssim\;&
\sqrt{\widetilde{\mathcal E}_{N,l}(t)} \mathcal D_N(t).
\end{align*}
Combining the above estimates, we obtain, for $|\alpha|=0$,
\begin{align}\label{I1-zero-estimate}
I_1 \leq \frac{\varepsilon}{2}\frac{\d}{\d t} \int_{\mathbb R^3}(a^\varepsilon+2c^\varepsilon)|b^\varepsilon|^2\d x
+ C\Big\{ \sqrt{\widetilde{\mathcal E}_{N,l}(t)}+\mathcal E_N(t)\Big\}\mathcal D_N(t).
\end{align}

\noindent{\bf Case 2.} $1\leq|\alpha|\leq N$.
By the macro-micro decomposition \eqref{f decomposition}, we have
\begin{align*}
I_1
=\;&\frac{\varepsilon}{2} \sum_{\alpha_1\leq\alpha} C_{\alpha}^{\alpha_1}
\Big( v\cdot\partial^{\alpha_1}E^\varepsilon \partial^{\alpha-\alpha_1}\mathbf Pf^\varepsilon,
\partial^\alpha f^\varepsilon\Big)\\
&+\frac{\varepsilon}{2} \sum_{\alpha_1\leq\alpha} C_{\alpha}^{\alpha_1}
\Big( v\cdot\partial^{\alpha_1}E^\varepsilon \partial^{\alpha-\alpha_1}\mathbf P^\perp f^\varepsilon,
\partial^\alpha\mathbf Pf^\varepsilon \Big)\\
&+\frac{\varepsilon}{2} \sum_{\alpha_1\leq\alpha} C_{\alpha}^{\alpha_1}
\Big( v\cdot\partial^{\alpha_1}E^\varepsilon \partial^{\alpha-\alpha_1}\mathbf P^\perp f^\varepsilon,
\partial^\alpha\mathbf P^\perp f^\varepsilon \Big)\\
:=\;&I_{1,1}+I_{1,2}+I_{1,3}.
\end{align*}
By \eqref{Pf define}, \eqref{norm inequality}, H\"older's inequality, and \eqref{Sobolev-ineq}, we obtain
\begin{align*}
I_{1,1}
\lesssim\;&\varepsilon \Big\{ \|E^\varepsilon\|_{L^\infty} \|\partial^\alpha\mathbf Pf^\varepsilon\|
+\|\partial^\alpha E^\varepsilon\| \|\mathbf Pf^\varepsilon\|_{L_x^\infty L_v^2}\Big\}
\|\partial^\alpha f^\varepsilon\|_D \\
&+ \varepsilon \sum_{1\leq|\alpha_1|\leq|\alpha|-1} \|\partial^{\alpha_1}E^\varepsilon\|_{L^3}
\|\partial^{\alpha-\alpha_1}\mathbf Pf^\varepsilon\|_{L_x^6L_v^2}
\|\partial^\alpha f^\varepsilon\|_D \\
\lesssim\;&
\varepsilon\|E^\varepsilon\|_{H^N}
\|\nabla_x\mathbf Pf^\varepsilon\|_{H_x^{N-1}L_v^2}
\Big\{\|\partial^\alpha\mathbf Pf^\varepsilon\|
+\|\partial^\alpha\mathbf P^\perp f^\varepsilon\|_D \Big\}\\
\lesssim\;& \sqrt{\mathcal E_N(t)}\mathcal D_N(t).
\end{align*}
By an analogous argument to that of $I_{1,1}$, we find
\begin{align*}
I_{1,2} \lesssim \sqrt{\mathcal E_N(t)}\mathcal D_N(t).
\end{align*}
By H\"older's inequality, \eqref{Sobolev-ineq}, \eqref{norm inequality}, and the definition of $\widetilde{\mathcal D}_{N,l}(t)$ in \eqref{dissipation functional}, we obtain
\begin{align*}
I_{1,3}
\lesssim\;&
\varepsilon\|E^\varepsilon\|_{L^\infty}
\left\|\langle v\rangle^{1-\frac{\gamma+2s}{2}}\partial^\alpha\mathbf P^\perp f^\varepsilon\right\|
\left\|\partial^\alpha\mathbf P^\perp f^\varepsilon\right\|_{D} \\
&+\varepsilon\Bigg\{\sum_{|\alpha_1|=1, \alpha_1\neq\alpha}
\|\partial^{\alpha_1}E^\varepsilon\|_{L^\infty}
\left\|\langle v\rangle^{1-\frac{\gamma+2s}{2}}\partial^{\alpha-\alpha_1}\mathbf P^\perp f^\varepsilon\right\|
+\|\partial^\alpha E^\varepsilon\|
\left\|\langle v\rangle^{1-\frac{\gamma+2s}{2}}\mathbf P^\perp f^\varepsilon\right\|_{L_x^\infty L_v^2}\\
&\qquad+
\sum_{2\leq|\alpha_1|\leq|\alpha|-1}\|\partial^{\alpha_1}E^\varepsilon\|_{L^6}
\left\|\langle v\rangle^{1-\frac{\gamma+2s}{2}}\partial^{\alpha-\alpha_1}\mathbf P^\perp f^\varepsilon\right\|_{L_x^3L_v^2}\Bigg\}
\|\partial^\alpha\mathbf P^\perp f^\varepsilon\|_D\\
\lesssim\;&
\bigg\{ (1+t)^{\frac{\vartheta}{2}}
\|E^\varepsilon\|_{L^\infty}\left\|w_l(\alpha,0)\partial^\alpha\mathbf P^\perp f^\varepsilon\right\|_D
+\|E^\varepsilon\|_{H^N}
\left\|\langle v\rangle^{1-\frac{\gamma+2s}{2}}\mathbf P^\perp f^\varepsilon\right\|_{H_x^{N-1}L_v^2}\bigg\}
\|\partial^\alpha\mathbf P^\perp f^\varepsilon\|_D\\
\lesssim\;&
\Big\{(1+t)^{\frac{1+\varepsilon_0}{4}+\frac{\vartheta}{2}}
\|E^\varepsilon\|_{L^\infty}+\|E^\varepsilon\|_{H^N}\Big\}\widetilde{\mathcal D}_{N,l}(t)\\
\lesssim\;&
\Big\{\sqrt{\mathcal E_N(t)}+\delta_0\Big\}\widetilde{\mathcal D}_{N,l}(t).
\end{align*}
Here, we have used the weight inequalities
\begin{equation*}
\begin{aligned}
&\langle v\rangle^{1-\frac{\gamma+2s}{2}}
\lesssim (1+t)^{\frac{\vartheta}{2}}
\langle v\rangle^{\frac{\gamma+2s}{2}}w_l(\alpha_1,0),
&\quad& \text{for} \;\;|\alpha_1|\leq N,\\
&\langle v\rangle^{1-\frac{\gamma+2s}{2}}
\leq \langle v\rangle^{\frac{\gamma+2s}{2}}w_l(\alpha_1,0),
&\quad& \text{for} \;\;|\alpha_1|\leq N-1,
\end{aligned}
\end{equation*}
which follow from \eqref{weight function}, \eqref{weight function2}, and $l\geq mN$, together with \eqref{priori-assumptions-decay} and the condition
\begin{align*}
\frac{1+\varepsilon_0}{4}+\frac{\vartheta}{2}
\leq\frac{3+2\varrho}{8}.
\end{align*}
Combining the above estimates with those of $I_{1,1}$ and $I_{1,2}$, we obtain, for $1\leq |\alpha|\leq N$,
\begin{align}\label{I1-positive-estimate}
I_1 \lesssim
\Big\{\sqrt{\mathcal E_N(t)}+\delta_0\Big\} \widetilde{\mathcal D}_{N,l}(t).
\end{align}
Hence, it follows from \eqref{I1-zero-estimate} and
\eqref{I1-positive-estimate} that, for $|\alpha|\leq N$,
\begin{align}\label{I1-estimate}
I_1\leq
\chi_{|\alpha|=0}\frac{\varepsilon}{2}\frac{\d}{\d t}\int_{\mathbb R^3}(a^\varepsilon+2c^\varepsilon)|b^\varepsilon|^2\d x
+C\Big\{\sqrt{\widetilde{\mathcal E}_{N,l}(t)}+\mathcal E_N(t)+\delta_0\Big\}\widetilde{\mathcal D}_{N,l}(t).
\end{align}

For the term $I_2$, it suffices to consider $1\leq|\alpha|\leq N$.
By \eqref{Pf define}, a direct calculation gives
\begin{align}\label{B important1}
(v\times\partial^{\alpha_1}B^\varepsilon)\cdot\nabla_v\partial^{\alpha-\alpha_1}\mathbf Pf^\varepsilon
=-\Big(\partial^{\alpha-\alpha_1}b^\varepsilon\times\partial^{\alpha_1}B^\varepsilon\Big)\cdot v\mu^{1/2}
\in\mathcal N(L).
\end{align}
Thus, by \eqref{f decomposition}, integration by parts in $v$, the orthogonality of $\mathbf P^\perp f^\varepsilon$ to $\mathcal N(L)$ and $\nabla_v\cdot(v\times B^\varepsilon)=0$, we obtain
\begin{align}\label{B important2}
I_2
=\;&-\sum_{0<\alpha_1\leq\alpha}C_\alpha^{\alpha_1}
\Big( (v\times\partial^{\alpha_1}B^\varepsilon)\cdot\nabla_v\partial^{\alpha-\alpha_1}f^\varepsilon,
\partial^\alpha f^\varepsilon\Big) \nonumber\\
=\;&
\sum_{0<\alpha_1\leq\alpha}C_\alpha^{\alpha_1}
\Big( \partial^{\alpha-\alpha_1}b^\varepsilon\times\partial^{\alpha_1}B^\varepsilon, \partial^\alpha b^\varepsilon \Big) \nonumber\\
&-\sum_{0<\alpha_1\leq\alpha}C_\alpha^{\alpha_1}
\Big((v\times\partial^{\alpha_1}B^\varepsilon)\cdot
\nabla_v\partial^{\alpha-\alpha_1}\mathbf P^\perp f^\varepsilon,\partial^\alpha\mathbf P^\perp f^\varepsilon\Big) \nonumber\\
:=\;&I_{2,1}+I_{2,2}.
\end{align}
An application of H\"{o}lder's inequality, together with \eqref{Sobolev-ineq}, yields
\begin{align*}
I_{2,1}
\lesssim\;& \bigg\{\|b^\varepsilon\|_{L_x^\infty}\|\partial^\alpha B^\varepsilon\|
+\sum_{1\leq|\alpha_1|\leq|\alpha|-1}
\|\partial^{\alpha_1}B^\varepsilon\|_{L^3}
\|\partial^{\alpha-\alpha_1}b^\varepsilon\|_{L^6}\bigg\}
\|\partial^\alpha b^\varepsilon\|\\
\lesssim\;& \|B^\varepsilon\|_{H^N}
\|\nabla_x b^\varepsilon\|_{H^{N-1}}
\|\partial^\alpha b^\varepsilon\|\\
\lesssim\;& \sqrt{\mathcal E_N(t)}\mathcal D_N(t).
\end{align*}
For the term $I_{2,2}$, we claim
\begin{align}\label{I22-duality}
I_{2,2}
\lesssim\;&\sum_{0<\alpha_1\leq\alpha}
\int_{\mathbb R^3}|\partial^{\alpha_1}B^\varepsilon|
\left|\langle v\rangle^{1-\frac{\gamma+2s}{2}+\sigma}\partial^{\alpha-\alpha_1}\mathbf P^\perp f^\varepsilon\right|_{H_v^{1-\sigma}}
|\partial^\alpha\mathbf P^\perp f^\varepsilon|_D\d x.
\end{align}
Here and below, we set $\sigma=\min\{s,1/2\}$. Indeed, using Sobolev duality in $v$, the fact that $(v\times\partial^{\alpha_1}B^\varepsilon)\cdot v=0$, and Lemma \ref{lem:D-fractional-control}, we obtain
\begin{align}\label{new B estimate}
\big\langle(v\times\partial^{\alpha_1}B^\varepsilon)\cdot\nabla_v g,h\big\rangle
=\;&\left\langle\frac{v\times\partial^{\alpha_1}B^\varepsilon}{\langle v\rangle}\cdot\nabla_v
\left(\langle v\rangle^{1-\frac{\gamma+2s}{2}+\sigma}g\right),\langle v\rangle^{\frac{\gamma+2s}{2}-\sigma}h\right\rangle \nonumber\\
\lesssim\;&
|\partial^{\alpha_1}B^\varepsilon|
\left|\langle v\rangle^{1-\frac{\gamma+2s}{2}+\sigma}g\right|_{H_v^{1-\sigma}}
\left|\langle v\rangle^{\frac{\gamma+2s}{2}-\sigma}h\right|_{H_v^\sigma}\\
\lesssim\;&
|\partial^{\alpha_1}B^\varepsilon|\left|\langle v\rangle^{1-\frac{\gamma+2s}{2}+\sigma}g\right|_{H_v^{1-\sigma}}|h|_D. \nonumber
\end{align}
We next estimate the second factor in the integrand on the right-hand side of \eqref{I22-duality}.
For $0<\sigma<1/2$, applying Lemma \ref{lem:weighted-interpolation} to $e^{\frac{q\langle v\rangle^2}{(1+t)^{\vartheta}}}g$ with
\begin{align*}
&\theta=1-2\sigma,\qquad r_1=1-\sigma,\qquad r_2=\sigma,\qquad r_3=1+\sigma,\\
&\ell_2=\frac{\gamma+2s}{2}-\sigma+l-(m+\gamma+2s)|\overline\alpha|,\qquad \ell_3=\ell_2-m,\\
&\ell_1=2\sigma\ell_2+(1-2\sigma)\ell_3=1-\frac{\gamma+2s}{2}+\sigma+l-(m+\gamma+2s)(|\overline\alpha|+1),
\end{align*}
we obtain, for $|\overline{\alpha}|\leq N-1$,
\begin{align}\label{I22-velocity-interpolation}
\left|\langle v\rangle^{1-\frac{\gamma+2s}{2}+\sigma}g\right|_{H_v^{1-\sigma}}
\lesssim\;& \left|\langle v\rangle^{1-\frac{\gamma+2s}{2}+\sigma} w_{l}(|\overline\alpha|+1,0) g\right|_{H_v^{1-\sigma}}
= \left|\langle v\rangle^{\ell_1}e^{\frac{q\langle v\rangle^2}{(1+t)^{\vartheta}}}g\right|_{H_v^{1-\sigma}}\nonumber\\
\lesssim\;&
\left|\langle v\rangle^{\ell_2}e^{\frac{q\langle v\rangle^2}{(1+t)^{\vartheta}}}g\right|_{H_v^\sigma}^{2\sigma}
\left|\langle v\rangle^{\ell_3}e^{\frac{q\langle v\rangle^2}{(1+t)^{\vartheta}}}g\right|_{H_v^{1+\sigma}}^{1-2\sigma}\nonumber\\
=\;&
\left|\langle v\rangle^{\frac{\gamma+2s}{2}-\sigma}w_l(\overline\alpha,0)g\right|_{H_v^\sigma}^{2\sigma}
\left|\langle v\rangle^{\frac{\gamma+2s}{2}-\sigma}w_l(\overline\alpha,e_i)g\right|_{H_v^{1+\sigma}}^{1-2\sigma}\nonumber\\
\lesssim\;&
|w_l(\overline\alpha,0)g|_D^{2\sigma}\Big\{|w_l(\overline\alpha,0) g|_D+|w_l(|\overline\alpha|,1)\nabla_v g|_D\Big\}^{1-2\sigma}.
\end{align}
Here we have used $2\sigma(m-1)=1-2(\gamma+2s)$ and $l-(m+\gamma+2s)(|\overline\alpha|+1)\geq0$, together with the estimate
\begin{align*}
\left|\langle v\rangle^{\frac{\gamma+2s}{2}-\sigma}w_l(\overline\alpha,e_i) g\right|_{H_v^{1+\sigma}}
\lesssim
|w_l(\overline\alpha,0) g|_D+|w_l(|\overline\alpha|,1)\nabla_v g|_D,
\end{align*}
which follows from Lemma \ref{lem:D-fractional-control} and $|\nabla_vw_l(\overline\alpha,e_i)|\lesssim\langle v\rangle w_l(\overline\alpha,e_i)\lesssim w_l(\overline\alpha,0)$.
For $\sigma=1/2$ and $|\overline\alpha|\leq N-1$, we obtain directly
\begin{align}\label{I22-velocity-interpolation2}
\!\!\!\left|\langle v\rangle^{\frac32-\frac{\gamma+2s}{2}} g\right|_{H_v^{\frac{1}{2}}}
\lesssim \left|\langle v\rangle^{\frac32-\frac{\gamma+2s}{2}} w_l(|\overline\alpha|+1,0) g\right|_{H_v^{\frac{1}{2}}}
\lesssim \left|\langle v\rangle^{\frac{\gamma+2s}{2}-\frac12} w_l(\overline\alpha,0) g\right|_{H_v^{\frac{1}{2}}}
\lesssim |w_l(\overline\alpha,0) g|_D.
\end{align}
Inserting \eqref{I22-velocity-interpolation} and \eqref{I22-velocity-interpolation2} into \eqref{I22-duality}, and applying H\"older's inequality together with \eqref{Sobolev-ineq}, we obtain
\begin{align*}
I_{2,2}
\lesssim\;&
\Bigg\{\sum_{1\leq|\alpha_1|\leq2}\|\partial^{\alpha_1}B^\varepsilon\|_{L^\infty}
\left\|\langle v\rangle^{1-\frac{\gamma+2s}{2}+\sigma}\partial^{\alpha-\alpha_1}\mathbf P^\perp f^\varepsilon\right\|_{L^2_x H^{1-\sigma}_v}
+\|\partial^\alpha B^\varepsilon\| \left\|\langle v\rangle^{1-\frac{\gamma+2s}{2}+\sigma}\mathbf P^\perp f^\varepsilon\right\|_{L_x^\infty H^{1-\sigma}_v}\\
&\quad+
\sum_{3\leq|\alpha_1|\leq|\alpha|-1}\|\partial^{\alpha_1}B^\varepsilon\|_{L^6}
\left\|\langle v\rangle^{1-\frac{\gamma+2s}{2}+\sigma}\partial^{\alpha-\alpha_1}\mathbf P^\perp f^\varepsilon\right\|_{L_x^3 H^{1-\sigma}_v}\Bigg\}
\|\partial^\alpha\mathbf P^\perp f^\varepsilon\|_D\\
\lesssim\;&
\|\nabla_xB^\varepsilon\|_{L^\infty}
\sum_{|\alpha_1|=1}\left\|\langle v\rangle^{1-\frac{\gamma+2s}{2}+\sigma}\partial^{\alpha-\alpha_1}\mathbf P^\perp f^\varepsilon\right\|_{L^2_x H^{1-\sigma}_v}
\|\partial^\alpha\mathbf P^\perp f^\varepsilon\|_D \\
&+\|B^\varepsilon\|_{H^N}
\left\|\langle v\rangle^{1-\frac{\gamma+2s}{2}+\sigma}\mathbf P^\perp f^\varepsilon\right\|_{H_x^{N-2} H^{1-\sigma}_v}
\|\partial^\alpha\mathbf P^\perp f^\varepsilon\|_D \\
\lesssim\;&
\Big\{(1+t)^{\frac{(1-2s)^+(1+\varepsilon_0)}4}\|\nabla_xB^\varepsilon\|_{L^\infty}+\|B^\varepsilon\|_{H^N}\Big\}
\widetilde{\mathcal D}_{N,l}(t)\nonumber\\
\lesssim\;&
\Big\{\sqrt{\mathcal E_N(t)}+\delta_0\Big\}\widetilde{\mathcal D}_{N,l}(t),
\end{align*}
where we have used $1-2\sigma=(1-2s)^+$ and \eqref{priori-assumptions-decay}.
Combining the estimates of $I_{2,1}$ and $I_{2,2}$, we conclude that, for $|\alpha|\leq N$,
\begin{align}\label{I2-estimate}
I_2 \lesssim \Big\{\sqrt{\mathcal E_N(t)}+\delta_0\Big\} \widetilde{\mathcal D}_{N,l}(t).
\end{align}

For the term $I_3$, it suffices to consider $1\leq|\alpha|\leq N$.
By \eqref{f decomposition} and integration by parts, we obtain
\begin{align*}
I_3
=\;&-\varepsilon\sum_{0<\alpha_1\leq\alpha}C_\alpha^{\alpha_1}
\Big(\partial^{\alpha_1}E^\varepsilon\cdot\nabla_v\partial^{\alpha-\alpha_1}f^\varepsilon,
\partial^\alpha f^\varepsilon\Big)\\
=\;&-\varepsilon\sum_{0<\alpha_1\leq\alpha}C_\alpha^{\alpha_1}
\Big(\partial^{\alpha_1}E^\varepsilon\cdot\nabla_v\partial^{\alpha-\alpha_1}\mathbf Pf^\varepsilon,
\partial^\alpha f^\varepsilon\Big)\\
&+\varepsilon\sum_{0<\alpha_1\leq\alpha}C_\alpha^{\alpha_1}
\Big(\partial^{\alpha_1}E^\varepsilon\cdot\nabla_v\partial^\alpha\mathbf Pf^\varepsilon,
\partial^{\alpha-\alpha_1}\mathbf P^\perp f^\varepsilon\Big)\\
&-\varepsilon\sum_{0<\alpha_1\leq\alpha}C_\alpha^{\alpha_1}
\Big(\partial^{\alpha_1}E^\varepsilon\cdot\nabla_v\partial^{\alpha-\alpha_1}\mathbf P^\perp f^\varepsilon,
\partial^\alpha\mathbf P^\perp f^\varepsilon\Big)\\
:=\;&I_{3,1}+I_{3,2}+I_{3,3}.
\end{align*}
By \eqref{Pf define}, \eqref{norm inequality}, H\"{o}lder's inequality and \eqref{Sobolev-ineq}, we obtain
\begin{align*}
I_{3,1}
\lesssim\;&
\varepsilon\bigg\{\|\partial^\alpha E^\varepsilon\| \|\mathbf Pf^\varepsilon\|_{L_x^\infty L_v^2}
+\sum_{1\leq|\alpha_1|\leq|\alpha|-1}
\|\partial^{\alpha_1}E^\varepsilon\|_{L^3}
\|\partial^{\alpha-\alpha_1}\mathbf Pf^\varepsilon\|_{L_x^6L_v^2}\bigg\}
\|\partial^\alpha f^\varepsilon\|_D\\
\lesssim\;&
\varepsilon\|E^\varepsilon\|_{H^N}\|\nabla_x\mathbf Pf^\varepsilon\|_{H_x^{N-1}L_v^2}
\Big\{\|\partial^\alpha\mathbf Pf^\varepsilon\|+\|\partial^\alpha\mathbf P^\perp f^\varepsilon\|_D\Big\}\\
\lesssim\;&
\sqrt{\mathcal E_N(t)}\mathcal D_N(t).
\end{align*}
The term $I_{3,2}$ is treated in exactly the same manner, yielding
\begin{align*}
I_{3,2}\lesssim\sqrt{\mathcal E_N(t)}\mathcal D_N(t).
\end{align*}
For the term $I_{3,3}$, arguing as in the estimate of $I_{2,2}$ and using
\eqref{priori-assumptions-decay}, we obtain
\begin{align*}
I_{3,3}
\lesssim \Big\{(1+t)^{\frac{(1-2s)^+(1+\varepsilon_0)}4}\|\nabla_xE^\varepsilon\|_{L^\infty}+\|E^\varepsilon\|_{H^N}\Big\}
\widetilde{\mathcal D}_{N,l}(t)
\lesssim\Big\{\sqrt{\mathcal E_N(t)}+\delta_0\Big\}\widetilde{\mathcal D}_{N,l}(t).
\end{align*}
Combining the estimates of $I_{3,1}$--$I_{3,3}$, we conclude that, for $|\alpha|\leq N$,
\begin{align}\label{I3-estimate}
I_3 \lesssim \Big\{\sqrt{\mathcal E_N(t)}+\delta_0\Big\} \widetilde{\mathcal D}_{N,l}(t).
\end{align}

For the term $I_4$, \eqref{soft gamma3} yields
\begin{align}\label{I4-estimate}
I_4\lesssim \sqrt{\mathcal E_N(t)}\mathcal D_N(t).
\end{align}

Substituting \eqref{I1-estimate}--\eqref{I4-estimate} into \eqref{unweighted-energy-basic} and summing over $|\alpha|\leq N$, we obtain \eqref{ss:unweighted-estimate}.
This completes the proof.
\end{proof}

\smallskip

\subsection{Weighted Energy Estimate}
\hspace*{\fill}

In this subsection, we derive energy estimates with the weight function $\widetilde{w}_l(\alpha,\beta)$ in the regime $\max\{-3,-\frac{3}{2}-2s\}<\gamma<-2s$.
We first establish the estimate for the highest-order mixed derivatives $|\alpha|+|\beta|=N$ with $|\beta|\geq1$.

\begin{lemma}\label{weighted 1}
Under the assumptions of Lemma \ref{ss:unweighted}, there holds
\begin{align}\label{weighted estimate1}
&\sum_{j=1}^{N}C_j
\sum_{\substack{|\alpha|+|\beta|=N\\|\beta|=j}}
\Bigg\{
\frac12\frac{\d}{\d t}
\left\|\widetilde w_l(\alpha,\beta)\partial_\beta^\alpha\mathbf P^\perp f^\varepsilon\right\|^2
+\frac{\lambda}{\varepsilon^2}\left\|\widetilde w_l(\alpha,\beta)\partial_\beta^\alpha\mathbf P^\perp f^\varepsilon\right\|_D^2\nonumber\\
&\qquad\qquad
+\frac{q\vartheta}{(1+t)^{1+\vartheta}}
\left\|\langle v\rangle\widetilde w_l(\alpha,\beta)\partial_\beta^\alpha\mathbf P^\perp f^\varepsilon\right\|^2
+\frac{1+\varepsilon_0}{4(1+t)}\left\|\widetilde w_l(\alpha,\beta)\partial_\beta^\alpha\mathbf P^\perp f^\varepsilon\right\|^2\Bigg\}\nonumber\\
\lesssim\;&
\frac1{\varepsilon^2}
\sum_{|\alpha|+|\beta|\leq N-1}\left\|\widetilde w_l(\alpha,\beta)\partial_\beta^\alpha\mathbf P^\perp f^\varepsilon\right\|_D^2
+\sum_{|\alpha|=N}\left\|\widetilde w_l(\alpha,0)\partial^\alpha\mathbf P^\perp f^\varepsilon\right\|_D^2\nonumber\\
&+\mathcal D_N(t)
+\Big\{\delta_0+\mathcal E_N(t)+\sqrt{\widetilde{\mathcal E}_{N,l}(t)}\Big\}\widetilde{\mathcal D}_{N,l}(t),
\end{align}
where $C_j>0$ are fixed constants satisfying $C_j\gg C_{j+1}$.
\end{lemma}

\begin{proof}
We apply the micro-projection operator $\mathbf P^\perp$ to the first equation of \eqref{rVMB} and use the macro-micro decomposition \eqref{f decomposition} to obtain
\begin{align}\label{rrVMB}
&\partial_t\mathbf P^\perp f^\varepsilon
+\frac1\varepsilon v\cdot\nabla_x\mathbf P^\perp f^\varepsilon
+\frac1{\varepsilon^2}L\mathbf P^\perp f^\varepsilon
+\big(\varepsilon E^\varepsilon+v\times B^\varepsilon\big)\cdot\nabla_v\mathbf P^\perp f^\varepsilon\nonumber\\
&=\frac{\varepsilon}{2}E^\varepsilon\cdot v\mathbf P^\perp f^\varepsilon
+\frac1\varepsilon\Gamma(f^\varepsilon,f^\varepsilon)
+[\mathbf P,\mathcal B]f^\varepsilon,
\end{align}
where $[\mathbf P,\mathcal B]:=\mathbf P\mathcal B-\mathcal B\mathbf P$ and the operator $\mathcal B$ is defined by
\begin{align*}
\mathcal B:=\frac1\varepsilon v\cdot\nabla_x+\big(\varepsilon E^\varepsilon+v\times B^\varepsilon\big)\cdot\nabla_v
-\frac{\varepsilon}{2}E^\varepsilon\cdot v.
\end{align*}
Let $|\alpha|+|\beta|=N$ and $|\beta|\geq1$.
Applying $\partial_\beta^\alpha$ to \eqref{rrVMB}, taking the inner product with $\widetilde w_l^2(\alpha,\beta)\partial_\beta^\alpha\mathbf P^\perp f^\varepsilon$, and using integration by parts and \eqref{L-coercive-weighted-2}, we obtain
\begin{align}\label{weighted:top-mixed-basic}
\begin{split}
&\frac12\frac{\d}{\d t}
\left\|\widetilde w_l(\alpha,\beta)\partial_\beta^\alpha\mathbf P^\perp f^\varepsilon\right\|^2
+\frac{\lambda}{\varepsilon^2}\left\|\widetilde w_l(\alpha,\beta)\partial_\beta^\alpha\mathbf P^\perp f^\varepsilon\right\|_D^2\\
&+\frac{q\vartheta}{(1+t)^{1+\vartheta}}
\left\|\langle v\rangle\widetilde w_l(\alpha,\beta)\partial_\beta^\alpha\mathbf P^\perp f^\varepsilon\right\|^2
+\frac{1+\varepsilon_0}{4(1+t)}
\left\|\widetilde w_l(\alpha,\beta)\partial_\beta^\alpha\mathbf P^\perp f^\varepsilon\right\|^2\\
\leq\;&
\frac{C}{\varepsilon^2}\sum_{|\beta'|<|\beta|}
\left\|\widetilde w_l(\alpha,\beta')\partial_{\beta'}^\alpha\mathbf P^\perp f^\varepsilon\right\|_D^2
-\frac1\varepsilon
\Big(\partial^\alpha_\beta( v\cdot \nabla_x \mathbf P^\perp f^\varepsilon),
\widetilde w_l^2(\alpha,\beta)\partial_\beta^\alpha\mathbf P^\perp f^\varepsilon\Big)\\
&-\Big(\partial_\beta^\alpha
\big((v\times B^\varepsilon)\cdot\nabla_v\mathbf P^\perp f^\varepsilon\big),
\widetilde w_l^2(\alpha,\beta)\partial_\beta^\alpha\mathbf P^\perp f^\varepsilon\Big)\\
&-\varepsilon\Big(\partial_\beta^\alpha\big(E^\varepsilon\cdot\nabla_v\mathbf P^\perp f^\varepsilon
-\frac{1}{2}E^\varepsilon\cdot v\mathbf P^\perp f^\varepsilon\big),
\widetilde w_l^2(\alpha,\beta)\partial_\beta^\alpha\mathbf P^\perp f^\varepsilon\Big)\\
&+\frac1\varepsilon
\Big(\partial_\beta^\alpha\Gamma(f^\varepsilon,f^\varepsilon),
\widetilde w_l^2(\alpha,\beta)\partial_\beta^\alpha\mathbf P^\perp f^\varepsilon\Big)
+\Big(\partial_\beta^\alpha[\mathbf P,\mathcal B]f^\varepsilon,
\widetilde w_l^2(\alpha,\beta)\partial_\beta^\alpha\mathbf P^\perp f^\varepsilon\Big).
\end{split}
\end{align}
Here we have used
\begin{align*}
\partial_t\widetilde w_l^2(\alpha,\beta)
=-\left\{ \frac{2q\vartheta\langle v\rangle^2}{(1+t)^{1+\vartheta}}
+\frac{1+\varepsilon_0}{2(1+t)} \right\}\widetilde w_l^2(\alpha,\beta),
\qquad \text{for}\;\; |\alpha|+|\beta|=N.
\end{align*}
We denote the last five terms on the right-hand side of \eqref{weighted:top-mixed-basic} by $I_i$ $(i=1,\cdots, 5)$, respectively.

The term $I_1$ can be handled by exploiting the weight relation
\begin{align*}
\widetilde w_l(\alpha,\beta)=\langle v\rangle^{\gamma+2s}\widetilde w_l(\alpha+e_i,\beta-e_i).
\end{align*}
Together with \eqref{norm inequality} and the Cauchy--Schwarz inequality, this implies that, for any $\eta>0$,
\begin{align*}
I_1
\lesssim\;&
\frac1\varepsilon
\left\|\widetilde w_l(\alpha+e_i,\beta-e_i)\partial_{\beta-e_i}^{\alpha+e_i}\mathbf P^\perp f^\varepsilon\right\|_D
\left\|\widetilde w_l(\alpha,\beta)\partial_\beta^\alpha\mathbf P^\perp f^\varepsilon\right\|_D\\
\lesssim\;&
\frac{\eta}{\varepsilon^2}
\left\|\widetilde w_l(\alpha,\beta)\partial_\beta^\alpha\mathbf P^\perp f^\varepsilon\right\|_D^2
+\left\|\widetilde w_l(\alpha+e_i,\beta-e_i)\partial_{\beta-e_i}^{\alpha+e_i}\mathbf P^\perp f^\varepsilon\right\|_D^2.
\end{align*}

Turning to the term $I_2$, we integrate by parts in $v$ and use $(v\times B^\varepsilon)\cdot\nabla_v\widetilde w^2_l(\alpha,\beta)=0$ to write
\begin{align*}
I_2=\;&
-\sum_{0<\alpha_1\leq\alpha}C_\alpha^{\alpha_1}
\Big((v\times\partial^{\alpha_1}B^\varepsilon)\cdot\nabla_v\partial_\beta^{\alpha-\alpha_1}\mathbf P^\perp f^\varepsilon,
\widetilde w_l^2(\alpha,\beta)\partial_\beta^\alpha\mathbf P^\perp f^\varepsilon\Big)\\
&-\sum_{\alpha_1\leq\alpha}C_\alpha^{\alpha_1} C_{\beta}^{e_i}
\Big((\partial_{e_i} v\times\partial^{\alpha_1}B^\varepsilon)\cdot\nabla_v\partial_{\beta-e_i}^{\alpha-\alpha_1}\mathbf P^\perp f^\varepsilon,
\widetilde w_l^2(\alpha,\beta)\partial_\beta^\alpha\mathbf P^\perp f^\varepsilon\Big)\\
:=\;&I_{2,1}+I_{2,2}.
\end{align*}
For the term $I_{2,1}$, it suffices to consider $1\leq|\alpha|\leq N-1$.
Since $(v\times \partial^{\alpha_1} B^\varepsilon)\cdot\nabla_v\widetilde w_l(\alpha,\beta)=0$, we can commute $\widetilde w_l(\alpha,\beta)$ with $(v\times \partial^{\alpha_1} B^\varepsilon)\cdot\nabla_v$ to obtain
\begin{align}\label{weighted:J1-start}
I_{2,1}=\;&-\sum_{0<\alpha_1\leq\alpha}C_\alpha^{\alpha_1}
\Big((v\times\partial^{\alpha_1}B^\varepsilon)\cdot\nabla_v \big(\widetilde{w}_l(\alpha,\beta)\partial_\beta^{\alpha-\alpha_1}\mathbf P^\perp f^\varepsilon\big),
\widetilde w_l(\alpha,\beta)\partial_\beta^\alpha\mathbf P^\perp f^\varepsilon\Big) \nonumber\\
\lesssim\;&
\sum_{0<\alpha_1\leq\alpha}\int_{\mathbb R^3}
|\partial^{\alpha_1}B^\varepsilon|
\left|\langle v\rangle^{1-\frac{\gamma+2s}{2}+\sigma}\widetilde{w}_l(\alpha,\beta)\partial_\beta^{\alpha-\alpha_1}\mathbf P^\perp f^\varepsilon\right|_{H_v^{1-\sigma}}
\left|\widetilde w_l(\alpha,\beta)\partial_\beta^\alpha\mathbf P^\perp f^\varepsilon\right|_D\d x.
\end{align}
Here, $\sigma=\min\{s,1/2\}$, and the last inequality follows from \eqref{new B estimate}.
We next estimate the second factor in the integrand on the right-hand side of \eqref{weighted:J1-start}.
For $0<\sigma<1/2$, applying Lemma \ref{lem:weighted-interpolation} to $e^{\frac{q\langle v\rangle^2}{(1+t)^\vartheta}}g$ with
\begin{align*}
&\theta=1-2\sigma,\qquad r_1=1-\sigma,\qquad r_2=\sigma,\qquad r_3=1+\sigma, \\
&\ell_2=\frac{\gamma+2s}{2}-\sigma+l-(m+\gamma+2s)(|\alpha|-1)-m|\beta|,\qquad\ell_3=\ell_2-m, \\
&\ell_1=2\sigma\ell_2+(1-2\sigma)\ell_3=1-\frac{\gamma+2s}{2}+\sigma+l-(m+\gamma+2s)|\alpha|-m|\beta|,
\end{align*}
we obtain
\begin{align}\label{weighted:top-mixed-magnetic-interpolation1}
&\left|\langle v\rangle^{1-\frac{\gamma+2s}{2}+\sigma}w_l(\alpha,\beta)g\right|_{H_v^{1-\sigma}}
=\left|\langle v\rangle^{\ell_1}e^{\frac{q\langle v\rangle^2}{(1+t)^\vartheta}}g\right|_{H_v^{1-\sigma}} \nonumber\\
\lesssim\;&
\left|\langle v\rangle^{\ell_2}e^{\frac{q\langle v\rangle^2}{(1+t)^\vartheta}}g\right|_{H_v^\sigma}^{2\sigma}
\left|\langle v\rangle^{\ell_3}e^{\frac{q\langle v\rangle^2}{(1+t)^\vartheta}}g\right|_{H_v^{1+\sigma}}^{1-2\sigma} \nonumber\\
=\;&
\left|\langle v\rangle^{\frac{\gamma+2s}{2}-\sigma}w_l(|\alpha|-1,|\beta|)g\right|_{H_v^\sigma}^{2\sigma}
\left|\langle v\rangle^{\frac{\gamma+2s}{2}-\sigma}w_l(|\alpha|-1,|\beta|+1)g\right|_{H_v^{1+\sigma}}^{1-2\sigma} \nonumber\\
\lesssim\;&
|w_l(|\alpha|-1,|\beta|)g|_D^{2\sigma}
\Big\{|w_l(|\alpha|-1,|\beta|)g|_D+|w_l(|\alpha|-1,|\beta|+1)\nabla_v g|_D\Big\}^{1-2\sigma}.
\end{align}
Here, we have used
$2\sigma(m-1)=1-2(\gamma+2s)$, together with the estimate
\begin{align*}
\left|\langle v\rangle^{\frac{\gamma+2s}{2}-\sigma}w_l(|\alpha|-1,|\beta|+1)g\right|_{H_v^{1+\sigma}}
\lesssim
|w_l(|\alpha|-1,|\beta|)g|_D+|w_l(|\alpha|-1,|\beta|+1)\nabla_{v}g|_D,
\end{align*}
which follows from Lemma \ref{lem:D-fractional-control} and
\begin{align*}
|\nabla_vw_l(|\alpha|-1,|\beta|+1)|
\lesssim \langle v\rangle w_l(|\alpha|-1,|\beta|+1)
\lesssim w_l(|\alpha|-1,|\beta|).
\end{align*}
For $\sigma=1/2$, the definition of the weight function in \eqref{weight function2} and Lemma \ref{lem:D-fractional-control} give directly
\begin{align*}
\left|\langle v\rangle^{\frac32-\frac{\gamma+2s}{2}}w_l(\alpha,\beta)g\right|_{H_v^{\frac{1}{2}}}
=\left|\langle v\rangle^{\frac{\gamma+2s}{2}-\frac12}w_l(|\alpha|-1,|\beta|)g\right|_{H_v^{\frac{1}{2}}}
\lesssim |w_l(|\alpha|-1,|\beta|)g|_D.
\end{align*}
Consequently, for $\sigma=\min\{s,1/2\}$, we have
\begin{align}\label{weighted:top-mixed-magnetic-interpolation}
&\left|\langle v\rangle^{1-\frac{\gamma+2s}{2}+\sigma}w_l(\alpha,\beta)g\right|_{H_v^{1-\sigma}}
\lesssim
\sum_{|\beta_1|\leq1}|w_l(|\alpha|-1,|\beta+\beta_1|)\partial_{\beta_1} g|_D.
\end{align}
Hence, by H\"older's inequality and \eqref{Sobolev-ineq}, \eqref{weighted:J1-start} and \eqref{weighted:top-mixed-magnetic-interpolation} yield
\begin{align*}
I_{2,1}
\lesssim\;&
(1+t)^{-\frac{1+\varepsilon_0}{4}}
\sum_{\substack{|\alpha_1|=1\\|\beta_1|\leq1}}
\|\partial^{\alpha_1}B^\varepsilon\|_{L^\infty}
\left\|w_l(|\alpha|-1,|\beta+\beta_1|)\partial_{\beta+\beta_1}^{\alpha-\alpha_1}\mathbf P^\perp f^\varepsilon\right\|_D
\left\|\widetilde w_l(\alpha,\beta)\partial_\beta^\alpha\mathbf P^\perp f^\varepsilon\right\|_D\\
&+(1+t)^{-\frac{1+\varepsilon_0}{4}}
\sum_{\substack{2\leq|\alpha_1|\leq|\alpha|\\|\beta_1|\leq 1}}
\|\partial^{\alpha_1}B^\varepsilon\|_{L^6}
\left\|w_l(|\alpha|-1,|\beta+\beta_1|)\partial_{\beta+\beta_1}^{\alpha-\alpha_1}\mathbf P^\perp f^\varepsilon\right\|_{L_x^3L_D^2}
\left\|\widetilde w_l(\alpha,\beta)\partial_\beta^\alpha\mathbf P^\perp f^\varepsilon\right\|_D\\
\lesssim\;&
\|B^\varepsilon\|_{H^N}\widetilde{\mathcal D}_{N,l}(t)
\lesssim \sqrt{\mathcal E_N(t)}\widetilde{\mathcal D}_{N,l}(t).
\end{align*}
Here, we have used $|\alpha+\alpha'|\leq N$ for $|\alpha'|=1$ since $|\alpha|\leq N-1$ and
\begin{align*}
(1+t)^{-\frac{1+\varepsilon_0}{4}} w_l(|\alpha|-1,|\beta+\beta_1|)\leq \widetilde{w}_l(\alpha-\alpha_1+\alpha',\beta+\beta_1),
\quad \text{for}\;\; |\alpha'|\leq 1\;\text{and} \; 2\leq |\alpha_1|\leq |\alpha|.
\end{align*}
As for the term $I_{2,2}$, recalling that $\partial_{e_i}v=e_i$, we express it as
\begin{align*}
I_{2,2}
=-\sum_{\alpha_1\leq\alpha}
C_\alpha^{\alpha_1}C_{\beta}^{e_i}
\Big( (e_i\times\partial^{\alpha_1}B^\varepsilon)\cdot\nabla_v\partial_{\beta-e_i}^{\alpha-\alpha_1}\mathbf P^\perp f^\varepsilon,
\widetilde w_l^2(\alpha,\beta)\partial_\beta^\alpha\mathbf P^\perp f^\varepsilon\Big).
\end{align*}
Sobolev duality in $v$, combined with Lemma \ref{lem:D-fractional-control}, leads to
\begin{align}\label{weighted:I22-duality}
&\left|\big\langle(e_i\times\partial^{\alpha_1}B^\varepsilon)\cdot\nabla_v g,
\widetilde w_l^2(\alpha,\beta)h\big\rangle\right|\nonumber\\
\lesssim\;&
\left|\left\langle(e_i\times\partial^{\alpha_1}B^\varepsilon)\cdot
 \nabla_v \big(\langle v\rangle^{-\frac{\gamma+2s}{2}+\sigma} \widetilde{w}_l(\alpha,\beta)g\big),
\langle v\rangle^{\frac{\gamma+2s}{2}-\sigma}\widetilde w_l(\alpha,\beta)h\right\rangle\right|\nonumber\\
&+\left|\left\langle(e_i\times\partial^{\alpha_1}B^\varepsilon)\cdot
\nabla_v \big(\langle v\rangle^{-\frac{\gamma+2s}{2}+\sigma} \widetilde{w}_l(\alpha,\beta) \big) g,
\langle v\rangle^{\frac{\gamma+2s}{2}-\sigma}\widetilde w_l(\alpha,\beta)h\right\rangle\right| \nonumber\\
\lesssim\;&
|\partial^{\alpha_1}B^\varepsilon|
\bigg\{\left|\langle v\rangle^{-\frac{\gamma+2s}{2}+\sigma}\widetilde w_l(\alpha,\beta)g\right|_{H_v^{1-\sigma}}
+\left|\langle v\rangle^{1-\frac{\gamma+2s}{2}}\widetilde w_l(\alpha,\beta)g\right|_{L_v^2}\bigg\}
|\widetilde w_l(\alpha,\beta)h|_D.
\end{align}
Here, we have used
\begin{align*}
\left|\nabla_v\log\Big(\langle v\rangle^{-\frac{\gamma+2s}{2}+\sigma}\widetilde w_l(\alpha,\beta)\Big)\right|\lesssim \langle v\rangle.
\end{align*}
The first velocity norm on the right-hand side of \eqref{weighted:I22-duality} is estimated by the same weighted interpolation argument as for $I_{2,1}$.
For the second velocity norm on the right-hand side of \eqref{weighted:I22-duality}, since $m\geq1-(\gamma+2s)$, the definition of the weight function in \eqref{weight function} gives
\begin{align*}
\langle v\rangle^{1-\frac{\gamma+2s}{2}} w_l(\alpha,\beta)
\leq \langle v\rangle^{\frac{\gamma+2s}{2}} w_l(|\alpha|,|\beta|-1).
\end{align*}
Consequently, applying \eqref{weighted:I22-duality} with $g=\partial_{\beta-e_i}^{\alpha-\alpha_1}\mathbf P^\perp f^\varepsilon$ and $h=\partial_\beta^\alpha\mathbf P^\perp f^\varepsilon$, \eqref{norm inequality}, H\"older's inequality, \eqref{Sobolev-ineq} and the definition of $\widetilde{\mathcal D}_{N,l}(t)$, we obtain
\begin{align*}
I_{2,2}
\lesssim
\|B^\varepsilon\|_{H^N}\widetilde{\mathcal D}_{N,l}(t)
\lesssim
\sqrt{\mathcal E_N(t)}\widetilde{\mathcal D}_{N,l}(t).
\end{align*}
Combining the estimates of $I_{2,1}$ and $I_{2,2}$, we conclude that
\begin{align*}
I_2\lesssim\sqrt{\mathcal E_N(t)}\widetilde{\mathcal D}_{N,l}(t).
\end{align*}

For the term $I_3$, by integration by parts and H\"{o}lder's inequality, we obtain
\begin{align}\label{weighted:mixed-electric-decomposition}
I_3
\lesssim\;&
\varepsilon\|E^\varepsilon\|_{L^\infty}
\left\|\langle v\rangle^{\frac{1}{2}}\widetilde w_l(\alpha,\beta)\partial_\beta^\alpha\mathbf P^\perp f^\varepsilon\right\|^2\nonumber\\
&+\varepsilon
\sum_{0<\alpha_1\leq\alpha}C_\alpha^{\alpha_1}
\left|\Big(\partial^{\alpha_1}E^\varepsilon\cdot\nabla_v\partial_\beta^{\alpha-\alpha_1}\mathbf P^\perp f^\varepsilon,
\widetilde w_l^2(\alpha,\beta)\partial_\beta^\alpha\mathbf P^\perp f^\varepsilon\Big)\right|\nonumber\\
&+\frac{\varepsilon}{2}
\sum_{\substack{ \alpha_1\leq\alpha,\ |\beta_1|\leq 1 \\|\alpha_1|+|\beta_1|\geq 1}}C_\alpha^{\alpha_1}C_\beta^{\beta_1}
\left|\Big(\partial^{\alpha_1}E^\varepsilon\cdot \partial_{\beta_1}v\partial_{\beta-\beta_1}^{\alpha-\alpha_1}\mathbf P^\perp f^\varepsilon,
\widetilde w_l^2(\alpha,\beta)\partial_\beta^\alpha\mathbf P^\perp f^\varepsilon\Big)\right|.
\end{align}
Repeating the Sobolev duality and weighted interpolation arguments used for $I_{2,2}$ and $I_{2,1}$, we can estimate the second term on the right-hand side of \eqref{weighted:mixed-electric-decomposition} in the same manner.
For the third term, we consider the cases $|\beta_1|=0$ and $|\beta_1|=1$ separately.
The definition of the weight function in \eqref{weight function} yields
\begin{equation*}
\begin{aligned}
\langle v\rangle^{1-\frac{\gamma+2s}{2}}w_l(\alpha,\beta)
&\leq \langle v\rangle^{\frac{\gamma+2s}{2}}w_l(|\alpha|-1,|\beta|),&\qquad &\text{for }|\alpha|\geq1,\\
\langle v\rangle^{-\frac{\gamma+2s}{2}}w_l(\alpha,\beta)
&\leq \langle v\rangle^{\frac{\gamma+2s}{2}}w_l(|\alpha|,|\beta|-1),&\qquad &\text{for }|\beta|\geq1,
\end{aligned}
\end{equation*}
since $m\geq1-2(\gamma+2s)$.
Consequently, we bound the last two terms on the right-hand side of \eqref{weighted:mixed-electric-decomposition} by
\begin{align*}
\|E^\varepsilon\|_{H^N}\widetilde{\mathcal D}_{N,l}(t) \lesssim \sqrt{\mathcal E_N(t)}\widetilde{\mathcal D}_{N,l}(t),
\end{align*}
where we have used H\"{o}lder's inequality, \eqref{norm inequality}, and \eqref{Sobolev-ineq}.
For the first term on the right-hand side of \eqref{weighted:mixed-electric-decomposition}, H\"older's inequality, \eqref{norm inequality} and the definition of $\widetilde{\mathcal D}_{N,l}(t)$ in \eqref{dissipation functional} yield
\begin{align}\label{E-new-estimate}
&\varepsilon\|E^\varepsilon\|_{L^\infty}
\left\|\langle v\rangle^{\frac{1}{2}}\widetilde w_l(\alpha,\beta)\partial_\beta^\alpha\mathbf P^\perp f^\varepsilon\right\|^2 \nonumber\\
\lesssim\;&
\|E^\varepsilon\|_{L^\infty}
\left\|\widetilde w_l(\alpha,\beta)\partial_\beta^\alpha\mathbf P^\perp f^\varepsilon\right\|_D^{\frac{2}{2-\gamma-2s}}
\left\|\langle v\rangle\widetilde w_l(\alpha,\beta)\partial_\beta^\alpha\mathbf P^\perp f^\varepsilon\right\|^{\frac{2(1-\gamma-2s)}{2-\gamma-2s}}
\nonumber\\
\lesssim\;&
(1+t)^{\frac{(1+\vartheta)(1-\gamma-2s)}{2-\gamma-2s}}
\|E^\varepsilon\|_{L^\infty}\widetilde{\mathcal D}_{N,l}(t)\nonumber\\
\lesssim\;&
\delta_0\widetilde{\mathcal D}_{N,l}(t).
\end{align}
Here, we have used the fact
\begin{align*}
\frac{(1+\vartheta)(1-\gamma-2s)}{2-\gamma-2s}\leq \frac{3+2\varrho}{8},
\end{align*}
for $\gamma+2s >-\frac{3}{2}$, $\max\left\{1,\frac{2-5(\gamma+2s)}{4-2(\gamma+2s)}\right\}<\varrho<\frac32$, and $0<\vartheta\leq \frac{(3+2\varrho)(2-\gamma-2s)}{8(1-\gamma-2s)}-1$, which, together with \eqref{priori-assumptions-decay}, implies
\begin{align*}
(1+t)^{\frac{(1+\vartheta)(1-\gamma-2s)}{2-\gamma-2s}}
\|E^\varepsilon(t)\|_{L^\infty}\leq\delta_0.
\end{align*}
Combining these estimates, we conclude that
\begin{align*}
I_3 \lesssim \Big\{\delta_0+\sqrt{\mathcal E_N(t)}\Big\}\widetilde{\mathcal D}_{N,l}(t).
\end{align*}

For the term $I_4$, applying \eqref{soft gamma6}, we obtain
\begin{align*}
I_4\lesssim\sqrt{\widetilde{\mathcal E}_{N,l}(t)}\widetilde{\mathcal D}_{N,l}(t).
\end{align*}

Regarding the commutator term $I_5$, it follows from \eqref{Pf define}, H\"{o}lder's inequality, and \eqref{Sobolev-ineq} that, for any $\eta>0$,
\begin{align*}
I_5
\lesssim\;& \frac{\eta}{\varepsilon^2} \left\| \widetilde{w}_l(\alpha,\beta)\partial^\alpha_\beta\mathbf{P}^{\perp}f^\varepsilon \right\|^2_{D}
+ \left\| \partial^\alpha \nabla_x \mathbf{P}^{\perp}f^\varepsilon \right\|_{D}^2
+ \left\| \partial^\alpha \nabla_x \mathbf{P} f^\varepsilon \right\|^2
+ \mathcal{E}_{N}(t)\mathcal{D}_{N}(t) \\
\lesssim\;& \frac{\eta}{\varepsilon^2} \left\| \widetilde{w}_l(\alpha,\beta)\partial^\alpha_\beta\mathbf{P}^{\perp}f^\varepsilon \right\|^2_{D}
+ \mathcal{D}_{N}(t) +\mathcal{E}_{N}(t)\mathcal{D}_{N}(t).
\end{align*}

Consequently, substituting the above estimates into \eqref{weighted:top-mixed-basic}, summing over $|\alpha|+|\beta|=N$ with $|\beta|=j$ for each $1\leq j\leq N$, and combining the resulting inequalities with suitably chosen constants $C_j\gg C_{j+1}>0$, we obtain the desired estimate \eqref{weighted estimate1} by taking $\eta>0$ sufficiently small.
This completes the proof.
\end{proof}
\smallskip

We next establish the weighted energy estimate for the highest-order pure spatial derivatives $|\alpha|=N$.

\begin{lemma}\label{weighted 2}
Under the assumptions of Lemma \ref{ss:unweighted}, there holds
\begin{align}\label{weighted estimate2}
&\frac{\varepsilon}{2}\frac{\d}{\d t}
\sum_{|\alpha|=N}\left\|\widetilde w_l(\alpha,0)\partial^\alpha f^\varepsilon\right\|^2
+\frac{\lambda}{\varepsilon}
\sum_{|\alpha|=N}\left\|\widetilde w_l(\alpha,0)\partial^\alpha\mathbf P^\perp f^\varepsilon\right\|_D^2\nonumber\\
&+\frac{\varepsilon q\vartheta}{(1+t)^{1+\vartheta}}
\sum_{|\alpha|=N}\left\|\langle v\rangle\widetilde w_l(\alpha,0)\partial^\alpha f^\varepsilon\right\|^2
+\frac{\varepsilon(1+\varepsilon_0)}{4(1+t)}
\sum_{|\alpha|=N}\left\|\widetilde w_l(\alpha,0)\partial^\alpha f^\varepsilon\right\|^2\nonumber\\
\lesssim\;&
\frac{\varepsilon^2}{(1+t)^{1+\varepsilon_0}}
\sum_{|\alpha|=N}\|\partial^\alpha E^\varepsilon\|^2
+\mathcal D_N(t)
+\Big\{\delta_0+\sqrt{\widetilde{\mathcal E}_{N,l}(t)}\Big\}\widetilde{\mathcal D}_{N,l}(t).
\end{align}
\end{lemma}

\begin{proof}
Applying $\partial^\alpha$ with $|\alpha|=N$ to the first equation in \eqref{rVMB}, taking the inner product with $\varepsilon\widetilde w_l^2(\alpha,0)\partial^\alpha f^\varepsilon$, and using integration by parts,
$L\mathbf Pf^\varepsilon=0$ and \eqref{L-coercive-weighted-1}, we obtain
\begin{align}\label{weighted:top-basic}
&\frac{\varepsilon}{2}\frac{\d}{\d t}
\left\|\widetilde w_l(\alpha,0)\partial^\alpha f^\varepsilon\right\|^2
+\frac{\lambda}{\varepsilon}
\left\|\widetilde w_l(\alpha,0)\partial^\alpha\mathbf P^\perp f^\varepsilon\right\|_D^2\nonumber\\
&+\frac{\varepsilon q\vartheta}{(1+t)^{1+\vartheta}}
\left\|\langle v\rangle\widetilde w_l(\alpha,0)\partial^\alpha f^\varepsilon\right\|^2
+\frac{\varepsilon(1+\varepsilon_0)}{4(1+t)}
\left\|\widetilde w_l(\alpha,0)\partial^\alpha f^\varepsilon\right\|^2\nonumber\\
\leq\;&
\frac{C}{\varepsilon}\|\partial^\alpha\mathbf P^\perp f^\varepsilon\|_D^2
-\frac1\varepsilon
\Big(L\partial^\alpha\mathbf P^\perp f^\varepsilon,\widetilde w_l^2(\alpha,0)\partial^\alpha\mathbf Pf^\varepsilon\Big)
+\varepsilon
\Big(\partial^\alpha E^\varepsilon\cdot v\mu^{1/2},\widetilde w_l^2(\alpha,0)\partial^\alpha f^\varepsilon\Big)\nonumber\\
&-\varepsilon
\Big(\partial^\alpha\big((v\times B^\varepsilon)\cdot\nabla_v f^\varepsilon\big),\widetilde w_l^2(\alpha,0)\partial^\alpha f^\varepsilon\Big)
-\varepsilon^2
\Big(\partial^\alpha(E^\varepsilon\cdot\nabla_v f^\varepsilon),\widetilde w_l^2(\alpha,0)\partial^\alpha f^\varepsilon\Big)\nonumber\\
&+\frac{\varepsilon^2}{2}
\Big(\partial^\alpha(E^\varepsilon\cdot vf^\varepsilon),\widetilde w_l^2(\alpha,0)\partial^\alpha f^\varepsilon\Big)
+\Big(\partial^\alpha\Gamma(f^\varepsilon,f^\varepsilon),\widetilde w_l^2(\alpha,0)\partial^\alpha f^\varepsilon\Big).
\end{align}
Here we have used
\begin{align*}
\partial_t\widetilde w_l^2(\alpha,0)
=-\left\{\frac{2q\vartheta\langle v\rangle^2}{(1+t)^{1+\vartheta}}
+\frac{1+\varepsilon_0}{2(1+t)}\right\}\widetilde w_l^2(\alpha,0),
\qquad \text{for}\;\;|\alpha|=N.
\end{align*}
We denote the last six terms on the right-hand side of \eqref{weighted:top-basic} by $I_i$ $(i=1,\cdots,6)$, respectively.

For the term $I_1$, by \eqref{L-upper-unweighted}, \eqref{Pf define}, and the Cauchy--Schwarz inequality, we obtain
\begin{align*}
I_1
\lesssim
\frac{1}{\varepsilon}
\|\partial^\alpha\mathbf P^\perp f^\varepsilon\|_D
\|\partial^\alpha\mathbf Pf^\varepsilon\|
\lesssim
\frac1{\varepsilon^2}\|\partial^\alpha\mathbf P^\perp f^\varepsilon\|_D^2+\|\partial^\alpha\mathbf Pf^\varepsilon\|^2
\lesssim\mathcal D_N(t).
\end{align*}

For the term $I_2$, by \eqref{Pf define}, \eqref{norm inequality}, and the Cauchy--Schwarz inequality, we have
\begin{align*}
I_2
\lesssim\;&
\varepsilon(1+t)^{-\frac{1+\varepsilon_0}{2}}
\|\partial^\alpha E^\varepsilon\|
\Big\{\|\partial^\alpha\mathbf Pf^\varepsilon\|
+\|\partial^\alpha\mathbf P^\perp f^\varepsilon\|_D\Big\}\\
\lesssim\;&
\frac{\varepsilon^2}{(1+t)^{1+\varepsilon_0}}
\|\partial^\alpha E^\varepsilon\|^2
+\|\partial^\alpha\mathbf Pf^\varepsilon\|^2
+\|\partial^\alpha\mathbf P^\perp f^\varepsilon\|_D^2\\
\lesssim\;&
\frac{\varepsilon^2}{(1+t)^{1+\varepsilon_0}}\|\partial^\alpha E^\varepsilon\|^2+\mathcal D_N(t).
\end{align*}

By \eqref{f decomposition}, \eqref{Pf define}, integration by parts in $v$, H\"older's inequality and \eqref{Sobolev-ineq}, we find that all terms in $I_3$--$I_5$ containing at least one macroscopic factor
are controlled by $\sqrt{\mathcal E_N(t)}\mathcal D_N(t)$.
Furthermore, following the arguments in the proof of Lemma \ref{weighted 1}, we can control the purely microscopic terms in $I_3$--$I_5$ by $\Big\{\delta_0+\sqrt{\mathcal E_N(t)}\Big\} \widetilde{\mathcal D}_{N,l}(t)$.

Finally, by \eqref{soft gamma5}, we have
\begin{align*}
I_6=\Big(\partial^\alpha\Gamma(f^\varepsilon,f^\varepsilon),\widetilde w_l^2(\alpha,0)\partial^\alpha f^\varepsilon\Big)
\lesssim
\sqrt{\widetilde{\mathcal E}_{N,l}(t)}\widetilde{\mathcal D}_{N,l}(t).
\end{align*}

Substituting the above estimates into \eqref{weighted:top-basic} and summing over $|\alpha|=N$, we obtain \eqref{weighted estimate2}.
This completes the proof.
\end{proof}
\smallskip

We next establish the weighted energy estimate for pure spatial derivatives $|\alpha|\leq N-1$.
In this case, $\widetilde w_l(\alpha,0)=w_l(\alpha,0)$.

\begin{lemma}\label{weighted 3}
Under the assumptions of Lemma \ref{ss:unweighted}, there holds
\begin{align}\label{weighted estimate3}
&\frac12\frac{\d}{\d t}
\sum_{|\alpha|\leq N-1}\left\|\widetilde w_l(\alpha,0)\partial^\alpha\mathbf P^\perp f^\varepsilon\right\|^2
+\frac{\lambda}{\varepsilon^2}\sum_{|\alpha|\leq N-1}
\left\|\widetilde w_l(\alpha,0)\partial^\alpha\mathbf P^\perp f^\varepsilon\right\|_D^2\nonumber\\
&+\frac{q\vartheta}{(1+t)^{1+\vartheta}}\sum_{|\alpha|\leq N-1}
\left\|\langle v\rangle\widetilde w_l(\alpha,0)\partial^\alpha\mathbf P^\perp f^\varepsilon\right\|^2\nonumber\\
\lesssim\;&
\mathcal D_N(t)+\Big\{\sqrt{\widetilde{\mathcal E}_{N,l}(t)}+\mathcal E_N(t)+\delta_0\Big\}\widetilde{\mathcal D}_{N,l}(t).
\end{align}
\end{lemma}
\begin{proof}
Applying $\partial^\alpha$ with $|\alpha|\leq N-1$ to \eqref{rrVMB}, taking the inner product with $\widetilde w_l^2(\alpha,0)\partial^\alpha\mathbf P^\perp f^\varepsilon$, and using integration by parts, \eqref{L-coercive-weighted-1}, and \eqref{norm inequality}, we obtain
\begin{align}\label{with weight 3}
&\frac12\frac{\d}{\d t}
\left\|\widetilde w_l(\alpha,0)\partial^\alpha\mathbf P^\perp f^\varepsilon\right\|^2
+\frac{\lambda}{\varepsilon^2}\left\|\widetilde w_l(\alpha,0)\partial^\alpha\mathbf P^\perp f^\varepsilon\right\|_D^2
+\frac{q\vartheta}{(1+t)^{1+\vartheta}}\left\|\langle v\rangle\widetilde w_l(\alpha,0)\partial^\alpha\mathbf P^\perp f^\varepsilon\right\|^2\nonumber\\
\leq\;&
\frac{C}{\varepsilon^2}\left\|\partial^\alpha\mathbf P^\perp f^\varepsilon\right\|_D^2
-\Big(\partial^\alpha\big((v\times B^\varepsilon)\cdot\nabla_v\mathbf P^\perp f^\varepsilon\big),
\widetilde w_l^2(\alpha,0)\partial^\alpha\mathbf P^\perp f^\varepsilon\Big)\nonumber\\
&-\varepsilon\Big(\partial^\alpha\big(E^\varepsilon\cdot\nabla_v\mathbf P^\perp f^\varepsilon\big),
\widetilde w_l^2(\alpha,0)\partial^\alpha\mathbf P^\perp f^\varepsilon\Big)
+\frac{\varepsilon}{2}\Big(\partial^\alpha\big(E^\varepsilon\cdot v\mathbf P^\perp f^\varepsilon\big),
\widetilde w_l^2(\alpha,0)\partial^\alpha\mathbf P^\perp f^\varepsilon\Big)\nonumber\\
&+\frac1\varepsilon\Big(\partial^\alpha\Gamma(f^\varepsilon,f^\varepsilon),
\widetilde w_l^2(\alpha,0)\partial^\alpha\mathbf P^\perp f^\varepsilon\Big)
+\Big(\partial^\alpha[\mathbf P,\mathcal B]f^\varepsilon,
\widetilde w_l^2(\alpha,0)\partial^\alpha\mathbf P^\perp f^\varepsilon\Big).
\end{align}
We denote the last five terms on the right-hand side of \eqref{with weight 3} by $I_i$ $(i=1,\cdots,5)$, respectively.

For the term $I_1$, it suffices to consider $1\leq|\alpha|\leq N-1$, since integration by parts and the identity $(v\times B^\varepsilon)\cdot v=0$ imply that $I_1=0$ when $|\alpha|=0$.
Repeating the weighted interpolation argument leading to \eqref{weighted:top-mixed-magnetic-interpolation1} and \eqref{weighted:top-mixed-magnetic-interpolation} with $\beta=0$, we obtain
\begin{align}\label{weighted:magnetic-interpolation}
\left|\langle v\rangle^{1-\frac{\gamma+2s}{2}+\sigma}w_l(\alpha,0)g \right|_{H_v^{1-\sigma}}
\lesssim \sum_{|\beta_1|\leq1} \big|w_l(|\alpha|-1,|\beta_1|)\partial_{\beta_1} g\big|_D.
\end{align}
Following the argument for $I_{2,1}$ in the proof of Lemma \ref{weighted 1}, and using \eqref{weighted:magnetic-interpolation}, H\"older's inequality and \eqref{Sobolev-ineq}, we obtain
\begin{align}\label{weighted:magnetic-start}
I_1
\lesssim\;&
\sum_{0<\alpha_1\leq\alpha}\int_{\mathbb R^3}|\partial^{\alpha_1}B^\varepsilon|
\left|\langle v\rangle^{1-\frac{\gamma+2s}{2}+\sigma}w_l(\alpha,0)\partial^{\alpha-\alpha_1}\mathbf P^\perp f^\varepsilon\right|_{H_v^{1-\sigma}}
\left|w_l(\alpha,0)\partial^\alpha\mathbf P^\perp f^\varepsilon\right|_D\d x \nonumber\\
\lesssim\;&
\sum_{\substack{0<\alpha_1\leq\alpha\\ |\beta_1|\leq 1}}
\int_{\mathbb R^3}|\partial^{\alpha_1}B^\varepsilon|
\left|w_l(|\alpha|-1,|\beta_1|)\partial_{\beta_1}^{\alpha-\alpha_1}\mathbf P^\perp f^\varepsilon\right|_D
\left|w_l(\alpha,0)\partial^\alpha\mathbf P^\perp f^\varepsilon\right|_D\d x \nonumber\\
\lesssim\;&
\|B^\varepsilon\|_{H^N} \widetilde{\mathcal D}_{N,l}(t)
\lesssim \sqrt{\mathcal E_N(t)} \widetilde{\mathcal D}_{N,l}(t).
\end{align}

For the terms $I_2$ and $I_3$, arguing as in the estimate of $I_3$ in the proof of Lemma \ref{weighted 1}, we obtain
\begin{align*}
I_2+I_3
\lesssim
\varepsilon\|E^\varepsilon\|_{L^\infty}
\left\|\langle v\rangle^{\frac{1}{2}}w_l(\alpha,0)\partial^\alpha\mathbf P^\perp f^\varepsilon\right\|^2
+\|E^\varepsilon\|_{H^N}\widetilde{\mathcal D}_{N,l}(t)
\lesssim
\Big\{\delta_0+\sqrt{\mathcal E_N(t)}\Big\}\widetilde{\mathcal D}_{N,l}(t),
\end{align*}
where the first term is controlled by the same velocity-weight interpolation and \eqref{E-new-estimate}.

For the term $I_4$, applying \eqref{soft gamma4} in Lemma \ref{soft Gamma}, we obtain
\begin{align*}
I_4 \lesssim\sqrt{\widetilde{\mathcal E}_{N,l}(t)}\widetilde{\mathcal D}_{N,l}(t).
\end{align*}

Finally, from the definition of $\mathbf{P}f^\varepsilon$ in \eqref{Pf define}, H\"{o}lder's inequality and \eqref{Sobolev-ineq}, we can derive that for any $\eta>0$,
\begin{align*}
I_5\lesssim \frac{\eta}{\varepsilon^2} \left\| \widetilde{w}_l(\alpha,0)\partial^\alpha\mathbf{P}^{\perp}f^\varepsilon \right\|^2_{D}
+ \mathcal{D}_{N}(t) +\mathcal{E}_{N}(t)\mathcal{D}_{N}(t).
\end{align*}

Substituting the above estimates into \eqref{with weight 3}, summing over $|\alpha|\leq N-1$, and choosing $\eta>0$ sufficiently small, we obtain \eqref{weighted estimate3}. This completes the proof.
\end{proof}
\smallskip

Finally, we establish the weighted energy estimate for mixed derivatives with $|\alpha|+|\beta|\leq N-1$ and $|\beta|\geq1$.

\begin{lemma}\label{weighted 4}
Under the assumptions of Lemma \ref{ss:unweighted}, there holds
\begin{align}\label{weighted estimate4}
&\sum_{j=1}^{N-1}\widetilde{C}_j \sum_{\substack{|\alpha|+|\beta|\leq N-1\\|\beta|=j}}
\Bigg\{
\frac12 \frac{\d}{\d t}\left\|\widetilde w_l(\alpha,\beta)\partial_\beta^\alpha\mathbf P^\perp f^\varepsilon\right\|^2
+\frac{\lambda}{\varepsilon^2}
\left\|\widetilde w_l(\alpha,\beta)\partial_\beta^\alpha\mathbf P^\perp f^\varepsilon\right\|_D^2\nonumber\\
&\qquad\qquad+\frac{q\vartheta}{(1+t)^{1+\vartheta}}
\left\|\langle v\rangle\widetilde w_l(\alpha,\beta)\partial_\beta^\alpha\mathbf P^\perp f^\varepsilon\right\|^2 \Bigg\} \nonumber\\
\lesssim\;&
\frac1{\varepsilon^2}\sum_{|\alpha|\leq N-1}
\left\|\widetilde w_l(\alpha,0)\partial^\alpha\mathbf P^\perp f^\varepsilon\right\|_D^2+\mathcal D_N(t)
+\Big\{\delta_0+\mathcal E_N(t)+\sqrt{\widetilde{\mathcal E}_{N,l}(t)}\Big\}\widetilde{\mathcal D}_{N,l}(t),
\end{align}
where $\widetilde{C}_j>0$ are fixed constants with $\widetilde{C}_j\gg \widetilde{C}_{j+1}$.
\end{lemma}

\begin{proof}
The proof is similar to that of Lemma \ref{weighted 1} and is omitted for brevity.
\end{proof}
\smallskip

\subsection{Energy Estimate in Negative Sobolev Spaces}
\hspace*{\fill}

This subsection is devoted to the energy estimates in negative Sobolev spaces.

\begin{lemma}\label{negative sobolev lemma}
Under the assumptions of Lemma \ref{ss:unweighted}, there holds
\begin{align}\label{negative sobolev estimate}
&\frac12\frac{\d}{\d t}\|\Lambda^{-\varrho}(f^\varepsilon,E^\varepsilon,B^\varepsilon)\|^2
+\frac{\lambda}{\varepsilon^2}\|\Lambda^{-\varrho}\mathbf P^\perp f^\varepsilon\|_D^2
\lesssim
\sqrt{\overline{\mathcal E}_{N,l}(t)}\overline{\mathcal D}_{N,l}(t).
\end{align}
\end{lemma}

\begin{proof}
Applying $\Lambda^{-\varrho}$ to the first equation in \eqref{rVMB}, taking the inner product with $\Lambda^{-\varrho} f^\varepsilon$, and using the Maxwell equations in \eqref{rVMB}, integration by parts and \eqref{L-coercive-unweighted}, we get
\begin{align}\label{negative:basic-identity}
&\frac12\frac{\d}{\d t}\|\Lambda^{-\varrho}(f^\varepsilon,E^\varepsilon,B^\varepsilon)\|^2
+\frac{\lambda}{\varepsilon^2}\|\Lambda^{-\varrho}\mathbf P^\perp f^\varepsilon\|_{D}^2\nonumber\\
\leq\;&-\Big(\Lambda^{-\varrho}\big((v\times B^\varepsilon)\cdot\nabla_v f^\varepsilon\big),\Lambda^{-\varrho}f^\varepsilon\Big)
-\varepsilon\Big(\Lambda^{-\varrho}(E^\varepsilon\cdot \nabla_v f^\varepsilon),\Lambda^{-\varrho}f^\varepsilon\Big)\\
&+\frac{\varepsilon}{2}\Big(\Lambda^{-\varrho}(E^\varepsilon\cdot v f^\varepsilon),\Lambda^{-\varrho}f^\varepsilon\Big)
+\frac1\varepsilon\big(\Lambda^{-\varrho}\Gamma(f^\varepsilon,f^\varepsilon),\Lambda^{-\varrho}\mathbf P^\perp f^\varepsilon\big)
:=\sum_{i=1}^{4}J_i.\nonumber
\end{align}

For the term $J_1$, arguing as in the derivation of \eqref{B important1} and \eqref{B important2}, one has
\begin{align*}
J_1=
\big(\Lambda^{-\varrho}(b^\varepsilon\times B^\varepsilon),\Lambda^{-\varrho}b^\varepsilon\big)
-\Big(\Lambda^{-\varrho}\big((v\times B^\varepsilon)\cdot\nabla_v\mathbf P^\perp f^\varepsilon\big),
\Lambda^{-\varrho}\mathbf P^\perp f^\varepsilon\Big):=J_{1,1}+J_{1,2}.
\end{align*}
For the term $J_{1,1}$, observing that
\begin{align*}
\big(\Lambda^{-\varrho}(b^\varepsilon\times B^\varepsilon),\Lambda^{-\varrho}b^\varepsilon\big)
=\big(b^\varepsilon\times B^\varepsilon,\Lambda^{-2\varrho}b^\varepsilon\big)
=\big(b^\varepsilon,  B^\varepsilon \times\Lambda^{-2\varrho}b^\varepsilon\big)
=\big(\Lambda^{1-\varrho}b^\varepsilon,\Lambda^{\varrho-1}(B^\varepsilon\times\Lambda^{-2\varrho}b^\varepsilon)\big),
\end{align*}
we have from Lemma \ref{kato-ponce} that
\begin{align*}
J_{1,1}
\lesssim
\|\Lambda^{1-\varrho}b^\varepsilon\|
\Big\{\|\Lambda^{\varrho-1}B^\varepsilon\|\|\Lambda^{-2\varrho}b^\varepsilon\|_{L^\infty}
+\|B^\varepsilon\|_{L^{\frac{6}{5-2\varrho}}}\|\Lambda^{-1-\varrho}b^\varepsilon\|_{L^{\frac{3}{\varrho-1}}}\Big\}.
\end{align*}
By Lemmas \ref{sobolev interpolation} and \ref{negative embedding theorem2}, together with $1<\varrho<3/2$, we obtain
\begin{align*}
&\|\Lambda^{-2\varrho}b^\varepsilon\|_{L^\infty}
\lesssim
\|\Lambda^{-\varrho}b^\varepsilon\|^{\varrho-\frac12}
\|\Lambda^{1-\varrho}b^\varepsilon\|^{\frac32-\varrho},\\
&\|\Lambda^{-1-\varrho}b^\varepsilon\|_{L^{\frac{3}{\varrho-1}}}
\lesssim
\|\Lambda^{\frac32-2\varrho}b^\varepsilon\|
\lesssim
\|\Lambda^{-\varrho}b^\varepsilon\|^{\varrho-\frac12}
\|\Lambda^{1-\varrho}b^\varepsilon\|^{\frac32-\varrho},\\
&\|B^\varepsilon\|_{L^{\frac{6}{5-2\varrho}}}
\lesssim
\|\Lambda^{\varrho-1}B^\varepsilon\|
\lesssim
\|\Lambda^{-\varrho}B^\varepsilon\|^{\frac32-\varrho}
\|\Lambda^{2-\varrho}B^\varepsilon\|^{\varrho-\frac12}.
\end{align*}
Consequently, combining these estimates and applying Young's inequality, we have
\begin{align*}
J_{1,1}
\lesssim\;&
\|\Lambda^{1-\varrho}b^\varepsilon\|
\Big(\|\Lambda^{-\varrho}B^\varepsilon\|\|\Lambda^{1-\varrho}b^\varepsilon\|\Big)^{\frac32-\varrho}
\Big(\|\Lambda^{2-\varrho}B^\varepsilon\|\|\Lambda^{-\varrho}b^\varepsilon\|\Big)^{\varrho-\frac12}\\
\lesssim\;&
\|\Lambda^{-\varrho}B^\varepsilon\|\|\Lambda^{1-\varrho}b^\varepsilon\|^2
+\|\Lambda^{2-\varrho}B^\varepsilon\|\|\Lambda^{-\varrho}b^\varepsilon\|\|\Lambda^{1-\varrho}b^\varepsilon\|\\
\lesssim\;&
\sqrt{\overline{\mathcal E}_{N,l}(t)}\overline{\mathcal D}_{N,l}(t).
\end{align*}
For the term $J_{1,2}$, Lemma \ref{negative embedding theorem}, H\"older's inequality, and \eqref{norm inequality} yield
\begin{align*}
&\left|\big(\Lambda^{-\varrho}
 ((v\times B^\varepsilon)\cdot\nabla_v\mathbf P^\perp f^\varepsilon),
 \Lambda^{-\varrho}\mathbf P^\perp f^\varepsilon\big)\right|\\
\lesssim\;&
\left\|B^\varepsilon\langle v\rangle^{1-\frac{\gamma+2s}{2}}
\nabla_v\mathbf P^\perp f^\varepsilon
\right\|_{L_x^{\frac6{3+2\varrho}}L_v^2}
\|\Lambda^{-\varrho}\mathbf P^\perp f^\varepsilon\|_D\\
\lesssim\;&
\|B^\varepsilon\|_{L^{\frac{3}{\varrho}}}
\left\|\langle v\rangle^{1-\frac{\gamma+2s}{2}}
\nabla_v\mathbf P^\perp f^\varepsilon\right\|
\|\Lambda^{-\varrho}\mathbf P^\perp f^\varepsilon\|_D\\
\lesssim\;&
\sqrt{\mathcal E_N(t)}\overline{\mathcal D}_{N,l}(t),
\end{align*}
where we have used $\|B^\varepsilon\|_{L^{\frac{3}{\varrho}}}\lesssim\|B^\varepsilon\|_{H^1}$ and
\begin{align*}
\left\|\langle v\rangle^{1-\frac{\gamma+2s}{2}}\nabla_v\mathbf P^\perp f^\varepsilon\right\|
\lesssim
\sum_{|\beta_1|=1}\|w_l(0,\beta_1)\partial_{\beta_1}\mathbf P^\perp f^\varepsilon\|_D
\lesssim \sqrt{\widetilde{\mathcal D}_{N,l}(t)}.
\end{align*}
Consequently, it follows that
\begin{align}\label{negative:J1-bound}
J_1\lesssim\sqrt{\overline{\mathcal E}_{N,l}(t)}\overline{\mathcal D}_{N,l}(t).
\end{align}

For the terms $J_2$ and $J_3$, a direct computation gives
\begin{align}\label{J23new}
&\mathbf P\bigg(\frac{1}{2} E^\varepsilon\cdot v f^\varepsilon
  -E^\varepsilon\cdot \nabla_v f^\varepsilon\bigg)
=\left\{a^\varepsilon E^\varepsilon\cdot v
 +\frac13(E^\varepsilon\cdot b^\varepsilon)(|v|^2-3)\right\}\mu^{1/2}, \nonumber\\
&\mathbf P\bigg(\frac{1}{2} E^\varepsilon\cdot v \mathbf P^\perp f^\varepsilon
  -E^\varepsilon\cdot \nabla_v \mathbf P^\perp f^\varepsilon\bigg)=0.
\end{align}
Therefore, we have
\begin{align}\label{J23identity}
J_2+J_3=\;&
\varepsilon\big(\Lambda^{-\varrho}(a^\varepsilon E^\varepsilon),\Lambda^{-\varrho}b^\varepsilon\big)
+2\varepsilon\big(\Lambda^{-\varrho}(E^\varepsilon\cdot b^\varepsilon),\Lambda^{-\varrho}c^\varepsilon\big) \nonumber\\
&+\varepsilon\bigg(\Lambda^{-\varrho}\bigg(\frac{1}{2}E^\varepsilon\cdot v\mathbf P^\perp f^\varepsilon
-E^\varepsilon\cdot\nabla_v\mathbf P^\perp f^\varepsilon\bigg),\Lambda^{-\varrho}\mathbf P^\perp f^\varepsilon\bigg) \nonumber\\
&+\varepsilon\bigg(\Lambda^{-\varrho}\mathbf P^\perp\bigg(\frac{1}{2}E^\varepsilon\cdot v\mathbf P f^\varepsilon
-E^\varepsilon\cdot\nabla_v\mathbf P f^\varepsilon\bigg),\Lambda^{-\varrho}\mathbf P^\perp f^\varepsilon\bigg).
\end{align}
Using Lemmas \ref{negative embedding theorem}, \ref{negative embedding theorem2} and \ref{sobolev interpolation}, H\"{o}lder's inequality, Young's inequality and the fact that $1-\varrho<3/4-\varrho/2<1$, we bound the first two terms on the right-hand side of \eqref{J23identity} by
\begin{align*}
&\|\Lambda^{-\varrho}(b^\varepsilon,c^\varepsilon)\|
\|E^\varepsilon\|_{L^{\frac{12}{3+2\varrho}}}
\|(a^\varepsilon,b^\varepsilon)\|_{L^{\frac{12}{3+2\varrho}}}\\
\lesssim\;&
\|\Lambda^{-\varrho}(b^\varepsilon,c^\varepsilon)\|
\|\Lambda^{\frac34-\frac{\varrho}{2}}E^\varepsilon\|
\|\Lambda^{\frac34-\frac{\varrho}{2}}
(a^\varepsilon,b^\varepsilon)\|\\
\lesssim\;&
\|\Lambda^{-\varrho}(b^\varepsilon,c^\varepsilon)\|
\Big\{\|\Lambda^{1-\varrho}(E^\varepsilon,a^\varepsilon,b^\varepsilon)\|^2
+\|\nabla_x(E^\varepsilon,a^\varepsilon,b^\varepsilon)\|^2\Big\}\\
\lesssim\;&
\sqrt{\overline{\mathcal E}_{N,l}(t)}\overline{\mathcal D}_{N,l}(t).
\end{align*}
Arguing as in the estimate of $J_{1,2}$, we bound the third term on the right-hand side of \eqref{J23identity} by $\sqrt{\mathcal E_N(t)}\overline{\mathcal D}_{N,l}(t)$.
Furthermore, the last term on the right-hand side of \eqref{J23identity} can be controlled by
\begin{align*}
&\|E^\varepsilon\|_{L^{\frac{3}{\varrho}}}
\|(a^\varepsilon,b^\varepsilon,c^\varepsilon)\|
\|\Lambda^{-\varrho}\mathbf P^\perp f^\varepsilon\|_D\\
\lesssim\;& \|E^\varepsilon\|_{L^{\frac{3}{\varrho}}}
\Big\{\|\Lambda^{1-\varrho}(a^\varepsilon,b^\varepsilon,c^\varepsilon)\|
+\|\nabla_x(a^\varepsilon,b^\varepsilon,c^\varepsilon)\|\Big\}
\|\Lambda^{-\varrho}\mathbf P^\perp f^\varepsilon\|_D\\
\lesssim\;&\sqrt{\mathcal E_N(t)}\overline{\mathcal D}_{N,l}(t).
\end{align*}
Collecting these estimates, we obtain
\begin{align}\label{negative:J23-bound}
J_2+J_3\lesssim\sqrt{\overline{\mathcal E}_{N,l}(t)}\overline{\mathcal D}_{N,l}(t).
\end{align}

For the term $J_4$, the collision invariant property, the Riesz potential representation of $\Lambda^{-\varrho}$, \eqref{Gamma-unweighted}, Lemmas \ref{negative embedding theorem} and \ref{negative embedding theorem2}, and H\"older's inequality yield
\begin{align}\label{negative:J4-bound}
J_4
\lesssim\;&
\frac1\varepsilon
\left\||f^\varepsilon|_{L_v^2}|f^\varepsilon|_D\right\|_{L_x^{\frac6{3+2\varrho}}}
\left\|\Lambda^{-\varrho}\mathbf P^\perp f^\varepsilon\right\|_D \nonumber\\
\lesssim\;&
\frac1\varepsilon\|f^\varepsilon\|_{L_x^{\frac{3}{\varrho}}L_v^2}\|f^\varepsilon\|_D
\left\|\Lambda^{-\varrho}\mathbf P^\perp f^\varepsilon\right\|_D \nonumber \\
\lesssim\;&
\sqrt{\mathcal E_N(t)} \overline{\mathcal D}_{N,l}(t).
\end{align}
In the last inequality, we have used
\begin{align*}
\|f^\varepsilon\|_{D}
\lesssim\|\mathbf P f^\varepsilon\|+ \|\mathbf P^\perp f^\varepsilon\|_{D}
\lesssim
\|\Lambda^{1-\varrho}\mathbf Pf^\varepsilon\|+\|\nabla_x\mathbf Pf^\varepsilon\|+ \|\mathbf P^\perp f^\varepsilon\|_{D}
\lesssim \sqrt{\overline{\mathcal D}_{N,l}(t)}.
\end{align*}

Substituting \eqref{negative:J1-bound}, \eqref{negative:J23-bound} and \eqref{negative:J4-bound} into \eqref{negative:basic-identity}, we obtain  \eqref{negative sobolev estimate}. This completes the proof.
\end{proof}
\smallskip

We next establish the macroscopic dissipation estimate in negative Sobolev spaces.
\begin{lemma}\label{negative:macro-lemma}
Under the assumptions of Lemma \ref{ss:unweighted}, there exists an interactive energy functional
\begin{align}\label{negative:macro-energy}
\mathcal E_{int}^{-\varrho}(t)
:=\;&
\overline C_4\Bigg\{\varepsilon\overline C_3\Bigg[
\sum_{i,j=1}^3
\Big(\Lambda^{-\varrho}\Theta_{ij}(\mathbf P^\perp f^\varepsilon),
\Lambda^{-\varrho}(\partial_i b_j^\varepsilon+\partial_j b_i^\varepsilon)\Big)\nonumber\\
&\qquad+ \sum_{i=1}^3
\Big(\Lambda^{-\varrho}\Lambda_i(\mathbf P^\perp f^\varepsilon),\Lambda^{-\varrho}\partial_i c^\varepsilon\Big)
+\kappa_4\Big(\Lambda^{-\varrho}b^\varepsilon,\Lambda^{-\varrho}\nabla_xa^\varepsilon\Big)\Bigg]\nonumber\\
&\qquad-\kappa_5
\Big(\Lambda^{1-\varrho}E^\varepsilon,\Lambda^{1-\varrho}b^\varepsilon\Big)
-\kappa_6\Big(\Lambda^{1-\varrho}E^\varepsilon,\Lambda^{1-\varrho}\nabla_x\times B^\varepsilon\Big)\Bigg\},
\end{align}
such that the following estimate holds
\begin{align}\label{negative:macro-estimate}
&\frac{\d}{\d t}\mathcal E_{int}^{-\varrho}(t)
+\|\Lambda^{1-\varrho}\mathbf Pf^\varepsilon\|^2
+\varepsilon^2\|\Lambda^{-\varrho}a^\varepsilon\|^2
+\|\Lambda^{1-\varrho}E^\varepsilon\|^2
+\|\Lambda^{2-\varrho}B^\varepsilon\|^2\nonumber\\
\lesssim\;&
\mathcal D_N(t)+\frac1{\varepsilon^2}\|\Lambda^{-\varrho}\mathbf P^\perp f^\varepsilon\|_D^2+\mathcal E_N(t)\overline{\mathcal D}_{N,l}(t).
\end{align}
Here, $0<\kappa_6\ll\kappa_5\ll\kappa_4\ll1$ are suitably chosen fixed constants and $\overline C_3,\overline C_4 >0$ are sufficiently large fixed constants.
\end{lemma}

\begin{proof}
Arguing as in the proof of Lemma \ref{macroscopic estimate} and using Lemmas \ref{negative embedding theorem} and \ref{negative embedding theorem2} to estimate the nonlinear terms in negative Sobolev norms, we obtain \eqref{negative:macro-estimate}.
We omit the details.
\end{proof}
\smallskip

\subsection{Uniform Energy Estimates}
\hspace*{\fill}

In this subsection, we combine the estimates established in Lemmas \ref{macroscopic estimate}--\ref{negative:macro-lemma} to derive uniform a priori energy estimates with respect to $\varepsilon\in(0,1]$.

\begin{proposition}\label{uniform-unweighted}
Under the a priori assumption \eqref{priori-assumptions-decay}, there exists an energy functional $\mathcal E_N(t)$ satisfying \eqref{without weight energy functional} such that, for any $t\in [0,T]$,
\begin{align}\label{uniform-unweighted-estimate}
\frac{\d}{\d t}\mathcal E_N(t)+\lambda\mathcal D_N(t)
\lesssim
\delta_0\widetilde{\mathcal D}_{N,l}(t).
\end{align}
\end{proposition}

\begin{proof}
Under the a priori assumption \eqref{priori-assumptions-decay}, the combination $\kappa_4 \times \eqref{macro estimate} + \eqref{ss:unweighted-estimate}$ gives \eqref{uniform-unweighted-estimate} for a fixed constant $0 < \kappa_4 \ll 1$. This concludes the proof of Proposition \ref{uniform-unweighted}.
\end{proof}
\smallskip

\begin{proposition}\label{uniform-weighted}
Under the a priori assumption \eqref{priori-assumptions-decay}, there exists an energy functional $\widetilde{\mathcal E}_{N,l}(t)$ satisfying \eqref{energy functional} such that, for any $t\in [0,T]$,
\begin{align}\label{uniform-weighted-estimate}
\frac{\d}{\d t}\widetilde{\mathcal E}_{N,l}(t)+\lambda\widetilde{\mathcal D}_{N,l}(t)\leq 0.
\end{align}
\end{proposition}

\begin{proof}
First, multiplying \eqref{uniform-unweighted-estimate} by $(1+t)^{-\varepsilon_0}$, we obtain
\begin{align}\label{uniform-time-weighted-estimate}
\frac{\d}{\d t}
\Big\{(1+t)^{-\varepsilon_0}\mathcal E_N(t)\Big\}
+\frac{\varepsilon_0}{(1+t)^{1+\varepsilon_0}}\mathcal E_N(t)
+\lambda(1+t)^{-\varepsilon_0}\mathcal D_N(t)
\lesssim
\delta_0\widetilde{\mathcal D}_{N,l}(t).
\end{align}
Next, we choose fixed constants $\overline C_5\gg\overline C_6\gg1$.
Under the a priori assumption \eqref{priori-assumptions-decay}, the combination
$\overline{C}_5 \times \left\{ \eqref{weighted estimate2} + \eqref{weighted estimate3} \right\} + \overline{C}_6 \times \eqref{weighted estimate4} + \eqref{weighted estimate1}$ yields
\begin{align}\label{uniform-weighted-intermediate}
\frac{\d}{\d t}
\widetilde{\mathcal E}_{N,l}^{h}(t)+\lambda\widetilde{\mathcal D}_{N,l}^{h}(t)
\lesssim
\mathcal D_N(t)+\delta_0\widetilde{\mathcal D}_{N,l}(t)
+\frac{\varepsilon^2}{(1+t)^{1+\varepsilon_0}}\sum_{|\alpha|=N}\|\partial^\alpha E^\varepsilon\|^2,
\end{align}
where
\begin{align*}
\widetilde{\mathcal E}_{N,l}^{h}(t)
:=\;&
\overline C_5\Bigg\{\sum_{|\alpha|\leq N-1}
\left\|\widetilde w_l(\alpha,0)\partial^\alpha\mathbf P^\perp f^\varepsilon\right\|^2
+\varepsilon\sum_{|\alpha|=N}
\left\|\widetilde w_l(\alpha,0)\partial^\alpha f^\varepsilon\right\|^2\Bigg\}\\
&+\overline C_6\sum_{j=1}^{N-1}\widetilde C_j
\sum_{\substack{|\alpha|+|\beta|\leq N-1\\|\beta|=j}}
\left\|\widetilde w_l(\alpha,\beta)\partial_\beta^\alpha\mathbf P^\perp f^\varepsilon\right\|^2
+\sum_{j=1}^{N}C_j\sum_{\substack{|\alpha|+|\beta|=N\\|\beta|=j}}
\left\|\widetilde w_l(\alpha,\beta)\partial_\beta^\alpha\mathbf P^\perp f^\varepsilon\right\|^2,
\end{align*}
and $\widetilde{\mathcal D}_{N,l}^{h}(t)$ consists of the corresponding weighted dissipation terms, so that
\begin{align*}
\mathcal D_N(t)+\widetilde{\mathcal D}_{N,l}^{h}(t)\sim\widetilde{\mathcal D}_{N,l}(t).
\end{align*}

Finally, since $\sum_{|\alpha|=N}\|\partial^\alpha E^\varepsilon\|^2\lesssim\mathcal E_N(t)$, we choose a sufficiently large fixed constant $\overline C_7>0$.
The combination
$\overline{C}_7 \times \left\{ \eqref{uniform-unweighted-estimate} + \eqref{uniform-time-weighted-estimate} \right\} + \eqref{uniform-weighted-intermediate}$
then yields the desired estimate \eqref{uniform-weighted-estimate} for sufficiently small $\delta_0>0$, with
\begin{align*}
\widetilde{\mathcal E}_{N,l}(t):=\overline C_7\big\{1+(1+t)^{-\varepsilon_0}\big\}\mathcal E_N(t)+\widetilde{\mathcal E}_{N,l}^{h}(t)
\end{align*}
satisfying \eqref{energy functional}.
This completes the proof of Proposition \ref{uniform-weighted}.
\end{proof}

\begin{proposition}\label{uniform-negative}
Under the a priori assumption \eqref{priori-assumptions-decay}, there exists an energy functional $\overline{\mathcal E}_{N,l}(t)$ satisfying \eqref{negative sobolev energy} such that, for any $t\in[0,T]$,
\begin{align}\label{uniform-negative-estimate}
\frac{\d}{\d t}\overline{\mathcal E}_{N,l}(t)
+\lambda\overline{\mathcal D}_{N,l}(t)\leq0.
\end{align}
\end{proposition}

\begin{proof}
Under the a priori assumption \eqref{priori-assumptions-decay}, a suitable linear combination of \eqref{negative sobolev estimate}, \eqref{uniform-weighted-estimate} and \eqref{negative:macro-estimate} yields \eqref{uniform-negative-estimate}.
This completes the proof.
\end{proof}

\section{Global Existence}\label{Global Existence}

In this section, we prove Theorem \ref{mainth1} by closing the a priori assumption \eqref{priori-assumptions-decay}.

\subsection{Time-Decay Estimates}
\hspace*{\fill}

Following the approach of \cite{GW2012CPDE}, we first establish a refined energy estimate for pure spatial derivatives.

\begin{lemma}\label{decay:refined-energy-lemma}
Under the a priori assumption \eqref{priori-assumptions-decay}, there exist an energy functional $\mathcal E_{N_0}^k(t)$ and the corresponding dissipation rate functional
$\mathcal D_{N_0}^k(t)$ satisfying \eqref{low k energy} and \eqref{low k dissipation}, respectively, such that, for $k=0,1,2$ and $t\in[0,T]$,
\begin{align}\label{decay:refined-energy}
\frac{\d}{\d t}\mathcal E_{N_0}^k(t)+\lambda\mathcal D_{N_0}^k(t)
\lesssim\;&
\overline{\mathcal D}_{N,l}(t)\mathcal E_{N_0}^k(t)
+\|E^\varepsilon\|_{L^\infty}^4\widetilde{\mathcal D}_{N,l}(t).
\end{align}
\end{lemma}

\begin{proof}
Applying $\partial^\alpha$ with $k\leq|\alpha|\leq N_0$ to the first equation in \eqref{rVMB}, taking the inner product with $\partial^\alpha f^\varepsilon$, and using the Maxwell equations in \eqref{rVMB} and \eqref{L-coercive-unweighted}, we obtain
\begin{align}\label{decay:kinetic-basic}
&\frac12\frac{\d}{\d t}
\left\|\partial^\alpha(f^\varepsilon,E^\varepsilon,B^\varepsilon)\right\|^2
+\frac{\lambda}{\varepsilon^2}\left\|\partial^\alpha\mathbf P^\perp f^\varepsilon\right\|_D^2\nonumber\\
\leq\;&
-\Big(\partial^\alpha\big((v\times B^\varepsilon)\cdot\nabla_v f^\varepsilon\big),\partial^\alpha f^\varepsilon\Big)
-\varepsilon\Big(\partial^\alpha(E^\varepsilon\cdot\nabla_v f^\varepsilon),\partial^\alpha f^\varepsilon\Big) \nonumber\\
&+\frac{\varepsilon}{2}\Big(\partial^\alpha (E^\varepsilon\cdot vf^\varepsilon), \partial^\alpha f^\varepsilon\Big)
+\frac1\varepsilon\Big(\partial^\alpha\Gamma(f^\varepsilon,f^\varepsilon),\partial^\alpha f^\varepsilon\Big)
:=\sum_{i=1}^{4}J_i.
\end{align}

For the term $J_1$, arguing as in the derivation of \eqref{B important1} and \eqref{B important2}, we obtain
\begin{align}\label{decay:magnetic-decomposition}
J_1
=\;&
\sum_{0<\alpha_1\leq\alpha}C_\alpha^{\alpha_1}
\Big( \partial^{\alpha-\alpha_1}b^\varepsilon\times\partial^{\alpha_1}B^\varepsilon,
\partial^\alpha b^\varepsilon \Big)\nonumber\\
&-\sum_{0<\alpha_1\leq\alpha}C_\alpha^{\alpha_1}
\Big( (v\times\partial^{\alpha_1}B^\varepsilon)\cdot\nabla_v\partial^{\alpha-\alpha_1}\mathbf P^\perp f^\varepsilon,
\partial^\alpha\mathbf P^\perp f^\varepsilon\Big):=J_{1,1}+J_{1,2}.
\end{align}
For the term $J_{1,1}$, we distinguish two cases.
For the case $k=0,1$,  H\"older's inequality, \eqref{Sobolev-ineq}, and the fact that $k\leq1\leq|\alpha_1|$ give
\begin{align*}
J_{1,1}
\lesssim\;&\sum_{1\leq |\alpha_1|\leq |\alpha| \leq N_0-1}
\|\partial^{\alpha_1}B^\varepsilon\|
\|\partial^{\alpha-\alpha_1}b^\varepsilon\|_{L^3}
\|\partial^\alpha b^\varepsilon\|_{L^6}\\
&+\sum_{1\leq |\alpha_1|\leq |\alpha|=N_0}
\|\partial^{\alpha_1}B^\varepsilon\|
\|\partial^{\alpha-\alpha_1}b^\varepsilon\|_{L^\infty}
\|\partial^\alpha b^\varepsilon\|\\
\lesssim\;&
\sqrt{\mathcal E_{N_0}^k(t)}\sqrt{\overline{\mathcal D}_{N,l}(t)}\sqrt{\mathcal D_{N_0}^k(t)}.
\end{align*}
For the case $k=2$, the terms with $|\alpha_1|\geq2$ can be estimated in the same manner as in the case $k=0,1$.
Each term with $|\alpha_1|=1$ and $3\leq|\alpha|\leq N_0$ is bounded by
\begin{align*}
\|\nabla_xB^\varepsilon\|_{L^3}
\|\partial^{\alpha-\alpha_1}b^\varepsilon\|_{L^6}
\|\partial^\alpha b^\varepsilon\|
\lesssim \sqrt{\mathcal E_N(t)}\mathcal D_{N_0}^2(t).
\end{align*}
By Lemma \ref{sobolev interpolation} and \eqref{Sobolev-ineq}, the remaining terms with $|\alpha_1|=1$ and $|\alpha|=2$ can be bounded by
\begin{align*}
\|\nabla_xB^\varepsilon\|_{L^3}\|\nabla_xb^\varepsilon\|\|\nabla_x^2b^\varepsilon\|_{L^6}
\lesssim
\|\Lambda^{-1}B^\varepsilon\|^{\frac{1}{2}}\|\nabla_x^4B^\varepsilon\|^{\frac{1}{2}}\|\Lambda^{-1}b^\varepsilon\|^{\frac{1}{2}}
\|\nabla_x^3b^\varepsilon\|^{\frac{3}{2}}
\lesssim
\sqrt{\overline{\mathcal E}_{N,l}(t)}\mathcal D_{N_0}^2(t),
\end{align*}
where we have used $\varrho>1$, $N_0\geq5$ and \eqref{low k dissipation}.
Combining the above estimates and applying the Cauchy--Schwarz inequality and \eqref{priori-assumptions-decay}, we conclude that, for any $\eta>0$,
\begin{align*}
J_{1,1}\lesssim \bigg\{ \eta+\sqrt{\overline{\mathcal E}_{N,l}(t)}\chi_{k=2}\bigg\}
\mathcal D_{N_0}^k(t)+\overline{\mathcal D}_{N,l}(t)\mathcal E_{N_0}^k(t), \qquad k=0,1,2.
\end{align*}
For the term $J_{1,2}$, H\"older's inequality, \eqref{Sobolev-ineq} and \eqref{norm inequality} give
\begin{align*}
J_{1,2}
\lesssim\;&
\sum_{|\alpha_1|=1}\|\partial^{\alpha_1}B^\varepsilon\|_{L^6}
\left\|\langle v\rangle^{1-\frac{\gamma+2s}{2}}\nabla_v\partial^{\alpha-\alpha_1}\mathbf P^\perp f^\varepsilon\right\|_{L_x^3L_v^2}
\left\|\partial^\alpha\mathbf P^\perp f^\varepsilon\right\|_D\\
&+\sum_{2\leq |\alpha_1|\leq N_0}\|\partial^{\alpha_1}B^\varepsilon\|
\left\|\langle v\rangle^{1-\frac{\gamma+2s}{2}}\nabla_v\partial^{\alpha-\alpha_1}\mathbf P^\perp f^\varepsilon\right\|_{L_x^\infty L_v^2}
\left\|\partial^\alpha\mathbf P^\perp f^\varepsilon\right\|_D\\
\lesssim\;&
\sqrt{\mathcal E_{N_0}^k(t)}\sqrt{\widetilde{\mathcal D}_{N,l}(t)}\sqrt{\mathcal D_{N_0}^k(t)}
\lesssim
\eta\mathcal D_{N_0}^k(t)+\widetilde{\mathcal D}_{N,l}(t)\mathcal E_{N_0}^k(t),
\qquad k=0,1,2.
\end{align*}
Here, we have used
\begin{equation*}
\begin{aligned}
|\alpha-\alpha_1+\alpha'|+1&\leq N_0+1, &\quad& \text{for}\;\; |\alpha_1|=1 \text{ and } |\alpha'|\leq 1,\\
|\alpha-\alpha_1+\alpha'|+1&\leq N_0+1, &\quad& \text{for}\;\; |\alpha_1|\geq 2 \text{ and } |\alpha'|\leq 2.
\end{aligned}
\end{equation*}
Combining the estimates for $J_{1,1}$ and $J_{1,2}$ with $k\leq|\alpha|\leq N_0$, and using \eqref{priori-assumptions-decay}, we obtain
\begin{align}\label{decay:J1-estimate}
J_1
\lesssim \bigg\{\eta+\sqrt{\overline{\mathcal E}_{N,l}(t)}\bigg\}\mathcal D_{N_0}^k(t)
+\overline{\mathcal D}_{N,l}(t)\mathcal E_{N_0}^k(t),
\qquad k=0,1,2.
\end{align}

We next estimate $J_2+J_3$. It follows from \eqref{J23new} that
\begin{align}\label{decay:electric-decomposition}
J_2+J_3
=\;&
\varepsilon\Big(\partial^\alpha(a^\varepsilon E^\varepsilon),\partial^\alpha b^\varepsilon\Big)
+2\varepsilon\Big(\partial^\alpha(E^\varepsilon\cdot b^\varepsilon),\partial^\alpha c^\varepsilon\Big)\nonumber\\
&+\varepsilon\Big(\partial^\alpha\Big(\frac12E^\varepsilon\cdot v\mathbf P^\perp f^\varepsilon
-E^\varepsilon\cdot\nabla_v\mathbf P^\perp f^\varepsilon\Big),
\partial^\alpha\mathbf P^\perp f^\varepsilon\Big)\nonumber\\
&+\varepsilon\Big(\partial^\alpha\mathbf P^\perp\Big(\frac12E^\varepsilon\cdot v\mathbf Pf^\varepsilon
-E^\varepsilon\cdot\nabla_v\mathbf Pf^\varepsilon\Big),\partial^\alpha\mathbf P^\perp f^\varepsilon\Big).
\end{align}
Arguing as in the estimates of $J_{1,1}$ and $J_{1,2}$ with $k\leq|\alpha|\leq N_0$, we obtain, for any $\eta>0$,
\begin{align*}
J_2+J_3\lesssim \|E^\varepsilon\|_{L^\infty}
\left\|\langle v\rangle^{\frac12}\partial^\alpha\mathbf P^\perp f^\varepsilon\right\|^2
+ \bigg\{\eta+\sqrt{\overline{\mathcal E}_{N,l}(t)}\bigg\} \mathcal D_{N_0}^k(t)
+\overline{\mathcal D}_{N,l}(t)\mathcal E_{N_0}^k(t).
\end{align*}
For the first term on the right-hand side of the above inequality, velocity-weight interpolation, \eqref{norm inequality}, and the Cauchy--Schwarz inequality give
\begin{align*}
\|E^\varepsilon\|_{L^\infty}
\left\|\langle v\rangle^{\frac12}\partial^\alpha\mathbf P^\perp f^\varepsilon\right\|^2
\lesssim\;&
\|E^\varepsilon\|_{L^\infty}
\left\|\partial^\alpha\mathbf P^\perp f^\varepsilon\right\|_D^{\frac32}
\left\|\langle v\rangle^{2-\frac{3(\gamma+2s)}2}\partial^\alpha\mathbf P^\perp f^\varepsilon\right\|^{\frac12}\\
\lesssim\;&
\|E^\varepsilon\|_{L^\infty}
\left\|\partial^\alpha\mathbf P^\perp f^\varepsilon\right\|_D^{\frac32}
\left\|w_l(\alpha,0)\partial^\alpha\mathbf P^\perp f^\varepsilon\right\|_D^{\frac12}\\
\leq\;&
\eta\left\|\partial^\alpha\mathbf P^\perp f^\varepsilon\right\|_D^2
+\|E^\varepsilon\|_{L^\infty}^4\widetilde{\mathcal D}_{N,l}(t).
\end{align*}
Combining these estimates for $k\leq|\alpha|\leq N_0$, we obtain
\begin{align}\label{decay:J23-estimate}
J_2+J_3
\lesssim \;&
\bigg\{\eta+\sqrt{\overline{\mathcal E}_{N,l}(t)}\bigg\} \mathcal D_{N_0}^k(t)
+\overline{\mathcal D}_{N,l}(t)\mathcal E_{N_0}^k(t)
+\|E^\varepsilon\|_{L^\infty}^4\widetilde{\mathcal D}_{N,l}(t),
\qquad k=0,1,2.
\end{align}

Finally, we estimate the term $J_4$. By the collision invariant property, \eqref{Gamma-unweighted}, and \eqref{Pf define}, we obtain
\begin{align*}
J_4=\;&
\frac1\varepsilon\Big(\partial^\alpha\Gamma(f^\varepsilon,f^\varepsilon),\partial^\alpha\mathbf P^\perp f^\varepsilon\Big)\\
\lesssim\;&
\frac1\varepsilon\sum_{\alpha_1\leq\alpha}
\int_{\mathbb R^3}\Big\{|\partial^{\alpha_1}\mathbf Pf^\varepsilon|_{L_v^2}|\partial^{\alpha-\alpha_1}\mathbf Pf^\varepsilon|_{L_v^2}
+|\partial^{\alpha_1}f^\varepsilon|_{L_v^2}|\partial^{\alpha-\alpha_1}\mathbf P^\perp f^\varepsilon|_D\Big\}
|\partial^\alpha\mathbf P^\perp f^\varepsilon|_D\d x.
\end{align*}
We distinguish two cases. For the case $k=0, 1$, since $\max\{ |\alpha_1|, |\alpha-\alpha_1|\} \geq k$, H\"older's inequality and \eqref{Sobolev-ineq} readily yield
\begin{align*}
J_4\lesssim\;&
\sqrt{\mathcal E_{N_0}^k(t)}\sqrt{\overline{\mathcal D}_{N,l}(t)}\sqrt{\mathcal D_{N_0}^k(t)}
+\varepsilon\sqrt{\mathcal E_N(t)}\mathcal D_{N_0}^k(t),
\end{align*}
where we have used $N_0+2\leq N$. For the case $k=2$, the same argument applies unless $|\alpha|=2$ and $|\alpha_1|=1$.
By \eqref{Gamma-unweighted}, H\"older's inequality, and \eqref{Sobolev-ineq}, the remaining terms can be bounded by
\begin{align*}
&\frac1\varepsilon\|\nabla_xf^\varepsilon\|_{L_x^6L_v^2}
\|\nabla_xf^\varepsilon\|_{L_x^3L_D^2}
\|\nabla_x^2\mathbf P^\perp f^\varepsilon\|_D\\
\lesssim\;&
\frac1\varepsilon\|\nabla_x^2f^\varepsilon\|
\Big\{\|\nabla_x\mathbf Pf^\varepsilon\|_{L_x^3L_v^2}+\|\nabla_x\mathbf P^\perp f^\varepsilon\|_{L_x^3L_D^2}\Big\}
\|\nabla_x^2\mathbf P^\perp f^\varepsilon\|_D\\
\lesssim\;&
\sqrt{\mathcal E_{N_0}^2(t)}\sqrt{\overline{\mathcal D}_{N,l}(t)}\sqrt{\mathcal D_{N_0}^2(t)}.
\end{align*}
Combining the above estimates for $k\leq|\alpha|\leq N_0$, and applying the Cauchy--Schwarz inequality, we obtain, for any $\eta>0$,
\begin{align}\label{decay:J4-estimate}
J_4\lesssim
\bigg\{\eta+\sqrt{\overline{\mathcal E}_{N,l}(t)}\bigg\}\mathcal D_{N_0}^k(t)
+\overline{\mathcal D}_{N,l}(t)\mathcal E_{N_0}^k(t), \qquad k=0,1,2.
\end{align}

Substituting \eqref{decay:J1-estimate}, \eqref{decay:J23-estimate}, and \eqref{decay:J4-estimate} into \eqref{decay:kinetic-basic} and summing over
$k\leq|\alpha|\leq N_0$, we obtain, for any $\eta>0$,
\begin{align}\label{decay:kinetic-summed}
&\frac12\frac{\d}{\d t}
\sum_{k\leq|\alpha|\leq N_0}
\left\|\partial^\alpha(f^\varepsilon,E^\varepsilon,B^\varepsilon)\right\|^2
+\frac{\lambda}{\varepsilon^2}\sum_{k\leq|\alpha|\leq N_0}\|\partial^\alpha\mathbf P^\perp f^\varepsilon\|_D^2\nonumber\\
\lesssim\;&
\bigg\{\eta+\sqrt{\overline{\mathcal E}_{N,l}(t)}\bigg\}\mathcal D_{N_0}^k(t)
+\overline{\mathcal D}_{N,l}(t)\mathcal E_{N_0}^k(t)+\|E^\varepsilon\|_{L^\infty}^4\widetilde{\mathcal D}_{N,l}(t),
\qquad k=0,1,2.
\end{align}

Furthermore, following the proof of Lemma \ref{macroscopic estimate} with suitable modifications to the nonlinear estimates, we obtain the following macroscopic estimate.
For $N_0\geq5$ and $k=0,1,2$, there exists an interactive energy functional $\mathcal E_{int}^{k\to N_0}(t)$ satisfying
\begin{align*}
\left|\mathcal E_{int}^{k\to N_0}(t)\right|
\lesssim\sum_{k\leq|\alpha|\leq N_0}
\left\|\partial^\alpha(f^\varepsilon,E^\varepsilon,B^\varepsilon)\right\|^2,
\end{align*}
such that
\begin{align}\label{decay:macroscopic-estimate}
&\frac{\d}{\d t}\mathcal E_{int}^{k\to N_0}(t)
+\sum_{k+1\leq|\alpha|\leq N_0}
\|\partial^\alpha\mathbf Pf^\varepsilon\|^2+\varepsilon^2\sum_{|\alpha|=k}\|\partial^\alpha a^\varepsilon\|^2\nonumber\\
&+\sum_{k+1\leq|\alpha|\leq N_0-1}
\|\partial^\alpha E^\varepsilon\|^2+\sum_{k+2\leq|\alpha|\leq N_0-1}\|\partial^\alpha B^\varepsilon\|^2\nonumber\\
\lesssim\;&
\frac1{\varepsilon^2}\sum_{k\leq|\alpha|\leq N_0}\|\partial^\alpha\mathbf P^\perp f^\varepsilon\|_D^2
+\sqrt{\overline{\mathcal E}_{N,l}(t)}\mathcal D_{N_0}^k(t)+\overline{\mathcal D}_{N,l}(t)\mathcal E_{N_0}^k(t).
\end{align}

Finally, for a sufficiently small fixed constant $\kappa_7>0$, define
\begin{align*}
\mathcal E_{N_0}^k(t) :=\frac12\sum_{k\leq|\alpha|\leq N_0}
\left\|\partial^\alpha(f^\varepsilon,E^\varepsilon,B^\varepsilon)\right\|^2+\kappa_7\mathcal E_{int}^{k\to N_0}(t),
\end{align*}
which satisfies \eqref{low k energy}.
Under the a priori assumption \eqref{priori-assumptions-decay}, the combination $\eqref{decay:kinetic-summed}+\kappa_7\times\eqref{decay:macroscopic-estimate}$ yields \eqref{decay:refined-energy}, provided that $\kappa_7$, $\eta$, and $\delta_0$ are chosen sufficiently small in turn.
This completes the proof.
\end{proof}

\begin{lemma}\label{decay:time-decay-lemma}
Under the a priori assumption \eqref{priori-assumptions-decay}, for $k=0,1,2$ and $t\in[0,T]$, we have
\begin{align}\label{decay:time-decay}
\mathcal E_{N_0}^k(t)\lesssim (1+t)^{-\frac{k+\varrho}{2}}\overline{\mathcal E}_{N,l}(0).
\end{align}
\end{lemma}

\begin{proof}
We first establish an interpolation estimate between $\mathcal E_{N_0}^k(t)$ and $\mathcal D_{N_0}^k(t)$.
By Lemma \ref{sobolev interpolation}, we have
\begin{align*}
\left\|\nabla_x^k(\mathbf Pf^\varepsilon,E^\varepsilon,B^\varepsilon)\right\|
\lesssim\left\|\Lambda^{-\varrho}(\mathbf Pf^\varepsilon,E^\varepsilon,B^\varepsilon)\right\|^{\frac{2}{k+\varrho+2}}
\left\|\nabla_x^{k+2}(\mathbf Pf^\varepsilon,E^\varepsilon,B^\varepsilon)\right\|^{\frac{k+\varrho}{k+\varrho+2}},
\end{align*}
and
\begin{align*}
\left\|\nabla_x^{N_0}(\mathbf Pf^\varepsilon,E^\varepsilon,B^\varepsilon)\right\|
\lesssim\left\|\nabla_x^{N_0-1}(\mathbf Pf^\varepsilon,E^\varepsilon,B^\varepsilon)\right\|^{\frac{k+\varrho}{k+\varrho+2}}
\left\|\Lambda^{N_0+\frac{k+\varrho}{2}}(\mathbf Pf^\varepsilon,E^\varepsilon,B^\varepsilon)\right\|^{\frac{2}{k+\varrho+2}}.
\end{align*}
Since $k+2\leq N_0-1$ and $N_0+(k+\varrho)/2\leq N$, the definitions of $\overline{\mathcal E}_{N,l}(t)$ and $\mathcal D_{N_0}^k(t)$ in \eqref{negative sobolev energy} and \eqref{low k dissipation}, respectively, together with the above two interpolation inequalities, yield
\begin{align*}
\sum_{k\leq|\alpha|\leq N_0}\left\|\partial^\alpha(\mathbf Pf^\varepsilon,E^\varepsilon,B^\varepsilon)\right\|^2
\lesssim
\Big\{\overline{\mathcal E}_{N,l}(t)\Big\}^{\frac{2}{k+\varrho+2}}\Big\{\mathcal D_{N_0}^k(t)\Big\}^{\frac{k+\varrho}{k+\varrho+2}}.
\end{align*}
Moreover, velocity-weight interpolation and \eqref{norm inequality} give, for $k\leq|\alpha|\leq N_0$,
\begin{align*}
\left\|\partial^\alpha\mathbf P^\perp f^\varepsilon\right\|
\lesssim\;&
\left\{\frac1\varepsilon\left\|\langle v\rangle^{\frac{\gamma+2s}{2}}\partial^\alpha\mathbf P^\perp f^\varepsilon\right\|\right\}^{\frac{k+\varrho}{k+\varrho+2}}
\left\|\langle v\rangle^{-\frac{(\gamma+2s)(k+\varrho)}4}\partial^\alpha\mathbf P^\perp f^\varepsilon\right\|^{\frac{2}{k+\varrho+2}}\\
\lesssim\;&
\Big\{\mathcal D_{N_0}^k(t)\Big\}^{\frac{k+\varrho}{2(k+\varrho+2)}}
\Big\{\overline{\mathcal E}_{N,l}(t)\Big\}^{\frac{1}{k+\varrho+2}}.
\end{align*}
Combining these estimates yields
\begin{align}\label{decay:energy-dissipation-interpolation}
\mathcal E_{N_0}^k(t)
\lesssim
\Big\{\overline{\mathcal E}_{N,l}(t)\Big\}^{\frac{2}{k+\varrho+2}}
\Big\{\mathcal D_{N_0}^k(t)\Big\}^{\frac{k+\varrho}{k+\varrho+2}}.
\end{align}

We next establish \eqref{decay:time-decay}.
Multiplying \eqref{decay:refined-energy} by $(1+t)^{1+\frac{k+\varrho}{2}}$ and using \eqref{decay:energy-dissipation-interpolation}, \eqref{priori-assumptions-decay}, and the Cauchy--Schwarz inequality, we obtain
\begin{align}
&\frac{\d}{\d t}
\Big\{(1+t)^{1+\frac{k+\varrho}{2}}\mathcal E_{N_0}^k(t)\Big\}
+\lambda(1+t)^{1+\frac{k+\varrho}{2}}\mathcal D_{N_0}^k(t) \nonumber\\
\lesssim\;&
(1+t)^{\frac{k+\varrho}{2}}\mathcal E_{N_0}^k(t)
+\overline{\mathcal D}_{N,l}(t)(1+t)^{1+\frac{k+\varrho}{2}}\mathcal E_{N_0}^k(t)
+(1+t)^{1+\frac{k+\varrho}{2}}\|E^\varepsilon(t)\|_{L^\infty}^4\widetilde{\mathcal D}_{N,l}(t)\nonumber\\
\lesssim\;&
\eta(1+t)^{1+\frac{k+\varrho}{2}}\mathcal D_{N_0}^k(t)+\overline{\mathcal E}_{N,l}(t)
+\overline{\mathcal D}_{N,l}(t)(1+t)^{1+\frac{k+\varrho}{2}}\mathcal E_{N_0}^k(t)+\delta_0^4\widetilde{\mathcal D}_{N,l}(t), \label{differential-inequality1}
\end{align}
where, in the last step, we have used \eqref{decay:energy-dissipation-interpolation} to obtain
\begin{align*}
(1+t)^{\frac{k+\varrho}{2}}\mathcal E_{N_0}^k(t)
\lesssim
\eta(1+t)^{1+\frac{k+\varrho}{2}}\mathcal D_{N_0}^k(t)
+\overline{\mathcal E}_{N,l}(t),
\end{align*}
and \eqref{priori-assumptions-decay} to deduce
\begin{align*}
(1+t)^{1+\frac{k+\varrho}{2}}\|E^\varepsilon(t)\|_{L^\infty}^4
\leq \delta_0^4(1+t)^{\frac{k-\varrho-1}{2}}
\leq\delta_0^4,
\end{align*}
since $k\leq2$ and $\varrho>1$.
Moreover, \eqref{uniform-negative-estimate} gives
\begin{align}\label{existence estimate 1}
\sup_{0\leq\tau\leq t}\overline{\mathcal E}_{N,l}(\tau)+\int_0^t\overline{\mathcal D}_{N,l}(\tau)\d\tau\lesssim \overline{\mathcal E}_{N,l}(0).
\end{align}
Integrating the differential inequality \eqref{differential-inequality1} in time and applying Gronwall's inequality along with \eqref{existence estimate 1} yields
\begin{align*}
(1+t)^{1+\frac{k+\varrho}{2}}\mathcal E_{N_0}^k(t)
\lesssim (1+t)\sup_{0\leq\tau\leq t}\overline{\mathcal E}_{N,l}(\tau),
\end{align*}
from which \eqref{decay:time-decay} follows immediately.
\end{proof}
\smallskip

\subsection{Proof of Global Existence}
\hspace*{\fill}

\begin{proof}[{\bf Proof of Theorem \ref{mainth1}.}]
By \eqref{existence estimate 1} and \eqref{decay:time-decay}, we have, for $0\leq t\leq T$,
\begin{align}\label{existence estimate 2}
\mathcal E_{N_0}^k(t)
\lesssim\overline{\mathcal E}_{N,l}(0)(1+t)^{-\frac{k+\varrho}{2}},\qquad k=0,1,2.
\end{align}
Consequently, \eqref{Sobolev-ineq} yields
\begin{align*}
\|E^\varepsilon(t)\|_{L^\infty}
&\lesssim
\|\nabla_xE^\varepsilon(t)\|^{\frac{1}{2}}
\|\nabla_x^2E^\varepsilon(t)\|^{\frac{1}{2}}
\lesssim
\sqrt{\overline{\mathcal E}_{N,l}(0)}(1+t)^{-\frac{3+2\varrho}{8}},\\
\|\nabla_x(E^\varepsilon,B^\varepsilon)(t)\|_{L^\infty}
&\lesssim
\|\nabla_x^2(E^\varepsilon,B^\varepsilon)(t)\|^{\frac{1}{2}}
\|\nabla_x^3(E^\varepsilon,B^\varepsilon)(t)\|^{\frac{1}{2}}
\lesssim \sqrt{\overline{\mathcal E}_{N,l}(0)} (1+t)^{-\frac{2+\varrho}{4}}.
\end{align*}
Together with \eqref{existence estimate 1}, these estimates verify the a priori assumption \eqref{priori-assumptions-decay}, provided that $\overline{\mathcal E}_{N,l}(0)$ is sufficiently small.

For each fixed $\varepsilon\in(0,1]$, a standard approximation argument based on the energy estimates above yields a unique local solution preserving the non-negativity of $F^\varepsilon$, and the uniform a priori estimate \eqref{existence estimate 1} allows this solution to be extended globally by a continuation argument.
Finally, the uniform energy estimate \eqref{thm1 estimate 1} and the time-decay estimate \eqref{thm1 estimate 2} follow from \eqref{existence estimate 1} and \eqref{existence estimate 2}, respectively.
This completes the proof.
\end{proof}
\medskip

\section{Limit to the Incompressible NSFM System}\label{Limit section}

In this section, we derive the incompressible NSFM system \eqref{INSFM limit} from the one-species VMB system \eqref{rVMB}.
We follow the moment method in \cite{WZGX2025SIMA}, while providing the additional arguments needed for the non-cutoff soft potentials and the acoustic components.

\begin{proof}[{\bf Proof of Theorem \ref{mainth2}.}]
By \eqref{thm1 estimate 1}, there exists a constant $C>0$, independent of $\varepsilon$, such that
\begin{align}\label{limit:uniform-bounds}
&\sup_{t\geq0}
\left\{\|f^\varepsilon(t)\|_{H_x^NL_v^2}^2+\|(E^\varepsilon,B^\varepsilon)(t)\|_{H^N_x}^2\right\}\leq C, \\
&\sum_{|\alpha|\leq N-1}\int_0^\infty\|\partial^\alpha\mathbf P^\perp f^\varepsilon(t)\|^2 \d t
\lesssim \sum_{|\alpha|\leq N-1}\int_0^\infty\|w_l(\alpha,0) \partial^\alpha\mathbf P^\perp f^\varepsilon(t)\|_D^2 \d t \leq C\varepsilon^2,
\end{align}
which implies that, up to a subsequence,
\begin{equation*}
\begin{aligned}
f^\varepsilon& \to  f
&&\quad \text{ weakly}\!-\!* \text{ in }L^\infty(\mathbb R^+;H_x^NL_v^2),\\
(E^\varepsilon,B^\varepsilon)& \to  (E,B)
&&\quad \text{ weakly}\!-\!* \text{ in }L^\infty(\mathbb R^+;H_x^N), \\
\mathbf{P}^\perp f^\varepsilon &\to 0
&&\quad \text{ strongly in }L^2(\mathbb R^+;H_x^{N-1} L_v^2)
\end{aligned}
\end{equation*}
as $\varepsilon \to 0$. From these convergence results, we immediately deduce that $\mathbf P^\perp f=0$.

As in \cite{WZGX2025SIMA}, we define the hydrodynamic variables
\begin{align}\label{rho-u-theta limit}
\rho^\varepsilon=\langle f^\varepsilon,\mu^{1/2}\rangle,\qquad
u^\varepsilon=\langle f^\varepsilon,v\mu^{1/2}\rangle,\qquad
\theta^\varepsilon=\frac13\langle f^\varepsilon,(|v|^2-3)\mu^{1/2}\rangle.
\end{align}
These variables satisfy $\rho^\varepsilon=a^\varepsilon$, $u^\varepsilon=b^\varepsilon$, and $\theta^\varepsilon=2c^\varepsilon$.
We also define the corresponding moments of the limiting function
\begin{align*}
\rho=\langle f,\mu^{1/2}\rangle,\qquad
u=\langle f,v\mu^{1/2}\rangle,\qquad
\theta=\frac13\langle f,(|v|^2-3)\mu^{1/2}\rangle.
\end{align*}
Thus we have
\begin{align*}
(\rho^\varepsilon,u^\varepsilon,\theta^\varepsilon)
\to(\rho,u,\theta)
\quad\text{ weakly}\!-\!* \text{ in }L^\infty(\mathbb R^+;H^N_x)
\end{align*}
as $\varepsilon \to 0$. Since $\mathbf P^\perp f=0$, the limiting function has the form
\begin{align*}
f=\left\{\rho+u\cdot v+\theta\frac{|v|^2-3}{2}\right\}\mu^{1/2}.
\end{align*}
Taking velocity moments of the first equation of \eqref{rVMB} and integrating by parts in $v$, we obtain
\begin{equation}
\left\{\begin{aligned}\label{limit:moment-equations}
&\partial_t\rho^\varepsilon+\frac1\varepsilon\nabla_x\cdot u^\varepsilon=0,\\
&\partial_tu^\varepsilon+\frac1\varepsilon\nabla_x(\rho^\varepsilon+\theta^\varepsilon)
+\frac1\varepsilon\nabla_x\cdot\langle \mathcal{A}\mu^{1/2},\mathbf P^\perp f^\varepsilon\rangle
=E^\varepsilon+u^\varepsilon\times B^\varepsilon+\varepsilon\rho^\varepsilon E^\varepsilon,\\
&\partial_t\theta^\varepsilon+\frac{2}{3\varepsilon}\nabla_x\cdot u^\varepsilon
+\frac{2}{3\varepsilon}\nabla_x\cdot\langle \mathcal{B}\mu^{1/2},\mathbf P^\perp f^\varepsilon\rangle
=\frac{2\varepsilon}{3}E^\varepsilon\cdot u^\varepsilon,
\end{aligned}\right.
\end{equation}
where
\begin{align*}
\mathcal{A}(v)=v\otimes v-\frac{|v|^2}{3}\mathbb I_{3\times3},
\qquad
\mathcal{B}(v)=v\left(\frac{|v|^2}{2}-\frac52\right).
\end{align*}
The functions $\widehat{\mathcal A}(v)$ and $\widehat{\mathcal B}(v)$ are uniquely determined by $L(\widehat{\mathcal A}\mu^{1/2})=\mathcal A\mu^{1/2}$ and $L(\widehat{\mathcal B}\mu^{1/2})=\mathcal B\mu^{1/2}$, respectively, with $\widehat{\mathcal A}\mu^{1/2},\widehat{\mathcal B}\mu^{1/2}\in\mathcal N^\perp(L)$.

Using the definitions of $\rho^\varepsilon$ and $u^\varepsilon$ in \eqref{rho-u-theta limit}, we rewrite the Maxwell equations in \eqref{rVMB} as
\begin{equation}\label{limit:maxwell-equations}
\left\{\begin{aligned}
&\partial_t E^\varepsilon-\nabla_x\times B^\varepsilon=-u^\varepsilon,
&&\nabla_x\cdot E^\varepsilon=\varepsilon\rho^\varepsilon,\\
&\partial_t B^\varepsilon+\nabla_x\times E^\varepsilon=0,
&&\nabla_x\cdot B^\varepsilon=0.
\end{aligned}\right.
\end{equation}
The moment equations \eqref{limit:moment-equations} and the Maxwell equations \eqref{limit:maxwell-equations} imply that the limiting fields satisfy
\begin{align*}
\nabla_x\cdot u=0,\qquad \rho+\theta=0,\qquad
\nabla_x\cdot E=\nabla_x\cdot B=0.
\end{align*}
Moreover, using \eqref{limit:uniform-bounds} and following the compactness argument in \cite{WZGX2025SIMA}, we obtain, along a subsequence, the following convergence for any $T>0$
\begin{align}\label{limit:strong-slow-components}
\left(\mathcal P u^\varepsilon,\frac35\theta^\varepsilon-\frac25\rho^\varepsilon,
E^\varepsilon,B^\varepsilon\right)
\to (u,\theta,E,B) \quad\text{strongly in } C([0,T];H^{N-1}_{loc})
\end{align}
as $\varepsilon\to0$.
We omit the details of this standard compactness argument and focus on the strong convergence of the acoustic and microscopic components.

From \eqref{limit:moment-equations}, we obtain
\begin{align}
&\partial_t\mathcal P u^\varepsilon
=\mathcal P\Big\{E^\varepsilon+u^\varepsilon\times B^\varepsilon
+\varepsilon\rho^\varepsilon E^\varepsilon
-\frac1\varepsilon\nabla_x\cdot\big\langle\mathcal A\mu^{1/2},\mathbf P^\perp f^\varepsilon\big\rangle\Big\}, \label{limit:slow-momentum}\\
&\partial_t\left(\frac35\theta^\varepsilon-\frac25\rho^\varepsilon\right)
=-\frac{2}{5\varepsilon}\nabla_x\cdot\big\langle\mathcal B\mu^{1/2},
\mathbf P^\perp f^\varepsilon\big\rangle+\frac{2\varepsilon}{5}E^\varepsilon\cdot u^\varepsilon, \label{limit:slow-temperature}\\
&\partial_t(\rho^\varepsilon+\theta^\varepsilon)+\frac{5}{3\varepsilon}\nabla_x\cdot\mathcal P^\perp u^\varepsilon
=-\frac{2}{3\varepsilon}\nabla_x\cdot\big\langle\mathcal B\mu^{1/2},
\mathbf P^\perp f^\varepsilon\big\rangle+\frac{2\varepsilon}{3}E^\varepsilon\cdot u^\varepsilon, \label{limit:acoustic-original-q}\\
&\partial_t\mathcal P^\perp u^\varepsilon+\frac1\varepsilon\nabla_x(\rho^\varepsilon+\theta^\varepsilon)
=\mathcal P^\perp\Big\{E^\varepsilon+u^\varepsilon\times B^\varepsilon+\varepsilon\rho^\varepsilon E^\varepsilon
-\frac1\varepsilon\nabla_x\cdot\big\langle\mathcal A\mu^{1/2},
\mathbf P^\perp f^\varepsilon\big\rangle\Big\},\label{limit:acoustic-original-w}
\end{align}
where $\mathcal P=I-\nabla_x\Delta_x^{-1}\nabla_x\cdot$ and $\mathcal P^\perp=I-\mathcal P$ denote the Leray projection onto divergence-free vector fields and its complementary projection, respectively.

To analyze the acoustic system \eqref{limit:acoustic-original-q}--\eqref{limit:acoustic-original-w}, we first expand the stress tensor and heat flux.
Taking the $L^2_v$ inner products of the first equation of \eqref{rVMB} with $\varepsilon\widehat{\mathcal A}\mu^{1/2}$ and $\varepsilon\widehat{\mathcal B}\mu^{1/2}$, respectively, we obtain
\begin{align}
\frac1\varepsilon
\big\langle\mathcal A\mu^{1/2},\mathbf P^\perp f^\varepsilon\big\rangle
=\;&u^\varepsilon\otimes u^\varepsilon-\frac{|u^\varepsilon|^2}{3}\mathbb I
-\nu\Big\{\nabla_xu^\varepsilon+(\nabla_xu^\varepsilon)^{\mathsf T}-\frac23(\nabla_x\cdot u^\varepsilon)\mathbb I\Big\}\nonumber \\
&-\varepsilon\partial_t\big\langle\widehat{\mathcal A}\mu^{1/2},\mathbf P^\perp f^\varepsilon\big\rangle+r_{\mathcal A}^\varepsilon,
\label{limit:stress-closure}\\
\frac1\varepsilon\big\langle\mathcal B\mu^{1/2},\mathbf P^\perp f^\varepsilon\big\rangle
=\;&\frac52u^\varepsilon\theta^\varepsilon-\frac{5\kappa}{2}\nabla_x\theta^\varepsilon
-\varepsilon\partial_t\big\langle\widehat{\mathcal B}\mu^{1/2},\mathbf P^\perp f^\varepsilon\big\rangle+r_{\mathcal B}^\varepsilon.
\label{limit:heat-closure}
\end{align}
Here, the viscosity coefficient $\nu$ and the heat conductivity coefficient $\kappa$ are defined by
\begin{align}\label{coefficient-define}
\nu=\frac1{10}\sum_{i,j=1}^3\big\langle\mathcal A_{ij}\mu^{1/2},\widehat{\mathcal A}_{ij}\mu^{1/2}\big\rangle, \qquad
\kappa=\frac2{15}\sum_{i=1}^3\big\langle\mathcal B_i\mu^{1/2},\widehat{\mathcal B}_i\mu^{1/2}\big\rangle,
\end{align}
and, for $\Xi\in\{\mathcal A,\mathcal B\}$, the remainder is given by
\begin{align}\label{r-define}
r_\Xi^\varepsilon=\Big\langle
&\Gamma(\mathbf Pf^\varepsilon,\mathbf P^\perp f^\varepsilon)
+\Gamma(\mathbf P^\perp f^\varepsilon,\mathbf Pf^\varepsilon)
+\Gamma(\mathbf P^\perp f^\varepsilon,\mathbf P^\perp f^\varepsilon)
-v\cdot\nabla_x\mathbf P^\perp f^\varepsilon \nonumber\\
&-\varepsilon(v\times B^\varepsilon)\cdot\nabla_vf^\varepsilon
-\varepsilon^2E^\varepsilon\cdot\nabla_vf^\varepsilon
+\frac{\varepsilon^2}{2}(E^\varepsilon\cdot v)f^\varepsilon, \, \widehat\Xi\mu^{1/2}\Big\rangle,
\end{align}
where we have used the identity
\begin{align*}
\Gamma(\mathbf Pf^\varepsilon,\mathbf Pf^\varepsilon)
=\frac12L\big\{\mu^{-1/2}(\mathbf Pf^\varepsilon)^2\big\}.
\end{align*}

We next rewrite the forcing terms involving only the slow components in \eqref{limit:acoustic-original-q} and \eqref{limit:acoustic-original-w}.
Using
\begin{align*}
u^\varepsilon=\mathcal P u^\varepsilon+\mathcal P^\perp u^\varepsilon,\qquad
\theta^\varepsilon=\Big(\frac35\theta^\varepsilon-\frac25\rho^\varepsilon\Big)+\frac25(\rho^\varepsilon+\theta^\varepsilon),
\end{align*}
we define
\begin{align}
\Phi^\varepsilon=\;&\Delta_x^{-1}\nabla_x\cdot\big(\mathcal P u^\varepsilon\times B^\varepsilon\big)
-\Delta_x^{-1}\nabla_x\cdot\nabla_x\cdot
\Big( \mathcal P u^\varepsilon\otimes \mathcal P u^\varepsilon-\frac{|\mathcal P u^\varepsilon|^2}{3}\mathbb I\Big),
\label{limit:acoustic-potential-phi}\\
\Psi^\varepsilon=\;&\kappa\nabla_x\Big(\frac35\theta^\varepsilon-\frac25\rho^\varepsilon\Big)
-\mathcal P^\perp\bigg\{\mathcal P u^\varepsilon\Big(\frac35\theta^\varepsilon-\frac25\rho^\varepsilon\Big)\bigg\}.
\label{limit:acoustic-potential-psi}
\end{align}
Then, we have
\begin{align*}
\nabla_x\Phi^\varepsilon
&=\mathcal P^\perp\bigg\{\mathcal P u^\varepsilon\times B^\varepsilon
-\nabla_x\cdot\Big(\mathcal P u^\varepsilon\otimes\mathcal P u^\varepsilon
-\frac{|\mathcal P u^\varepsilon|^2}{3}\mathbb I\Big)\bigg\},\\
\frac53\nabla_x\cdot\Psi^\varepsilon
&=\frac{5\kappa}{3}\Delta_x\Big(\frac35\theta^\varepsilon-\frac25\rho^\varepsilon\Big)
-\frac53\nabla_x\cdot\bigg\{\mathcal P u^\varepsilon\Big(\frac35\theta^\varepsilon-\frac25\rho^\varepsilon\Big)\bigg\}.
\end{align*}
Substituting \eqref{limit:stress-closure} and \eqref{limit:heat-closure} into \eqref{limit:acoustic-original-q}--\eqref{limit:acoustic-original-w} and using these identities, we obtain
\begin{align}
&\partial_t(\rho^\varepsilon+\theta^\varepsilon)+\frac5{3\varepsilon}\nabla_x\cdot \mathcal P^\perp u^\varepsilon
-\frac{2\kappa}{3}\Delta_x(\rho^\varepsilon+\theta^\varepsilon) \nonumber\\
&\qquad\qquad\qquad =\frac53\nabla_x\cdot\Psi^\varepsilon+\mathcal N_1^\varepsilon
+\frac{2\varepsilon}{3}\partial_t\nabla_x\cdot \big\langle\widehat{\mathcal B}\mu^{1/2},\mathbf P^\perp f^\varepsilon\big\rangle
+R_1^\varepsilon, \label{limit:acoustic-parabolic-q}\\
&\partial_t\mathcal P^\perp u^\varepsilon+\frac1\varepsilon\nabla_x(\rho^\varepsilon+\theta^\varepsilon)
-\frac{4\nu}{3}\Delta_x\mathcal P^\perp u^\varepsilon \nonumber\\
&\qquad\qquad\qquad = \nabla_x\Phi^\varepsilon+\mathcal N_2^\varepsilon
+\varepsilon\partial_t\mathcal P^\perp\nabla_x\cdot \big\langle\widehat{\mathcal A}\mu^{1/2},\mathbf P^\perp f^\varepsilon\big\rangle
+R_2^\varepsilon, \label{limit:acoustic-parabolic-w}
\end{align}
where
\begin{align} \label{N1-define}
\mathcal N_1^\varepsilon
=\;&-\frac53\nabla_x\cdot
\Big( \mathcal P^\perp u^\varepsilon\theta^\varepsilon+\frac25\mathcal P u^\varepsilon(\rho^\varepsilon+\theta^\varepsilon)\Big),\\
\mathcal N_2^\varepsilon
=\;&\mathcal P^\perp\bigg\{ \mathcal P^\perp u^\varepsilon\times B^\varepsilon
-\nabla_x\cdot\Big(u^\varepsilon\otimes u^\varepsilon
-\mathcal P u^\varepsilon\otimes\mathcal P u^\varepsilon-\frac{|u^\varepsilon|^2-|\mathcal P u^\varepsilon|^2}{3}\mathbb I\Big) \bigg\},\label{N2-define}
\end{align}
and
\begin{align}\label{R12-define}
R_1^\varepsilon=-\frac23\nabla_x\cdot r_{\mathcal B}^\varepsilon+\frac{2\varepsilon}{3}E^\varepsilon\cdot u^\varepsilon,\qquad
R_2^\varepsilon=\mathcal P^\perp\left\{E^\varepsilon+\varepsilon\rho^\varepsilon E^\varepsilon-\nabla_x\cdot r_{\mathcal A}^\varepsilon\right\}.
\end{align}

We now collect the estimates for the quantities in \eqref{limit:acoustic-parabolic-q}--\eqref{limit:acoustic-parabolic-w}.
First, the definition of $r_\Xi^\varepsilon$ in \eqref{r-define}, the uniform energy estimate \eqref{thm1 estimate 1} and the time-decay estimate \eqref{thm1 estimate 2} imply that
\begin{align}
&\left\|\big\langle\widehat\Xi\mu^{1/2},\mathbf P^\perp f^\varepsilon\big\rangle\right\|_{L^\infty(\mathbb R^+;H_x^N)}\leq C,\nonumber\\
&\left\|\big\langle\widehat\Xi\mu^{1/2},\mathbf P^\perp f^\varepsilon\big\rangle\right\|_{L^2(\mathbb R^+;H_x^N)}
+\|r_\Xi^\varepsilon\|_{L^2(\mathbb R^+;H_x^{N-1})}\leq C\varepsilon, \qquad \Xi\in\{\mathcal A,\mathcal B\}.
\label{limit:closure-remainders}
\end{align}
Moreover, \eqref{limit:slow-momentum}--\eqref{limit:slow-temperature}, the Maxwell equations \eqref{limit:maxwell-equations}, the uniform energy estimate \eqref{thm1 estimate 1} and the time-decay estimate \eqref{thm1 estimate 2} yield
\begin{align}
&\bigg\|\partial_t\Big(\mathcal P u^\varepsilon, \frac35\theta^\varepsilon-\frac25\rho^\varepsilon\Big)\bigg\|_{L^2(\mathbb{R}^+;H_x^{N-1})}
+\|\partial_tB^\varepsilon\|_{L^2(\mathbb{R}^+;H_x^{N-2})} \leq C,
\label{limit:slow-time-bounds}
\end{align}
where we have used the bound
\begin{align*}
\|\mathcal{P}E^\varepsilon\|_{H^{N-1}} \lesssim \|E^\varepsilon\|_{H^{N-1}}
\lesssim \|\Lambda^{1-\varrho} E^\varepsilon\| +\sum_{1\leq |\alpha|\leq N-1} \|\partial^\alpha E^\varepsilon\|
\lesssim \sqrt{\overline{\mathcal{D}}_{N,l}(t)}.
\end{align*}
By the Maxwell equations \eqref{limit:maxwell-equations}, we have
\begin{align}\label{E-estimate1}
\|\mathcal P^\perp E^\varepsilon\|_{H^{N-1}}
=\varepsilon\|\nabla_x\Delta_x^{-1}\rho^\varepsilon\|_{H^{N-1}}
\lesssim \varepsilon \|\Lambda^{-\varrho}\rho^\varepsilon\|+ \varepsilon\|\rho^\varepsilon\|_{H^{N-2}}
\leq C\varepsilon.
\end{align}
The definitions of $\Phi^\varepsilon$ and $\Psi^\varepsilon$ in \eqref{limit:acoustic-potential-phi}--\eqref{limit:acoustic-potential-psi}, together with \eqref{limit:slow-time-bounds}, imply
\begin{align}
\|(\Phi^\varepsilon,\Psi^\varepsilon)\|_{L^\infty(\mathbb{R}^+;H_x^{N-1})} \leq C,\qquad
\|\partial_t(\Phi^\varepsilon,\Psi^\varepsilon)\|_{L^2(\mathbb{R}^+;H_x^{N-2})}\leq C,
\label{limit:acoustic-potential-bounds}
\end{align}
where we have used Lemma \ref{negative embedding theorem}, together with H\"{o}lder's inequality, to obtain
\begin{align*}
\left\|\Delta_x^{-1}\nabla_x\cdot(\mathcal P u^\varepsilon\times B^\varepsilon)\right\|_{L^2}
\lesssim\|\mathcal P u^\varepsilon\times B^\varepsilon\|_{L^{\frac{6}{5}}}
\lesssim\|\mathcal P u^\varepsilon\| \|B^\varepsilon\|_{L^3},
\end{align*}
and
\begin{align*}
&\left\|\partial_t\Delta_x^{-1}\nabla_x\cdot(\mathcal P u^\varepsilon\times B^\varepsilon)\right\|_{L^2(\mathbb{R}^+;L_x^2)}\\
\lesssim\;&
\|\partial_t\mathcal P u^\varepsilon\|_{L^2(\mathbb{R}^+;L_x^2)}\|B^\varepsilon\|_{L^\infty(\mathbb{R}^+;L_x^3)}
+\|\mathcal P u^\varepsilon\|_{L^\infty(\mathbb{R}^+;L_x^3)}\|\partial_tB^\varepsilon\|_{L^2(\mathbb{R}^+;L_x^2)}.
\end{align*}
Furthermore, using \eqref{N1-define}--\eqref{R12-define}, H\"older's inequality, \eqref{Sobolev-ineq} and \eqref{limit:closure-remainders}, we obtain
\begin{align}
&\|(\mathcal N_1^\varepsilon,\mathcal N_2^\varepsilon)\|_{H_x^{N-3}}
\leq C\|(\rho^\varepsilon+\theta^\varepsilon,
\mathcal P^\perp u^\varepsilon)\|_{H_x^{N-2}},
\label{limit:acoustic-nonlinearity-bound}\\
&\|R_1^\varepsilon\|_{L^2(\mathbb R^+;H_x^{N-2})}
+\|R_2^\varepsilon-\mathcal P^\perp E^\varepsilon\|_{L^2(\mathbb R^+;H_x^{N-2})}
\leq C\varepsilon.
\label{limit:acoustic-remainder-global-parts}
\end{align}
Together with \eqref{E-estimate1}, this yields, for any $T>0$,
\begin{align}
\|(R_1^\varepsilon,R_2^\varepsilon)\|_{L^2(0,T;H_x^{N-2})}
\leq C\varepsilon\left(1+\sqrt T\right).
\label{limit:acoustic-remainder-bound}
\end{align}

We now define the corrected acoustic variables
\begin{align}
\widetilde q^\varepsilon
&:=\rho^\varepsilon+\theta^\varepsilon-\frac{2\varepsilon}{3}\nabla_x\cdot\big\langle\widehat{\mathcal B}\mu^{1/2},
\mathbf P^\perp f^\varepsilon\big\rangle-\varepsilon\Phi^\varepsilon, \nonumber \\
\widetilde w^\varepsilon
&:=\mathcal P^\perp u^\varepsilon-\varepsilon\mathcal P^\perp\nabla_x\cdot\big\langle\widehat{\mathcal A}\mu^{1/2},
\mathbf P^\perp f^\varepsilon\big\rangle-\varepsilon\Psi^\varepsilon. \label{limit:corrected-acoustic-variables}
\end{align}
By \eqref{limit:closure-remainders} and \eqref{limit:acoustic-potential-bounds}, we have
\begin{align}
\big\|(\rho^\varepsilon+\theta^\varepsilon-\widetilde q^\varepsilon,
\mathcal P^\perp u^\varepsilon-\widetilde w^\varepsilon)\big\|_{L^\infty(\mathbb R^+;H_x^{N-1})}
\leq C\varepsilon.
\label{limit:acoustic-correction-small}
\end{align}
Substituting \eqref{limit:corrected-acoustic-variables} into \eqref{limit:acoustic-parabolic-q}--\eqref{limit:acoustic-parabolic-w} yields
\begin{align}
&\partial_t\widetilde q^\varepsilon+\frac5{3\varepsilon}\nabla_x\cdot\widetilde w^\varepsilon-\frac{2\kappa}{3}\Delta_x\widetilde q^\varepsilon
=\mathcal N_1^\varepsilon+\widetilde R_1^\varepsilon,\nonumber\\
&\partial_t\widetilde w^\varepsilon+\frac1\varepsilon\nabla_x\widetilde q^\varepsilon-\frac{4\nu}{3}\Delta_x\widetilde w^\varepsilon
=\mathcal N_2^\varepsilon+\widetilde R_2^\varepsilon, \label{limit:corrected-acoustic-system}
\end{align}
where
\begin{align*}
\widetilde R_1^\varepsilon
=\;&R_1^\varepsilon
-\frac53\nabla_x\cdot\mathcal P^\perp\nabla_x\cdot\big\langle\widehat{\mathcal A}\mu^{1/2},\mathbf P^\perp f^\varepsilon\big\rangle
-\varepsilon\partial_t\Phi^\varepsilon+\frac{2\kappa\varepsilon}{3}\Delta_x
\Big(\frac23\nabla_x\cdot\big\langle\widehat{\mathcal B}\mu^{1/2},\mathbf P^\perp f^\varepsilon\big\rangle+\Phi^\varepsilon\Big),\\
\widetilde R_2^\varepsilon
=\;&R_2^\varepsilon
-\frac23\nabla_x\nabla_x\cdot\big\langle\widehat{\mathcal B}\mu^{1/2},\mathbf P^\perp f^\varepsilon\big\rangle
-\varepsilon\partial_t\Psi^\varepsilon+\frac{4\nu\varepsilon}{3}\Delta_x
\Big(\mathcal P^\perp\nabla_x\cdot\big\langle\widehat{\mathcal A}\mu^{1/2},\mathbf P^\perp f^\varepsilon\big\rangle+\Psi^\varepsilon\Big).
\end{align*}
By the definitions of $\Phi^\varepsilon$ and $\Psi^\varepsilon$ in \eqref{limit:acoustic-potential-phi}--\eqref{limit:acoustic-potential-psi}, H\"older's inequality and \eqref{Sobolev-ineq}, we have
\begin{align*}
\|\Delta_x\Phi^\varepsilon\|_{H_x^{N-3}}
&\lesssim
\|\mathcal P u^\varepsilon\times B^\varepsilon\|_{H_x^{N-2}}+\|\mathcal P u^\varepsilon\otimes\mathcal P u^\varepsilon\|_{H_x^{N-1}}
\lesssim \sqrt{\overline{\mathcal D}_{N,l}(t)},\\
\|\Delta_x\Psi^\varepsilon\|_{H_x^{N-3}}
&\lesssim
\left\|\frac35\theta^\varepsilon-\frac25\rho^\varepsilon\right\|_{H_x^N}
+\left\|\mathcal P u^\varepsilon\Big(\frac35\theta^\varepsilon-\frac25\rho^\varepsilon\Big)\right\|_{H_x^{N-1}}
\lesssim \sqrt{\overline{\mathcal D}_{N,l}(t)},
\end{align*}
which, together with \eqref{thm1 estimate 1}, gives
\begin{align}
\|\Delta_x(\Phi^\varepsilon,\Psi^\varepsilon)\|_{L^2(\mathbb R^+;H_x^{N-3})}\leq C.
\label{limit:acoustic-potential-laplacians}
\end{align}
Using \eqref{limit:closure-remainders}, \eqref{limit:acoustic-potential-bounds}, \eqref{limit:acoustic-potential-laplacians}, and \eqref{limit:acoustic-remainder-bound}, we obtain, for any $T>0$,
\begin{align}
\|(\widetilde R_1^\varepsilon, \widetilde R_2^\varepsilon)\|_{L^2(0,T;H_x^{N-3})}
\leq C\varepsilon \left(1+\sqrt T\right).
\label{limit:corrected-acoustic-remainders-global-parts}
\end{align}
Applying $\partial^\alpha$ with $|\alpha|\leq N-2$ to the two equations in \eqref{limit:corrected-acoustic-system}, taking the $L_x^2$ inner products with $\frac35\partial^\alpha\widetilde q^\varepsilon$ and $\partial^\alpha\widetilde w^\varepsilon$, respectively, and summing over $|\alpha|\leq N-2$, we obtain, by integration by parts, \eqref{limit:acoustic-nonlinearity-bound}, \eqref{limit:acoustic-correction-small}, and Young's inequality,
\begin{align*}
&\frac12\frac{\d}{\d t}
\Big\{\frac35\|\widetilde q^\varepsilon\|_{H_x^{N-2}}^2+\|\widetilde w^\varepsilon\|_{H_x^{N-2}}^2\Big\}
+\frac{\kappa}{5}\|\nabla_x\widetilde q^\varepsilon\|_{H_x^{N-2}}^2
+\frac{2\nu}{3}\|\nabla_x\widetilde w^\varepsilon\|_{H_x^{N-2}}^2\\
\lesssim\;&
\|(\widetilde q^\varepsilon,\widetilde w^\varepsilon)\|_{H_x^{N-2}}^2
+\varepsilon^2
+\|(\widetilde R_1^\varepsilon,\widetilde R_2^\varepsilon)\|_{H_x^{N-3}}^2.
\end{align*}
Gronwall's inequality, together with \eqref{limit:acoustic-correction-small} and \eqref{limit:corrected-acoustic-remainders-global-parts}, yields
\begin{align}
\sup_{0\leq t\leq T}\|(\rho^\varepsilon+\theta^\varepsilon,\mathcal P^\perp u^\varepsilon)(t)\|_{H_x^{N-2}}^2
\leq Ce^{CT}
\left\{\|(\rho^\varepsilon+\theta^\varepsilon)(0)\|_{H_x^{N-2}}^2
+\|\mathcal P^\perp u^\varepsilon(0)\|_{H_x^{N-2}}^2+\varepsilon^2\right\}. \label{limit:acoustic-uniform-smallness}
\end{align}
The strong convergence of the initial data and the compatibility conditions $\rho_0+\theta_0=0$ and $\mathcal P u_0=u_0$ imply
\begin{align*}
(\rho^\varepsilon+\theta^\varepsilon)(0)\to0,\qquad\mathcal P^\perp u^\varepsilon(0)\to0 \quad\text{strongly in }H_x^{N-2}
\end{align*}
as $\varepsilon \to 0$. Consequently, for any fixed $T>0$, we have
\begin{align*}
\sup_{0\leq t\leq T}\|(\rho^\varepsilon+\theta^\varepsilon,\mathcal P^\perp u^\varepsilon)(t)\|_{H_x^{N-2}}\to 0
\end{align*}
as $\varepsilon \to 0$, which gives
\begin{align*}
&\sup_{0\leq t\leq T}\|(\rho^\varepsilon+\theta^\varepsilon, \mathcal P^\perp u^\varepsilon)(t)\|_{H_x^{N-1}}\\
\lesssim\;&
\Big\{\sup_{0\leq t\leq T}\|(\rho^\varepsilon+\theta^\varepsilon,\mathcal P^\perp u^\varepsilon)(t)\|_{H_x^{N-2}}\Big\}^{\frac{1}{2}}
\Big\{\sup_{0\leq t\leq T}\|(\rho^\varepsilon+\theta^\varepsilon,\mathcal P^\perp u^\varepsilon)(t)\|_{H_x^N}\Big\}^{\frac{1}{2}}
\to 0
\end{align*}
as $\varepsilon \to 0$.
Hence, using the time continuity of the solutions, we conclude that, for any $T>0$,
\begin{align}
(\rho^\varepsilon+\theta^\varepsilon,\mathcal P^\perp u^\varepsilon)
\to(0,0) \quad\text{strongly in }C([0,T];H_x^{N-1})
\label{limit:acoustic-strong-C}
\end{align}
as $\varepsilon \to 0$.
Combining \eqref{limit:acoustic-strong-C} with the local strong convergence \eqref{limit:strong-slow-components} yields, as $\varepsilon \to 0$,
\begin{align*}
(\rho^\varepsilon,u^\varepsilon,\theta^\varepsilon) \to(\rho,u,\theta) \quad\text{strongly in } C([0,T];H_{loc}^{N-1}).
\end{align*}

We next establish the strong convergence of the microscopic component $\mathbf P^\perp f^\varepsilon$.
Taking the inner product of \eqref{rrVMB} with $\mathbf P^\perp f^\varepsilon$ and integrating by parts, we obtain
\begin{align*}
&\frac12\frac{\d}{\d t}
\|\mathbf P^\perp f^\varepsilon\|^2+\frac1{\varepsilon^2}\big(L\mathbf P^\perp f^\varepsilon,\mathbf P^\perp f^\varepsilon\big)\\
=\;&-\frac1\varepsilon
\big(v\cdot\nabla_x\mathbf Pf^\varepsilon,\mathbf P^\perp f^\varepsilon\big)
+\frac1\varepsilon\big(\Gamma(f^\varepsilon,f^\varepsilon),\mathbf P^\perp f^\varepsilon\big)\\
&-\varepsilon
\big(E^\varepsilon\cdot\nabla_v\mathbf Pf^\varepsilon,\mathbf P^\perp f^\varepsilon\big)
+\frac{\varepsilon}{2}\big((E^\varepsilon\cdot v)f^\varepsilon,\mathbf P^\perp f^\varepsilon\big),
\end{align*}
where we have used \eqref{Pf define} to obtain
\begin{align*}
\big((v\times B^\varepsilon)\cdot\nabla_v\mathbf Pf^\varepsilon, \mathbf P^\perp f^\varepsilon\big)=0.
\end{align*}
By \eqref{norm inequality}, \eqref{Gamma-unweighted}, \eqref{Sobolev-ineq}, and the uniform energy estimate \eqref{thm1 estimate 1}, the first two terms satisfy
\begin{align*}
&-\frac1\varepsilon\big(v\cdot\nabla_x\mathbf Pf^\varepsilon,\mathbf P^\perp f^\varepsilon\big)
+\frac1\varepsilon\big(\Gamma(f^\varepsilon,f^\varepsilon),\mathbf P^\perp f^\varepsilon\big)\\
\lesssim\;&
\frac1\varepsilon
\left\{\|\nabla_x\mathbf Pf^\varepsilon\|
+\|f^\varepsilon\|_{H_x^2L_v^2}\big(\|f^\varepsilon\|_D+\|f^\varepsilon\|\big)\right\}
\|\mathbf P^\perp f^\varepsilon\|_D
\leq \frac{\lambda}{2\varepsilon^2}\|\mathbf P^\perp f^\varepsilon\|_D^2+C_T.
\end{align*}
The electric terms are bounded by
\begin{align*}
&-\varepsilon
\big(E^\varepsilon\cdot\nabla_v\mathbf Pf^\varepsilon,\mathbf P^\perp f^\varepsilon\big)
+\frac{\varepsilon}{2}\big((E^\varepsilon\cdot v)f^\varepsilon,\mathbf P^\perp f^\varepsilon\big)\\
\lesssim\;&
\varepsilon\|E^\varepsilon\|_{L^\infty}
\left\{\|\mathbf Pf^\varepsilon\|\|\mathbf P^\perp f^\varepsilon\|
+\|\langle v\rangle^{\frac{1}{2}}f^\varepsilon\|\|\langle v\rangle^{\frac{1}{2}}\mathbf P^\perp f^\varepsilon\|\right\}
\leq C_T.
\end{align*}
Combining these estimates with \eqref{L-coercive-unweighted} gives
\begin{align*}
\frac{\d}{\d t}\|\mathbf P^\perp f^\varepsilon\|^2
+\frac{\lambda}{\varepsilon^2} \|\mathbf P^\perp f^\varepsilon\|_D^2
\leq C_T.
\end{align*}
Multiplying this inequality by $2\|\mathbf P^\perp f^\varepsilon(t)\|_{L_{x,v}^2}^2$ and integrating in time, we deduce that
\begin{align*}
\sup_{0\leq t\leq T}\|\mathbf P^\perp f^\varepsilon(t)\|^4
\leq \|\mathbf P^\perp f_0^\varepsilon\|^4
+C_T\int_0^T\|\mathbf P^\perp f^\varepsilon(t)\|^2 \d t
\leq \|\mathbf P^\perp f_0^\varepsilon\|^4+C_T\varepsilon^2
\to 0
\end{align*}
as $\varepsilon \to 0$, since $f_0^\varepsilon\to f_0$ and $\mathbf P^\perp f_0=0$.
Interpolating with the uniform energy estimate \eqref{thm1 estimate 1} and using the time continuity of the solutions, we conclude that, as $\varepsilon \to 0$,
\begin{align}
\mathbf P^\perp f^\varepsilon\to0
\quad\text{strongly in }C([0,T];H_x^{N-1}L_v^2).
\label{limit:microscopic-strong-C}
\end{align}
Combining the macroscopic and microscopic convergence results, we obtain that, as $\varepsilon \to 0$,
\begin{align*}
f^\varepsilon\to f \quad\text{strongly in } C([0,T];H_{loc}^{N-1}L_v^2).
\end{align*}

We now identify the limiting system. By \eqref{limit:closure-remainders} and integration by parts in time, we have
\begin{align*}
\varepsilon\partial_t\big\langle\widehat\Xi\mu^{1/2},
\mathbf P^\perp f^\varepsilon\big\rangle\longrightarrow0
\quad\text{in the sense of distributions},\qquad \Xi\in\{\mathcal A,\mathcal B\}
\end{align*}
as $\varepsilon \to 0$.
Using \eqref{limit:stress-closure}--\eqref{limit:heat-closure}, \eqref{limit:closure-remainders}, and the strong convergences above, we pass to the limit in \eqref{limit:slow-momentum}, \eqref{limit:slow-temperature}, and \eqref{limit:maxwell-equations} and conclude that $(\rho,u,\theta,E,B)$ satisfies the incompressible NSFM system \eqref{INSFM limit}.
The initial conditions
\begin{align*}
u(0)=u_0,\qquad \theta(0)=\theta_0,\qquad E(0)=E_0,\qquad B(0)=B_0
\end{align*}
follow from \eqref{limit:strong-slow-components} and the assumed convergence of the initial data.
Moreover, the uniform bound \eqref{limit:uniform-bounds} and weak lower semicontinuity, together with the limiting system, yield
\begin{align*}
(\rho,u,\theta,E,B)\in C([0,\infty);H_x^{N-1})\cap L^\infty(0,\infty;H_x^N).
\end{align*}
An energy estimate for the difference of two solutions then yields uniqueness in this class.

We next strengthen the above local convergence to strong convergence in the whole space $\mathbb R^3$.
For this purpose, define the corrected differences
\begin{align}
U^\varepsilon
&=\mathcal P u^\varepsilon-u-\varepsilon\mathcal P\nabla_x\cdot\big\langle\widehat{\mathcal A}\mu^{1/2},
\mathbf P^\perp f^\varepsilon\big\rangle,\nonumber\\
\Theta^\varepsilon
&=\frac35\theta^\varepsilon-\frac25\rho^\varepsilon-\theta-\frac{2\varepsilon}{5}\nabla_x\cdot
\big\langle\widehat{\mathcal B}\mu^{1/2}, \mathbf P^\perp f^\varepsilon\big\rangle.
\label{limit:global-slow-errors}
\end{align}
Then $\nabla_x\cdot U^\varepsilon=0$, and the following identities hold
\begin{align}
u^\varepsilon-u
&=U^\varepsilon+\mathcal P^\perp u^\varepsilon
+\varepsilon\mathcal P\nabla_x\cdot \big\langle\widehat{\mathcal A}\mu^{1/2},\mathbf P^\perp f^\varepsilon\big\rangle,\nonumber\\
\theta^\varepsilon-\theta
&=\Theta^\varepsilon+\frac25(\rho^\varepsilon+\theta^\varepsilon)
+\frac{2\varepsilon}{5}\nabla_x\cdot \big\langle\widehat{\mathcal B}\mu^{1/2},\mathbf P^\perp f^\varepsilon\big\rangle.
\label{limit:global-error-decomposition}
\end{align}
Substituting \eqref{limit:stress-closure}--\eqref{limit:heat-closure} into \eqref{limit:slow-momentum}--\eqref{limit:slow-temperature} and subtracting the limiting system \eqref{INSFM limit}, we obtain
\begin{align}
&\partial_tU^\varepsilon-\nu\Delta_xU^\varepsilon
+\mathcal P\nabla_x\cdot
\big(u^\varepsilon\otimes u^\varepsilon-u\otimes u\big)
=\mathcal P(E^\varepsilon-E) +\mathcal P\big(u^\varepsilon\times B^\varepsilon-u\times B\big) \nonumber\\
&\qquad\qquad\qquad\qquad +\varepsilon\mathcal P(\rho^\varepsilon E^\varepsilon)
-\mathcal P\nabla_x\cdot r_{\mathcal A}^\varepsilon
+\nu\varepsilon\Delta_x\mathcal P\nabla_x\cdot \big\langle\widehat{\mathcal A}\mu^{1/2}, \mathbf P^\perp f^\varepsilon\big\rangle,\label{limit:global-velocity-error} \\
&\partial_t\Theta^\varepsilon-\kappa\Delta_x\Theta^\varepsilon
+\nabla_x\cdot\big(u^\varepsilon\theta^\varepsilon-u\theta\big)
=\frac{2\kappa}{5}\Delta_x(\rho^\varepsilon+\theta^\varepsilon) -\frac25\nabla_x\cdot r_{\mathcal B}^\varepsilon \nonumber\\
&\qquad\qquad\qquad\qquad+\frac{2\varepsilon}{5}E^\varepsilon\cdot u^\varepsilon
+\frac{2\kappa\varepsilon}{5}\Delta_x\nabla_x\cdot \big\langle\widehat{\mathcal B}\mu^{1/2}, \mathbf P^\perp f^\varepsilon\big\rangle, \label{limit:global-temperature-error}\\
&\partial_t(E^\varepsilon-E)-\nabla_x\times(B^\varepsilon-B)
=-U^\varepsilon-\mathcal P^\perp u^\varepsilon
-\varepsilon\mathcal P\nabla_x\cdot \big\langle\widehat{\mathcal A}\mu^{1/2}, \mathbf P^\perp f^\varepsilon\big\rangle, \label{limit:global-electric-error}\\
&\partial_t(B^\varepsilon-B)+\nabla_x\times(E^\varepsilon-E)=0. \label{limit:global-magnetic-error}
\end{align}
Taking the $L_x^2$ inner products of \eqref{limit:global-velocity-error}--\eqref{limit:global-magnetic-error} with $U^\varepsilon$, $\Theta^\varepsilon$, $E^\varepsilon-E$, and $B^\varepsilon-B$, respectively, summing the resulting identities, and integrating by parts, we obtain
\begin{align*}
&\frac12\frac{\d}{\d t}
\big\|(U^\varepsilon,\Theta^\varepsilon,E^\varepsilon-E,B^\varepsilon-B)\big\|^2
+\nu\|\nabla_xU^\varepsilon\|^2+\kappa\|\nabla_x\Theta^\varepsilon\|^2\\
=\;&
\big(u^\varepsilon\otimes u^\varepsilon-u\otimes u,\nabla_xU^\varepsilon\big)
+\big(u^\varepsilon\theta^\varepsilon-u\theta,\nabla_x\Theta^\varepsilon\big)
+\big(u^\varepsilon\times B^\varepsilon-u\times B,U^\varepsilon\big)\\
&+\varepsilon\big(\rho^\varepsilon E^\varepsilon,U^\varepsilon\big)
+\frac{2\varepsilon}{5}\big(E^\varepsilon\cdot u^\varepsilon,\Theta^\varepsilon\big)
+\big(r_{\mathcal A}^\varepsilon,\nabla_xU^\varepsilon\big)+\frac25\big(r_{\mathcal B}^\varepsilon,\nabla_x\Theta^\varepsilon\big)\\
&-\frac{2\kappa}{5}\big(\nabla_x(\rho^\varepsilon+\theta^\varepsilon),\nabla_x\Theta^\varepsilon\big)
-\nu\varepsilon\Big(\nabla_x\mathcal P\nabla_x\cdot\big\langle\widehat{\mathcal A}\mu^{1/2},\mathbf P^\perp f^\varepsilon\big\rangle,\nabla_xU^\varepsilon\Big)\\
&-\frac{2\kappa\varepsilon}{5}\Big(\nabla_x\nabla_x\cdot\big\langle\widehat{\mathcal B}\mu^{1/2}, \mathbf P^\perp f^\varepsilon\big\rangle,\nabla_x\Theta^\varepsilon\Big)
-\Big(\mathcal P^\perp u^\varepsilon +\varepsilon\mathcal P\nabla_x\cdot \big\langle\widehat{\mathcal A}\mu^{1/2},\mathbf P^\perp f^\varepsilon\big\rangle, E^\varepsilon-E\Big).
\end{align*}
By \eqref{limit:global-error-decomposition}, \eqref{limit:uniform-bounds}, \eqref{limit:closure-remainders}, and \eqref{Sobolev-ineq}, one has
\begin{align*}
\|u^\varepsilon\otimes u^\varepsilon-u\otimes u\|
+\|u^\varepsilon\theta^\varepsilon-u\theta\|
\lesssim
\|(U^\varepsilon,\Theta^\varepsilon)\|
+\|(\rho^\varepsilon+\theta^\varepsilon,\mathcal P^\perp u^\varepsilon)\| +\varepsilon.
\end{align*}
Moreover, using
$\big(U^\varepsilon\times B^\varepsilon,U^\varepsilon\big)=0$,
we have
\begin{align*}
\big(u^\varepsilon\times B^\varepsilon-u\times B,U^\varepsilon\big)
\lesssim \|U^\varepsilon\|^2+\|B^\varepsilon-B\|^2+\|\mathcal P^\perp u^\varepsilon\|^2+\varepsilon^2.
\end{align*}
Using the uniform $H_x^2$ bounds for the microscopic moments in \eqref{limit:closure-remainders} and Young's inequality to estimate the remaining terms, we obtain
\begin{align}
&\frac12\frac{\d}{\d t}
\big\|(U^\varepsilon,\Theta^\varepsilon,E^\varepsilon-E,B^\varepsilon-B)\big\|^2
+\frac\nu2\|\nabla_xU^\varepsilon\|^2 +\frac\kappa2\|\nabla_x\Theta^\varepsilon\|^2\nonumber\\
\lesssim\;&
\big\|(U^\varepsilon,\Theta^\varepsilon, E^\varepsilon-E,B^\varepsilon-B)\big\|^2
+\big\|(\rho^\varepsilon+\theta^\varepsilon, \mathcal P^\perp u^\varepsilon)\big\|_{H_x^1}^2
+\big\|(r_{\mathcal A}^\varepsilon,r_{\mathcal B}^\varepsilon)\big\|_{L_x^2}^2+\varepsilon^2.
\label{limit:global-error-energy}
\end{align}
The initial differences satisfy
\begin{align*}
\big\|(U^\varepsilon,\Theta^\varepsilon,
E^\varepsilon-E,B^\varepsilon-B)(0)\big\|^2
\lesssim\|f_0^\varepsilon-f_0\|^2
+\|(E_0^\varepsilon-E_0,B_0^\varepsilon-B_0)\|^2
+\varepsilon^2.
\end{align*}
Hence, Gronwall's inequality and \eqref{limit:closure-remainders} imply that, for any $T>0$,
\begin{align}
&\sup_{0\leq t\leq T}
\left\|\Big(\mathcal P u^\varepsilon-u,
\frac35\theta^\varepsilon-\frac25\rho^\varepsilon-\theta,
E^\varepsilon-E,B^\varepsilon-B\Big)(t)\right\|^2\label{limit:global-slow-L2}\\
\leq\;& C_T\Bigg\{\|f_0^\varepsilon-f_0\|^2
+\|(E_0^\varepsilon-E_0,B_0^\varepsilon-B_0)\|^2
+\varepsilon^2
+\int_0^T\big\|(\rho^\varepsilon+\theta^\varepsilon,
\mathcal P^\perp u^\varepsilon)(t)\big\|_{H_x^1}^2 \d t\Bigg\}\to 0 \nonumber
\end{align}
as $\varepsilon \to 0$, where we have used the assumed convergence of the initial data and \eqref{limit:acoustic-strong-C}.

Combining \eqref{limit:global-slow-L2} with \eqref{limit:acoustic-strong-C} and \eqref{limit:microscopic-strong-C}, and interpolating with the uniform $H_x^N$ bound \eqref{limit:uniform-bounds}, we obtain
\begin{align}
f^\varepsilon&\to f &&\quad\text{strongly in }C([0,T];H_x^{N-1}L_v^2),\nonumber\\
(E^\varepsilon,B^\varepsilon)&\to(E,B) &&\quad\text{strongly in }C([0,T];H_x^{N-1}) \label{limit:strong-global-space}
\end{align}
as $\varepsilon \to 0$.
The same argument gives
\begin{align*}
\left(\mathcal P u^\varepsilon, \frac35\theta^\varepsilon-\frac25\rho^\varepsilon\right)
\to(u,\theta) \qquad\text{strongly in }C([0,T];H_x^{N-1})
\end{align*}
as $\varepsilon \to 0$.
Since \eqref{limit:global-slow-L2} applies to the entire family, these convergences hold without extracting a subsequence.

Finally, the convergence is uniform for all $t\geq0$.
Indeed, \eqref{thm1 estimate 1}--\eqref{thm1 estimate 2} and Lemma \ref{sobolev interpolation} yield
\begin{align*}
\|f^\varepsilon(t)\|_{H_x^{N-1}L_v^2}
\lesssim \|f^\varepsilon(t)\|^{\frac{1}{N}} \|f^\varepsilon(t)\|_{H_x^NL_v^2}^{1-\frac{1}{N}}
\leq C(1+t)^{-\frac{\varrho}{4N}}.
\end{align*}
By \eqref{limit:strong-global-space}, the limiting function $f$ satisfies the same decay estimate.
Consequently, for any $T>0$,
\begin{align*}
&\sup_{t\geq0}\|f^\varepsilon(t)-f(t)\|_{H^{N-1}_xL_v^2}
\leq\sup_{0\leq t\leq T}\|f^\varepsilon(t)-f(t)\|_{H^{N-1}_xL_v^2}+C(1+T)^{-\frac{\varrho}{4N}}.
\end{align*}
Letting first $\varepsilon\to0$ and then $T\to\infty$, we obtain
\begin{align}
\sup_{t\geq0}
\|f^\varepsilon(t)-f(t)\|_{H^{N-1}_xL_v^2}
\to 0.
\label{limit:strong-global-time}
\end{align}
The convergence of $\big(\mathcal P u^\varepsilon,\frac35\theta^\varepsilon-\frac25\rho^\varepsilon,E^\varepsilon,B^\varepsilon\big)$ to $(u,\theta, E,B)$ is also uniform for all $t\geq0$ by the same argument.
This completes the proof.
\end{proof}
\bigskip

\appendix
\section{Some Auxiliary Lemmas}\label{appendix:auxiliary}

This appendix collects some auxiliary results needed for the exponentially weighted estimates of $L$ and $\Gamma$ in Section \ref{new-weighted-estimate}.
We first recall some estimates for $L$ and $\Gamma$ without exponential velocity weights.
\begin{lemma}\label{lem:unweighted-L-Gamma}
Let $0<s<1$ and $\ell\geq0$.
\begin{itemize}
\item[(1)] For $\gamma>-3$, the following estimates hold
\begin{align}\label{L-coercive-unweighted}
\left\langle Lf,f\right\rangle&\gtrsim|\mathbf P^\perp f|_D^2,\\
\left|\left\langle Lf,g\right\rangle\right|&\lesssim|f|_D|g|_D,\label{L-upper-unweighted}\\
\left\langle Lf,\langle v\rangle^{2\ell}f\right\rangle
&\gtrsim|\langle v\rangle^\ell f|_D^2
-C|f|_{L^2(B_C)}^2.\label{polynomial-L-coercivity}
\end{align}

\item[(2)] For $\gamma>\max\{-3,-\frac32-2s\}$, there hold
\begin{align}
\left|\left\langle\Gamma(f,g),h\right\rangle\right|
\lesssim\;&\big\{|f|_{L^2_{\gamma/2+s}}|g|_D
+|g|_{L^2_{\gamma/2+s}}|f|_D \big\}|h|_D \nonumber\\
&+\min\big\{|f|_{L^2}|g|_{L^2_{\gamma/2+s}},
|f|_{L^2_{\gamma/2+s}}|g|_{L^2}\big\}|h|_D,
\label{Gamma-unweighted}
\end{align}
and
\begin{align}\label{polynomial-Gamma}
\begin{split}
\left|\left\langle\langle v\rangle^\ell\Gamma(f,g),h\right\rangle\right|
\lesssim\;&
\big\{|f|_{L^2_{\gamma/2+s}}|\langle v\rangle^\ell g|_D+|\langle v\rangle^\ell g|_{L^2_{\gamma/2+s}}|f|_D \big\}|h|_D \\
&+\min\big\{|f|_{L^2}|\langle v\rangle^\ell g|_{L^2_{\gamma/2+s}},|f|_{L^2_{\gamma/2+s}}|\langle v\rangle^\ell g|_{L^2}\big\}|h|_D.
\end{split}
\end{align}

\item[(3)] Let $\gamma>-3$ and let
$\zeta(v)\in\mathcal S(\mathbb R^3)$ satisfy
$\sum_{|\beta|\leq2}|\partial_v^\beta\zeta(v)|
\leq C e^{-c|v|^2}$ for some $C,c>0$.
For any $m\geq0$, there holds
\begin{equation}\label{Gamma-zeta}
\left|\left\langle\Gamma(f,g),\zeta\right\rangle\right|
\lesssim|f|_{L^2_{-m}}|g|_{L^2_{-m}}.
\end{equation}
\end{itemize}
\end{lemma}


\begin{proof}
The estimates \eqref{L-coercive-unweighted}, \eqref{L-upper-unweighted}, \eqref{polynomial-L-coercivity} and \eqref{Gamma-unweighted} have been shown in
\cite{AMUXY2012JFA, GS2011}.
The estimate \eqref{Gamma-zeta} can be found in \cite{GS2011}. To prove \eqref{polynomial-Gamma}, we first apply \eqref{Gamma-unweighted} with $g$ replaced by $\langle v\rangle^\ell g$ to obtain
\begin{align*}
\left|\left\langle\Gamma(f,\langle v\rangle^\ell g),h\right\rangle\right|
\lesssim\;&
\Big\{|f|_{L^2_{\gamma/2+s}}|\langle v\rangle^\ell g|_D+|f|_D|\langle v\rangle^\ell g|_{L^2_{\gamma/2+s}}\Big\}|h|_D\\
&+\min\Big\{|f|_{L^2}|\langle v\rangle^\ell g|_{L^2_{\gamma/2+s}},
|f|_{L^2_{\gamma/2+s}}|\langle v\rangle^\ell g|_{L^2}\Big\}|h|_D.
\end{align*}
By Lemma \ref{polynomial commutator}, we have
\begin{align*}
&\left|\left\langle\langle v\rangle^\ell\Gamma(f,g)-\Gamma(f,\langle v\rangle^\ell g),h\right\rangle\right|\\
\lesssim\;&
\Big\{|f|_{L^2_{\gamma/2+s}}|\langle v\rangle^\ell g|_D
+\min\big\{|f|_{L^2}|\langle v\rangle^\ell g|_{L^2_{\gamma/2+s}},
|f|_{L^2_{\gamma/2+s}}|\langle v\rangle^\ell g|_{L^2}\big\}\Big\}|h|_D.
\end{align*}
Combining these estimates with
\begin{align*}
\left\langle\langle v\rangle^\ell\Gamma(f,g),h\right\rangle
=\left\langle\Gamma(f,\langle v\rangle^\ell g),h\right\rangle
+\left\langle\langle v\rangle^\ell\Gamma(f,g)-\Gamma(f,\langle v\rangle^\ell g),h\right\rangle
\end{align*}
proves \eqref{polynomial-Gamma}.
\end{proof}
\smallskip

We next recall the following commutator estimate for polynomial velocity weights from \cite[Proposition 3.13]{AMUXY2012JFA}.

\begin{lemma}\label{polynomial commutator}
Let $0<s<1$ and
$\gamma>\max\{-3,-\frac{3}{2}-2s\}$.
Then, for any $\ell\geqslant0$, one has
\begin{align*}
&\left|\left\langle\langle v\rangle^\ell\Gamma(f,g) -\Gamma(f,\langle v\rangle^\ell g),h\right\rangle\right|\\
\lesssim\;&
\Big\{|f|_{L^2_{\gamma/2+s}}\big|\langle v\rangle^{\ell-s}g\big|_{L^2_{\gamma/2+s}}
+\min\big\{|f|_{L^2}\big|\langle v\rangle^{\ell-s}g\big|_{L^2_{\gamma/2+s}},
|f|_{L^2_{\gamma/2+s}}\big|\langle v\rangle^{\ell-s}g\big|_{L^2}\big\}\Big\}|h|_D.
\end{align*}
\end{lemma}
\smallskip

For $\eta>0$ and $\delta>-3$, define
\begin{align}
J_{1,\eta}^{\Phi_\delta}(f)
&=\iiint_{\mathbb{R}^6\times\mathbb{S}^2} |v-v_*|^\delta b(\cos\theta)\mu(v_*)^\eta
|f(v')-f(v)|^2\d\sigma\d v_*\d v, \label{J1define}\\
J_{2,\eta}^{\Phi_\delta}(f)
&=\iiint_{\mathbb{R}^6\times\mathbb{S}^2} |v-v_*|^\delta b(\cos\theta)|f(v_*)|^2
\left|\mu(v')^{\eta/2}-\mu(v)^{\eta/2}\right|^2
\d\sigma\d v_*\d v.\label{J2define}
\end{align}
Let $|\cdot|_{D_{\mu^\eta}}$ denote $|\cdot|_D$
with $\mu$ replaced by $\mu^\eta$. Thus, we have
\begin{align*}
|f|_{D_{\mu^\eta}}^2=J_{1,\eta}^{\Phi_\gamma}(f)+J_{2,\eta}^{\Phi_\gamma}(f),
\qquad |f|_{D_{\mu}}=|f|_D.
\end{align*}

We next present three related lemmas.
The first two follow from the estimates and comparison argument in \cite[Section 2.3, pp.~941--942]{AMUXY2012JFA}, while the third follows from \cite[Lemma 2.17]{AMUXY2012JFA}.

\begin{lemma}\label{lem:Gaussian-J2}
Let $0<s<1$ and $\gamma>-3$.
For any fixed $\eta>0$, we have
\begin{equation}\label{Gaussian-J2}
J_{2,\eta}^{\Phi_\gamma}(f) \approx |f|_{L^2_{\gamma/2+s}}^2.
\end{equation}
\end{lemma}
\smallskip

\begin{lemma}\label{lem:Gaussian-kinetic-comparison}
Let $0<s<1$ and $\gamma>-3$.
For any fixed $\eta>0$, we have
\begin{align}
J_{1,\eta}^{\Phi_\gamma}(f)
&\lesssim J_{1,\eta/2}^{\Phi_0}
\big(\langle v\rangle^{\gamma/2}f\big)+|f|_{L^2_{\gamma/2+s}}^2,
\label{Gaussian-kinetic-upper}\\
J_{1,2\eta}^{\Phi_0}
\big(\langle v\rangle^{\gamma/2}f\big)
&\lesssim J_{1,\eta}^{\Phi_\gamma}(f)+|f|_{L^2_{\gamma/2+s}}^2.
\label{Gaussian-kinetic-lower}
\end{align}
\end{lemma}
\smallskip

\begin{lemma}\label{lem:Gaussian-Maxwell-comparison}
Let $0<s<1$ and $\eta>0$.
Then
\begin{equation}\label{Gaussian-Maxwell-comparison}
\iiint_{\mathbb{R}^6\times\mathbb{S}^2} b(\cos\theta)|f(v_*)|^2|g(v')-g(v)|^2\d\sigma\d v_*\d v
\lesssim |f|_{L^2_s}^2 \left(J_{1,\eta}^{\Phi_0}(g)+|g|_{L^2_s}^2\right).
\end{equation}
\end{lemma}
\smallskip

The following lemma establishes the equivalence of the non-isotropic norms associated with different Gaussian weights.

\begin{lemma}\label{lem:Gaussian-norm-comparison}
Let $0<s<1$ and $\gamma>-3$.
For $\eta\in[1/2,2]$, we have
\begin{equation}\label{Gaussian-norm-equivalence}
|f|_{D_{\mu^\eta}}\approx |f|_D,
\end{equation}
where the comparison constants are independent of $\eta$.
\end{lemma}

\begin{proof}
By Lemmas \ref{lem:Gaussian-J2}, \ref{lem:Gaussian-kinetic-comparison}, and \ref{lem:Gaussian-Maxwell-comparison}, we have
\begin{align*}
|f|_{D_{\mu^{1/2}}}^2
&\lesssim J_{1,1/4}^{\Phi_0}\big(\langle v\rangle^{\gamma/2}f\big)+|f|_{L^2_{\gamma/2+s}}^2\\
&\lesssim J_{1,4}^{\Phi_0}\big(\langle v\rangle^{\gamma/2}f\big)+|f|_{L^2_{\gamma/2+s}}^2\\
&\lesssim J_{1,2}^{\Phi_\gamma}(f)+|f|_{L^2_{\gamma/2+s}}^2\\
&\lesssim |f|_{D_{\mu^2}}^2.
\end{align*}
For $\eta\in[1/2,2]$, a direct calculation yields
\begin{gather*}
\mu^2\leq\mu^\eta\leq\mu^{1/2}, \\
\frac14|\mu(v')-\mu(v)| \leq \left|\mu(v')^{\eta/2}-\mu(v)^{\eta/2}\right|
\leq 4\left|\mu(v')^{1/4}-\mu(v)^{1/4}\right|.
\end{gather*}
Applying these inequalities to the two terms in the definition of $|\cdot|_{D_{\mu^\eta}}$, we obtain
\begin{align*}
\frac1{16}|f|_{D_{\mu^2}}^2 \leq |f|_{D_{\mu^\eta}}^2 \leq16|f|_{D_{\mu^{1/2}}}^2 \lesssim |f|_{D_{\mu^2}}^2.
\end{align*}
It follows that
\begin{align*}
|f|_{D_{\mu^\eta}}\approx |f|_{D_{\mu^2}},
\end{align*}
with comparison constants independent of $\eta\in[1/2,2]$.
Taking $\eta=1$ gives $|f|_D\approx|f|_{D_{\mu^2}}$, which proves \eqref{Gaussian-norm-equivalence}.
\end{proof}
\smallskip

Finally, we establish a multiplier estimate for $|\cdot|_D$.
\begin{lemma}\label{lem:multiplier}
Let $0<s<1$ and $\gamma>-3$.
For any function $m(v)\in W^{1,\infty}(\mathbb R^3)$, we have
\begin{align}\label{mf-estimate}
|mf|_D\lesssim\Big\{|m|_{L^\infty}+|m|_{L^\infty}^{1-s}|\nabla_v m|_{L^\infty}^{s}\Big\}|f|_D.
\end{align}
The implicit constant is independent of $m$ and $f$.
\end{lemma}

\begin{proof}
The assertion is immediate when $m$ is constant.
Using
\begin{align*}
m(v')f(v')-m(v)f(v)=m(v')\big(f(v')-f(v)\big)+\big(m(v')-m(v)\big)f(v)
\end{align*}
and $|m(v_*)f(v_*)|\leq|m|_{L^\infty}|f(v_*)|$ in the definition of $|\cdot|_D$, we obtain
\begin{align*}
|mf|_D^2&\leq 2|m|_{L^\infty}^2|f|_D^2+2\iiint_{\mathbb R^6\times\mathbb S^2}B(v-v_*,\sigma)\mu(v_*) |m(v')-m(v)|^2|f(v)|^2\d\sigma\d v_*\d v.
\end{align*}
Since $|v'-v|=|v-v_*|\sin(\theta/2)$, we have
\begin{align*}
|m(v')-m(v)|\leq\min\Big\{2|m|_{L^\infty},|\nabla_v m|_{L^\infty}|v-v_*|\sin(\theta/2)\Big\}.
\end{align*}
The angular kernel assumption therefore gives
\begin{align*}
\int_{\mathbb S^2}b(\cos\theta)|m(v')-m(v)|^2\d\sigma
\lesssim\;& \int_0^{\frac{\pi}{2}}\theta^{-1-2s}\min\Big\{|m|_{L^\infty}^2,|\nabla_v m|_{L^\infty}^2|v-v_*|^2\theta^2\Big\}\d\theta\\
\lesssim\;& |m|_{L^\infty}^{2-2s}|\nabla_v m|_{L^\infty}^{2s}|v-v_*|^{2s}.
\end{align*}
Consequently, by \eqref{norm inequality}, we have
\begin{align*}
&\iiint_{\mathbb R^6\times\mathbb S^2}B(v-v_*,\sigma)\mu(v_*)|m(v')-m(v)|^2|f(v)|^2\d\sigma\d v_*\d v\\
\lesssim\;& |m|_{L^\infty}^{2-2s}|\nabla_v m|_{L^\infty}^{2s}\int_{\mathbb R^3}|f(v)|^2
\left(\int_{\mathbb R^3}\mu(v_*)|v-v_*|^{\gamma+2s}\d v_*\right)\d v\\
\lesssim\;&|m|_{L^\infty}^{2-2s}|\nabla_v m|_{L^\infty}^{2s}|f|_{L^2_{\gamma/2+s}}^2\\
\lesssim\;&|m|_{L^\infty}^{2-2s}|\nabla_v m|_{L^\infty}^{2s}|f|_D^2,
\end{align*}
where we have used the fact that
\begin{align*}
\int_{\mathbb R^3}\mu(v_*)|v-v_*|^{\gamma+2s}\d v_*
\lesssim \langle v\rangle^{\gamma+2s} \quad \text{for}\;\; \gamma+2s>-3.
\end{align*}
Thus, one has
\begin{align*}
|mf|_D^2\lesssim\Big\{|m|_{L^\infty}^2+|m|_{L^\infty}^{2-2s}|\nabla_v m|_{L^\infty}^{2s}\Big\}|f|_D^2,
\end{align*}
which gives \eqref{mf-estimate}.
\end{proof}
\medskip

\section{Some Basic Inequalities}\label{appendix:inequalities}

In what follows, we collect some basic inequalities to be used throughout this paper.
Firstly, the following Sobolev interpolation inequalities have been shown in \cite{GW2012CPDE}.

\begin{lemma}\label{sobolev interpolation}
Let $2 \leq p<\infty$ and $k, \ell, m \in \mathbb{R}$. Then for any $f\in \dot H^\ell(\mathbb R^3)\cap\dot H^m(\mathbb R^3)$, we have
\begin{align}\label{sobolev interpolation1}
\left\| f\right\|_{\dot{W}^{k,p}} \lesssim\left\|f\right\|_{\dot{H}^{\ell}}^\theta\left\|f\right\|_{\dot{H}^{m}}^{1-\theta},
\end{align}
where $0 \leq \theta \leq 1$ satisfies
\begin{align*}
\frac{1}{p}-\frac{k}{3}=\left(\frac{1}{2}-\frac{\ell}{3}\right) \theta+\left(\frac{1}{2}-\frac{m}{3}\right)(1-\theta).
\end{align*}
For the case $p=+\infty$, we have
\begin{align}\label{sobolev interpolation2}
\left\| f\right\|_{\dot{W}^{k,\infty}} \lesssim\left\| f\right\|_{\dot{H}^{\ell}}^\theta\left\| f\right\|_{\dot{H}^{m}}^{1-\theta},
\end{align}
where $\ell < k+3/2 < m$, and $0 \leq \theta \leq 1$ satisfies
\begin{align*}
-\frac{k}{3}=\left(\frac{1}{2}-\frac{\ell}{3}\right) \theta+\left(\frac{1}{2}-\frac{m}{3}\right)(1-\theta).
\end{align*}
\end{lemma}
\smallskip

As a direct consequence of Lemma \ref{sobolev interpolation}, we have the following standard Sobolev embeddings.

\begin{lemma}\label{sobolev embedding 1}
For any $f \in H^2(\mathbb{R}^3)$, we have
\begin{align}\label{Sobolev-ineq}
\left\|f\right\|_{L^3}\lesssim \left\|f\right\|_{H^1},
\quad \left\|f\right\|_{L^6}\lesssim \left\| \nabla f \right\|_{L^2},
\quad \left\|f\right\|_{L^\infty}\lesssim \left\| \nabla f\right\|_{H^1}.
\end{align}
\end{lemma}
\smallskip

We next recall the following Kato--Ponce product estimates, whose proofs can be found in \cite{Kato-Ponce-CPAM-1988}.

\begin{lemma}\label{kato-ponce}
Let $r>0$ and $1 < p < +\infty$. Assume that $q_j, r_j \in (1, +\infty]$ for $j=1, 2$ and that
\begin{equation*}
\frac{1}{p} = \frac{1}{q_1} + \frac{1}{r_1} = \frac{1}{q_2} + \frac{1}{r_2}.
\end{equation*}
Then, for $f\in W^{r, q_1} \cap L^{r_2} $ and $g\in W^{r, q_2}\cap L^{r_1}$, we have
\begin{align*}
\| fg\|_{\dot{W}^{r, p}} &\leq C (\| f\|_{\dot{W}^{r, q_1}}\|g\|_{L^{r_1}} + \|f\|_{L^{r_2}}\| g\|_{\dot{W}^{r, q_2}}),\\
\| fg\|_{W^{r, p}} &\leq C (\| f\|_{W^{r, q_1}}\|g\|_{L^{r_1}} + \|f\|_{L^{r_2}}\| g\|_{W^{r, q_2}}).
\end{align*}
\end{lemma}
\smallskip

If $r \in (0,3)$, then $\Lambda^{-r}$ is the Riesz potential operator.
The Hardy--Littlewood--Sobolev theorem implies the following inequality for the Riesz potential operator $\Lambda^{-r}$.

\begin{lemma}\label{negative embedding theorem}
Let $0<r<3$ and $1< p < q < \infty$ satisfy $1/q + r/3=1/p$. Then
\begin{align}\label{negative embed 1}
\left\|\Lambda^{-r} f\right\|_{L^q} \lesssim\|f\|_{L^p}.
\end{align}
\end{lemma}
\smallskip


The following consequence of Lemma \ref{sobolev interpolation} is used frequently in the analysis.

\begin{lemma}\label{negative embedding theorem2}
Let $0<r<3/2$. Then we have
\begin{align}\label{negative embed 2}
\|f\|_{L^{\frac{3}{r}}} \lesssim \left\|\Lambda^{\frac{3}{2}-r} f\right\|.
\end{align}
\end{lemma}
\smallskip

We use the following form of Minkowski's inequality to compare mixed norms with different orders of integration (see \cite{GW2012CPDE}).

\begin{lemma}\label{minkowski theorem}
For $1 \leq p \leq q \leq \infty$, there holds
\begin{align}\label{minkowski}
\|f\|_{L_x^q L_v^p} \leq\|f\|_{L_v^p L_x^q}.
\end{align}
\end{lemma}
\smallskip

The following equivalence of weighted Sobolev norms is proved in \cite{AMSY2023,HMUY2008}.

\begin{lemma}\label{equivalent norm}
Let $1 \leq p \leq \infty$ and $\ell, \theta \in \mathbb{R}$. Then there exists a constant $C>1$ independent of $f$ such that
\begin{align}\label{equivalent norm1}
\frac{1}{C}\left|\langle v\rangle^{\ell}\left\langle D_v\right\rangle^\theta f\right|_{L^p} \leq\left|\left\langle D_v\right\rangle^\theta\langle v\rangle^{\ell} f\right|_{L^p} \leq C\left|\langle v\rangle^{\ell}\left\langle D_v\right\rangle^\theta f\right|_{L^p},
\end{align}
that is, these two norms are equivalent. Here, $\left\langle D_v\right\rangle$ denotes the Fourier multiplier with symbol $\langle\xi\rangle$.
\end{lemma}
\smallskip

Next, we present a weighted interpolation inequality that follows from \cite[Lemma 4.15]{CHJ2024SIAM}.

\begin{lemma}\label{lem:weighted-interpolation}
Let $r_i,\ell_i \in\mathbb R$, $i=1,2,3$, and let $0<\theta<1$.
Suppose that
\begin{align*}
r_1=(1-\theta)r_2+\theta r_3,
\qquad
\ell_1=(1-\theta)\ell_2+\theta\ell_3.
\end{align*}
Then, for any function $f(v)$, we have
\begin{equation}\label{interpolation-general}
\left|\langle v\rangle^{\ell_1} f\right|_{H_v^{r_1}}
\lesssim\left|\langle v\rangle^{\ell_2}f\right|_{H_v^{r_2}}^{1-\theta}
\left|\langle v\rangle^{\ell_3}f\right|_{H_v^{r_3}}^\theta.
\end{equation}
\end{lemma}


The next lemma shows that the dissipation norm controls weighted fractional Sobolev norms.

\begin{lemma}\label{lem:D-fractional-control}
Let $0<s<1$ and $\gamma>-3$.
For any $0\leq\sigma\leq s$, there holds
\begin{equation}\label{D-fractional-control}
\left|\langle v\rangle^{\frac{\gamma+2s}{2}-\sigma}g\right|_{H_v^\sigma}
\lesssim |g|_D.
\end{equation}
\end{lemma}

\begin{proof}
For $0<\sigma<s$, applying Lemma \ref{lem:weighted-interpolation} with
\begin{align*}
r_1=\sigma,\quad r_2=0,\quad r_3=s,\quad
\ell_1=\frac{\gamma+2s}{2}-\sigma,\quad \ell_2=\frac{\gamma+2s}{2},\quad
\ell_3=\frac{\gamma}{2},\quad \theta=\frac{\sigma}{s},
\end{align*}
and using \eqref{norm inequality}, we obtain
\begin{align*}
\left|\langle v\rangle^{\frac{\gamma+2s}{2}-\sigma}g\right|_{H_v^\sigma}
\lesssim
\left|\langle v\rangle^{\frac{\gamma+2s}{2}}g\right|_{L_v^2}^{1-\frac{\sigma}{s}}
\left|\langle v\rangle^{\frac{\gamma}{2}}g\right|_{H_v^s}^{\frac{\sigma}{s}}
\lesssim |g|_D.
\end{align*}
The endpoint cases $\sigma=0$ and $\sigma=s$ follow directly from \eqref{norm inequality}.
This completes the proof.
\end{proof}

\bigskip

\noindent\textbf{Acknowledgements:}
SQL was supported by the National Natural Science Foundation of China under grant numbers 12325107 and 12471217.
FJZ was supported by the National Natural Science Foundation of China under grant number 12271179.
\medskip

\noindent\textbf{Data Availability Statement:}
Data sharing is not applicable to this article as no datasets were generated or analysed during the current study.
\medskip

\noindent\textbf{Conflict of Interest:}
The authors declare that they have no conflict of interest.
\medskip

\noindent\textbf{AI Usage Disclosure:}
The authors used AI-assisted tools for writing assistance and to help check the correctness of the mathematical arguments.
The authors have carefully reviewed all AI-assisted revisions and take full responsibility for the content, arguments, and conclusions presented in the manuscript.


\medskip

\begin{thebibliography}{99}

\bibitem{AMUXY2011AA} Alexandre, R.; Morimoto, Y.; Ukai, S.; Xu, C.-J.; Yang, T. {\it The Boltzmann equation without angular cutoff in the whole space: II, Global existence for hard potential.} Anal. Appl. (Singap.) {\bf 9} (2011), no. 2, 113--134.

\bibitem{AMUXY2012JFA} Alexandre, R.; Morimoto, Y.; Ukai, S.; Xu, C.-J.; Yang, T. {\it The Boltzmann equation without angular cutoff in the whole space: I, Global existence for soft potential.} J. Funct. Anal. {\bf 262} (2012), no. 3, 915--1010.

\bibitem{AMSY2023} Alonso, R.; Morimoto, Y.; Sun, W.; Yang, T. {\it De Giorgi argument for weighted $L^2 \cap L^\infty$ solutions to the non-cutoff Boltzmann equation.} J. Stat. Phys. {\bf 190} (2023), no. 2, Paper No. 38, 98 pp.

\bibitem{AS2019} Ars\'{e}nio, D.; Saint-Raymond, L. {\it From the Vlasov--Maxwell--Boltzmann system to incompressible viscous electro-magneto-hydrodynamics. Vol. 1.} EMS Monographs in Mathematics. European Mathematical Society, Z\"urich, 2019. xii+406 pp.


\bibitem{BGL91} Bardos, C.; Golse, F.; Levermore, D. {\it Fluid dynamic limits of kinetic equations. I. Formal derivations.} {J. Stat. Phys.} {\bf 63} (1991), no. 1--2, 323--344.

\bibitem{BGL93} Bardos, C.; Golse, F.; Levermore, D. {\it Fluid dynamic limits of kinetic equations. II. Convergence proofs for the Boltzmann equation.} {Comm. Pure Appl. Math.} {\bf 46} (1993),  no. 5, 667--753.


\bibitem{BU91} Bardos, C.; Ukai, S. {\it The classical incompressible Navier--Stokes limit of the Boltzmann equation.}  Math. Models Methods Appl. Sci. {\bf 1} (1991), no. 2, 235--257.

\bibitem{Ca} Caflisch, R. {\it The fluid dynamic limit of the nonlinear Boltzmann equation.} {Comm. Pure Appl. Math.} {\bf 33} (1980), no. 5,  651--666.

\bibitem{CIP} Cercignani, C.; Illner, R.; Pulvirenti, M. {\it The Mathematical Theory of Dilute Gases.} Applied Mathematical Sciences, {\bf 106}.  Springer-Verlag, New York, 1994. viii+347 pp.

\bibitem{CYZ2026JFA} Chen, K.; Yang, A.; Zhong, M. {\it Diffusion limit and the optimal convergence rate of the classical solution to the one-species Vlasov--Maxwell--Boltzmann system.} J. Funct. Anal. {\bf 290} (2026), no. 11, Paper No. 111413.

\bibitem{CHJ2024SIAM} Cao, C.; He, L.-B.; Ji, J. {\it Propagation of moments and sharp convergence rate for inhomogeneous noncutoff Boltzmann equation with soft potentials.} SIAM J. Math. Anal. {\bf 56} (2024), no. 1, 1321--1426.

\bibitem{DEL}De Masi, A.; Esposito, R.; Lebowitz, J. {\it Incompressible Navier--Stokes and Euler limits of the Boltzmann equation.} {Comm. Pure Appl. Math.} {\bf 42} (1989), no. 8, 1189--1214.

\bibitem{Duan2011} Duan, R. {\it Dissipative property of the Vlasov-Maxwell-Boltzmann system with a uniform ionic background.} SIAM J. Math. Anal. {\bf 43} (2011), no. 6, 2732--2757.

\bibitem{DLYZ2017} Duan, R.; Lei, Y.; Yang, T.; Zhao, H. {\it The Vlasov--Maxwell--Boltzmann system near Maxwellians in the whole space with very soft potentials.} Comm. Math. Phys. {\bf 351} (2017), no. 1, 95--153.

\bibitem{DL2013} Duan, R.; Liu, S. {\it The Vlasov--Poisson--Boltzmann system without angular cutoff.} Comm. Math. Phys. {\bf 324} (2013), no. 1, 1--45.

\bibitem{DLYZ2013} Duan, R.; Liu, S.; Yang, T.; Zhao, H. {\it Stability of the nonrelativistic Vlasov--Maxwell--Boltzmann system for angular non-cutoff potentials.} Kinet. Relat. Models {\bf 6} (2013), no. 1, 159--204.

\bibitem{DS2011} Duan, R.;  Strain, R. {\it Optimal large-time behavior of the Vlasov--Maxwell--Boltzmann system in the whole space.} Comm. Pure Appl. Math. {\bf 64} (2011), no. 11, 1497--1546.

\bibitem{DYY2023M3AS} Duan, R.;  Yang, D.; Yu, H. {\it Compressible Euler--Maxwell limit for global smooth solutions to the Vlasov--Maxwell--Boltzmann system.} Math. Models Methods Appl. Sci. {\bf 33} (2023), no. 10, 2157--2221.


\bibitem{FLLZ2018} Fan, Y.; Lei, Y.; Liu, S.; Zhao, H. {\it The non-cutoff Vlasov--Maxwell--Boltzmann system with weak angular singularity.} Sci. China Math. {\bf 61} (2018), no. 1, 111--136.

\bibitem{GS} Golse, F.; Saint-Raymond, L.  {\it The Navier--Stokes limit of the Boltzmann equation for bounded collision kernels.} Invent. Math. {\bf 155} (2004), no. 1, 81--161.



\bibitem{GS2011} Gressman, P.; Strain, R. {\it Global classical solutions of the Boltzmann equation without angular cut-off.} J. Amer. Math. Soc. {\bf 24} (2011), no. 3, 771--847.

\bibitem{Guo2003} Guo, Y. {\it The Vlasov--Maxwell--Boltzmann system near Maxwellians.} Invent. Math. {\bf153} (2003), no. 3,  593--630.

\bibitem{guo2006} Guo, Y. {\it Boltzmann diffusive limit beyond the Navier--Stokes approximation.} Comm. Pure Appl. Math. {\bf 59} (2006), no. 5,  626--687.

\bibitem{guo2010cmp} Guo, Y.; Jang, J.  {\it Global Hilbert expansion for the Vlasov--Poisson--Boltzmann system.} Comm. Math. Phys. {\bf 299} (2010), no. 2, 469--501.

\bibitem{GJJ} Guo, Y.; Jang, J.; Jiang, N. {\it Acoustic limit for the Boltzmann equation in optimal scaling.} Comm. Pure Appl. Math. {\bf 63} (2010), no. 3, 337--361.

\bibitem{GW2012CPDE}Guo, Y.; Wang, Y. {\it Decay of dissipative equations and negative Sobolev spaces.} {Comm. Partial Differential Equations} {\bf 37} (2012), 2165--2208.



\bibitem{HMUY2008} Huo, Z.; Morimoto, Y.; Ukai, S.; Yang, T. {\it Regularity of solutions for spatially homogeneous Boltzmann equation without angular cutoff.}  Kinet. Relat. Models {\bf 1} (2008), no. 3, 453--489.

\bibitem{J2009} Jang, J. {\it  Vlasov--Maxwell--Boltzmann diffusive limit.} Arch. Ration. Mech. Anal. {\bf 194} (2009), no. 2, 531--584.

\bibitem{JL2023ARXIV} Jiang, N.; Lei, Y. {\it The incompressible Navier--Stokes--Fourier--Maxwell system limits of the Vlasov--Maxwell--Boltzmann system for soft potentials: the noncutoff cases and cutoff cases.} SIAM J. Math. Anal. {\bf 58} (2026), no. 2, 1100--1170.

\bibitem{JL2019} Jiang, N.; Luo, Y.-L. {\it From Vlasov--Maxwell--Boltzmann system to two-fluid incompressible Navier--Stokes--Fourier--Maxwell system with Ohm's law: convergence for classical solutions.} Ann. PDE  {\bf 8} (2022), no. 1, Paper No. 4, 126 pp.

\bibitem{JLZ2023ARMA} Jiang, N.; Luo, Y.-L.; Zhang, T.-F. {\it   Hydrodynamic limit of the incompressible Navier--Stokes--Fourier--Maxwell system with Ohm's law from the Vlasov--Maxwell--Boltzmann system: Hilbert expansion approach.} Arch. Ration. Mech. Anal. {\bf 247} (2023), no. 3, Paper No. 55, 85 pp.

\bibitem{JXZ2018} Jiang, N.; Xu, C.-J.; Zhao, H. {\it Incompressible Navier--Stokes--Fourier limit from the Boltzmann equation: classical solutions.} Indiana Univ. Math. J. {\bf 67} (2018), no. 5, 1817--1855.

\bibitem{Kato-Ponce-CPAM-1988} Kato, T.; Ponce, G. {\it Commutator estimates and the {E}uler and {N}avier--{S}tokes equations.} Comm. Pure Appl. Math. {\bf 41} (1988), no. 7, 891--907.


\bibitem{LLXZ2023ARXIV} Lei, Y.; Liu, S.; Xiao, Q.; Zhao, H. {\it Global Hilbert expansion for some nonrelativistic kinetic equations.} J. Lond. Math. Soc. (2) {\bf 110} (2024), no. 5, Paper No. e70016, 82 pp.

\bibitem{LZ2016} Lei, Y.; Zhao, H. {\it The Vlasov--Maxwell--Boltzmann system with a uniform ionic background near Maxwellians.} J. Differential Equations {\bf 260} (2016), no. 3, 2830--2897.

\bibitem{LYZ2016-1} Li, H.-L.; Yang, T.; Zhong, M. {\it Spectrum structure and behaviors of the Vlasov--Maxwell--Boltzmann systems.} SIAM J. Math. Anal. {\bf 48} (2016), no. 1, 595--669.

\bibitem{Lions1994-1}Lions, P.-L. {\it Compactness in Boltzmann's equation via Fourier integral operators and applications. I, II.} J. Math. Kyoto Univ. {\bf 34} (1994), no. 2, 391--427, 429--461.

\bibitem{Lions1994-2}Lions, P.-L. {\it Compactness in Boltzmann's equation via Fourier integral operators and applications. III.} J. Math. Kyoto Univ. {\bf 34} (1994), no. 3, 539--584.

\bibitem{LM} Lions, P.-L.; Masmoudi, N.  {\it From the Boltzmann equations to the equations of incompressible fluid mechanics. I, II.} Arch. Ration. Mech. Anal. {\bf158} (2001), no. 3, 173--193, 195--211.

\bibitem{Saint} Saint-Raymond, L. {\it Hydrodynamic limits of the Boltzmann equation.} Lecture Notes in Mathematics, {\bf 1971}. Springer-Verlag, Berlin, 2009. xii+188 pp.

\bibitem{Strain2006CMP} Strain, R.  {\it The Vlasov--Maxwell--Boltzmann system in the whole space.} {Comm. Math. Phys.} {\bf 268} (2006), no. 2, 543--567.


\bibitem{V} Villani, C. {\it A review of mathematical topics in collisional kinetic theory.} Handbook of Mathematical Fluid Dynamics, Vol. I, 71--305, North-Holland, Amsterdam, 2002.

\bibitem{WZGX2025SIMA} Wu, W.; Zhou, F.; Gong, W.; Xu, Y. {\it Diffusive limit of the one-species Vlasov--Maxwell--Boltzmann system for cutoff hard potentials.} SIAM J. Math. Anal. {\bf 57} (2025), no. 2, 1526--1585.

\bibitem{XZGW2026SIMA} Xu, Y.; Zhou, F.; Gong, W.; Wu, W. {\it Diffusive limit of the noncutoff Vlasov--Maxwell--Boltzmann system.} SIAM J. Math. Anal. {\bf 58} (2026), no. 2, 1187--1231.

\bibitem{YZ24} Yang, T.; Zhong, M. {\it Diffusion limit with optimal convergence rate of classical solutions to the Vlasov--Maxwell--Boltzmann system.} Adv. Math. {\bf 489} (2026), Paper No. 110800, 95 pp.


\end{thebibliography}
\end{document}